\documentclass[10pt]{amsart}

\usepackage{color,graphicx,amssymb}
\usepackage{amsthm}
\usepackage{dsfont}
\usepackage{tikz-cd}

\theoremstyle{definition}

\newtheorem{theorem}{Theorem}[section]
\newtheorem{definition}[theorem]{Definition}
\newtheorem{lemma}[theorem]{Lemma}

\newtheorem{proposition}[theorem]{Proposition}

\newtheorem{remark}[theorem]{Remark}

\newtheorem{acknowledgments}[theorem]{Acknowledgments}

\DeclareMathOperator{\Con}{Con}
\DeclareMathOperator{\Aut}{Aut}

\DeclareMathOperator{\Hom}{Hom}

\DeclareMathOperator{\id}{id}
\DeclareMathOperator{\im}{im}
\DeclareMathOperator{\ar}{ar}
\DeclareMathOperator{\Id}{Id}

\begin{document}

\title[Nonabelian extensions of multiplicative Lie algebras]{Correspondence, Wells and Hochschild-Serre sequences for nonabelian extensions of multiplicative Lie algebras}
\author{Dev Karan Singh, Alexander Wires, Shiv Datt Kumar}
\address{School of Mathematical Sciences, University of Electronic Science and Technology of China, Chengdu, 611731, Sichuan, PRC}
\address{Department of Mathematics, Motilal Nehru National Institute of Technology Allahabad, 211004, Prayagraj, Uttar
Pradesh, India}
\email{dev.2020rma01@mnnit.ac.in (Dev Karan Singh), awires81@uestc.edu.cn (Alexander Wires), sdt@mnnit.ac.in (Shiv Datt Kumar)}
\date{April 11, 2024}

\begin{abstract}
For nonabelian $2^{\mathrm{nd}}$-cohomology of multiplicative Lie algebras, we properly generalize from the group case three classic results. We prove a Correspondence theorem which compares $2^{\mathrm{nd}}$-cohomology associated to a realized abstract kernel to the abelian $2^{\mathrm{nd}}$-cohomology group over the algebraic center. For arbitrary extensions, we prove a Wells's Theorem characterizing ideal-preserving automorphisms and establish the 1-dimensional Lyndon-Hochschild-Serre exact sequence. Several previously established results are recovered when restricted to extensions with group-abelian or Lie-trivial ideals.
\end{abstract}


\maketitle


\section{Introduction}\label{introduction}
\vspace{0.3cm}

The notion of a multiplicative Lie algebra was introduced by G. Ellis \cite{GJ} in order formulate a non-commutative version of the Magnus-Witt Theorem \cite{magwitt} which provides an isomorphism between the free Lie algebra $L(F_{\mathrm{ab}})$ generated over the abelianization of a free group $F$ and the Lie algebra formed from the lower central series of $F$. As pointed out by G. Ellis, a non-commutative version of the isomorphism can also be interpreted as determining a basis for the universal commutator identities for any fixed weight $n$. This is very closely related to the results of C. Miller \cite{cmiller} which can be seen as providing an equational basis for the universal commutator identities of weight 2 which is also well-known to be intimately connected to the Schur multiplier of finite groups. This provides a fertile connection between Lie algebras, nilpotency and central extensions in groups, and universal algebra.

There have been several works exploring this connection which for space must by necessity be an incomplete list. In \cite{FPW}, Point and Wantiez develop the theory of nilpotency in multiplicative Lie algebras. In \cite{AGNM}, Bak et al. consider a homology theory for multiplicative Lie algebras and define an appropriate Schur multiplier which in Lal and Upadhyay \cite{RLS} we find its characterization in terms of the second-cohomology of central extensions. This connection with the classical theory of central extensions in groups continues in Donadze et al. \cite{GNMA} which explore structural results deduced from the non-abelian tensor product and non-abelian exterior square for multiplicative Lie algebras developed in Donadze, Inassaridze and Ladra \cite{GDM}. Lie perfect multiplicative Lie algebra and local nilpotency are considered in Singh et al. \cite{DMS}. In Pandey and Upadhyay \cite{MSS}, we observe an intriguing characterization of what Lie products can be supported by a finite group.

The present manuscript is concerned with the general problems related to the analysis of extensions: given multiplicative Lie algebras $I$ and $Q$, classify all multiplicative Lie algebras $A$ such that $I$ is an ideal and the factor algebra $A/I$ is isomorphic to $Q$. The analogous problem for groups is addressed through the theory of group extensions and the cohomology of groups. In Schreier \cite{Schreier} and Baer \cite{Baer} we read the seminal results in the theory of group extensions; in particular, in \cite{Baer} we see that every extensions of $Q$ by $H$ yields a homomorphism from \( Q \) to \( \operatorname{Out}(H) \) (called an abstract kernel). If \( H \) is not abelian, then it is well known that not every abstract kernel gives rise to an extension that realizes the same abstract kernel; for this case, a whole theory of obstructions to group extensions is developed in Eilenberg and MacLane \cite{Elienberg}.

There have been several manuscripts concerning the reconstruction of extensions of multiplicative Lie algebras in terms of the machinery associated to $2^{\mathrm{nd}}$-cohomology. Previous work has often developed $2^{\mathrm{nd}}$-cohomology for extensions where the kernel has a restricted type such as abelian groups and contained in center and in null center \cite{MS1} or abelian groups with trivial brackets and contained in algebraic center \cite{RLS}. For these restricted types of kernels, analogues of Wells's theorem \cite{wells} has been shown in \cite{DK_MS_Sd}. We should make special note that the characterization of extensions for arbitrary kernels in terms of factor systems was previously elucidated in \cite{notedimath}.

With the future goal of formulating higher-cohomologies and obstructions in the manner of Eilenberg and MacLane, we refine the discussion around the category of arbitrary extensions and the $2^{\mathrm{nd}}$-cohomology category of compatible 2-cocycles (or factor systems). We prove a Correspondence Theorem (Theorem~\ref{thm:corresp1}) which takes the category of equivalence classes of extensions for realized action terms and relates it to the abelian $2^{\mathrm{nd}}$-cohomology group associated to the algebraic center. We establish a Wells's-type theorem (Theorem~\ref{thm:wells}) for general extensions. We also prove an analogue for multiplicative Lie algebras (Theorem~\ref{thm:hoschserre}) of the low-dimensional Hochschild-Serre exact sequence for $2^{\mathrm{nd}}$-cohomology and a general extension (see Hochschild and Serre \cite{hochserre}). The Hochschild--Serre spectral sequence is a fundamental tool in algebraic topology and homological algebra, particularly in the study of group extensions and their associated cohomology. This spectral sequence provides a bridge between the cohomology of the kernel and the cohomology of the quotient group in an extension. It is an essential result when analyzing the cohomology of a group from the knowledge of its subgroups.

The research in this manuscript is also related to the larger project of affine and general extensions of universal algebras elucidated in Wires \cite{wiresI}. Since multiplicative Lie algebras form a congruence-permutable variety, the notion of an abelian ideal (see Definition~\ref{ideal_def}(8)) agrees with the term-condition commutator (see Proposition~\ref{prop:commutator}) and so extensions with abelian ideals corresponds to extensions with affine datum; in this way, the machinery and basic constructions of $2^{\mathrm{nd}}$-cohomology, central extensions and the Schur multiplier fall within the purview of the theory expounded in \cite{wiresI} and Wires \cite{wiresIV}. It is hoped that the results in the present manuscript give evidence for future development in arbitrary varieties of universal algebras in general and for the continued development of the theory of groups with multiple operators in Higgins \cite{higgins}, in particular.

\vspace{0.2cm}

\section{Preliminaries}\label{preliminaries}


\begin{definition}\label{multiplicative} \cite{GJ}
A \emph{multiplicative Lie algebra} is an algebra $\left\langle A, \cdot, 1, \{~,~\} \right\rangle$ where $\left\langle A, \cdot, 1 \right\rangle$ is a group with identity 1 and the binary operation $\{~,~\}: A \times A \to A$ satisfies the identities 
\begin{equation}\label{(2.1)}
			\{x,x\} = 1
\end{equation}
\begin{equation}\label{(2.2)}
			\{x , y\cdot z\} = \{x,y\}\cdot {\{x,z\}^y}
\end{equation}
\begin{equation}\label{(2.3)} 
			\{ x\cdot y, z\} = {\{y, z\}^x}\cdot \{x, z\}
\end{equation}
\begin{equation}\label{(2.4)}
			\{\{x, y\},{z^y}\} \cdot \{ \{ y, z\}, {x^z}\} \cdot \{ \{ z, x\},{y^x}\} = 1
\end{equation}
\begin{equation}\label{(2.5)}
			\{x, y\}^z = \{ {x^z},{y^z} \}
\end{equation}
for all $x,y,z \in A$. The binary operation $\{~,~\}$ is called the multiplicative Lie product (or Lie bracket) and $x^z = zxz^{-1}$. 
\end{definition}


\noindent \textbf{Note:} If $\left\langle A, \cdot, 1 \right\rangle$ is an abelian group, then $\left\langle  A, \cdot, 1, \{~,~\} \right\rangle$ is called a \emph{Lie ring}.


\begin{proposition}\label{Identity}
Multiplicative Lie algebras satisfy the following additional identities: 
\begin{equation}\label{(2.6)}
\{x, y\}\cdot \{y, x\} = 1
\end{equation}
\begin{equation}\label{(2.7)}
\{u, v\}^{\{x, y\}} = \{ u, v\}^{^{[x,y]}}
\end{equation}
\begin{equation}\label{(2.8)}
[\{ x, y\},z] = \{[x,y], z\}
\end{equation}
\begin{equation}\label{(2.9)}
\{ x^{-1}, y \} = (\{x, y\}^{-1}){^{x^{-1}}}
\end{equation}
\begin{equation}\label{(2.10)}
\{ 1, x\} = 1 = \{ x, 1\}
\end{equation}
\begin{equation}\label{(2.11)}
\{x, y\} = \{ y^{-1}, {x^y}\} 
\end{equation}
\begin{equation}\label{(2.12)}
\{x, y_{1}\cdots y_{n}\} = \{ x, y_{1}\} \cdot \{x,y_{2}\}^{y_{1}} \cdots \{x,y_{n} \}^{y_{1} \cdots y_{n-1}} 
\end{equation}
In addition, we have the following implication:
\begin{equation}\label{(2.13)}
\{ x, y\} = \{x, z\} \rightarrow \{ x , y^{-1} \cdot z \} = 1
\end{equation}
\end{proposition}
\begin{proof}
Eq~\eqref{(2.6)}-\eqref{(2.10)} were established in \cite{GJ} and Eq~\eqref{(2.11)} follows from Eq~\eqref{(2.9)}. Repeated application of Eq~\eqref{(2.2)} yields Eq~\eqref{(2.12)}. For Eq~\eqref{(2.13)}, assume $\{x, y\} = \{ x, z \}$ and observe
\begin{align*}
\{ x, y^{-1} \cdot z \} = {\{ x, y^{-1}\}} \cdot \{x, z\}{^{y^{-1}}} = {\{ x, y^{-1} \}} \cdot \{ x, y\}{^{y^{-1}}} &= \{ x, y^{-1} \} \cdot (\{ y, x \}^{-1})^{y^{-1}} \\
&= \{ x, y^{-1} \} \cdot \{ y^{-1}, x \} = 1 .
\end{align*} 
\end{proof}


The variety of multiplicative Lie algebras is denoted by $\mathcal M \mathcal L$. It follows from Eq~\eqref{(2.10)} that the group identity is idempotent for the bracket and so multiplicative Lie algebras are examples of groups with multiple operators in the vernacular of Higgins \cite{higgins}. As a variety of algebras which are expansions of groups, $\mathcal M \mathcal L$ is a Mal'cev variety. Given $A \in \mathcal M \mathcal L$, let $A^{\mathrm{gr}}$ denote the reduct to the group operations. Note the additional equation $\{x, y\} = 1$ determines the subvariety $\mathcal M\mathcal L^{\mathrm{gr}} = \{ A \in \mathcal M \mathcal L: A \vDash \{ x , y \} = 1 \}$ bi-interpretable with the variety of groups; informally, we recover groups as multiplicative Lie algebras with trivial brackets and in this way many of our results properly generalize the group case. The subvariety $\mathcal M\mathcal L^{\mathrm{Lie}} = \{ A \in \mathcal M \mathcal L: A \vDash [ x , y ] = 1 \}$ recovers $\mathds{Z}$-module Lie algebras.


\begin{definition}\label{ideal_def}
Let $\left\langle  A, \cdot, 1, \{~,~\} \right\rangle$ be a multiplicative Lie algebra. 
\begin{enumerate}
	
	\item[\textit{(1)}] A \emph{subalgebra} $H \leq A$ is a subgroup of $A$ which is closed under the bracket.
	
	\item[\textit{(2)}] A subalgebra $H$ of $A$ is said to be an \emph{ideal} of $G$ if it is a normal subgroup of $A$ and $\{h, x\} \in H$ for all $h \in H$ and $x \in A$.
	
	\item[\textit{(3)}] The \emph{group center} of $A$ is $Z(A^{\mathrm{gr}})$ the standard center of the group reduct.
	
	\item[\textit{(4)}] The \emph{Lie center} of $A$ is defined by 
\begin{align*}
LZ(A)=\{x \in A: \{x, y\} =1~\text{for all} ~ y \in A\} .
\end{align*}
			
	\item[\textit{(5)}] The \emph{null} ideal of $A$ is $\mathrm{Null} \, A = \{ a \in A : \{ a,x\} = \{ x,a\} = 0 \text{ for all } x \in A \}$. 

	\item[\textit{(6)}] The \emph{center} (or \emph{algebraic center}) of $A$ is defined by 
\begin{align*}
\zeta(A) &= Z(A^{\mathrm{gr}}) \wedge \mathrm{Null} \, A \\
&= \{ x \in A : x \in Z(A^{\mathrm{gr}}) \text{ and } \{ x,a\} = \{ a,x\} = 0 \text{ for all } a \in A \} .
\end{align*}

	\item[\textit{(7)}] For subalgebra $H \leq A$, the \emph{normalizer} of $H$ denoted by $N_{A}(H)$ is the largest subalgebra of $A$ which contains $H$ as an ideal.
	
	\item[\textit{(8)}] An ideal $I \triangleleft A$ is \emph{abelian} if the group operation is commutative and the bracket operation is trivial in the ideal.

\end{enumerate}
\end{definition}


For subset $I \subseteq A$, define $\alpha_{I} = \{(a,b) \in A : b^{-1} \cdot a \in I \}$. For a congruence $\alpha \in \Con A$, write $I_{\alpha} = 1/\alpha$ for the congruence class containing the neutral element $1 \in A$. For multiplicative Lie algebras, congruences are ideal determined.


\begin{lemma}\label{lem:congruence}
Let $A \in \mathcal M\mathcal L$. Then $\alpha \in \Con A$ if and only if $\alpha = \alpha_{I}$ for an ideal $I \triangleleft A$; furthermore, $\alpha = \alpha_{I_{\alpha}}$ and $I = I_{\alpha_{I}}$.
\end{lemma}
\begin{proof}
Let $\alpha \in \Con A$. Then $\alpha$ is congruence of the reduct $A^{\mathrm{gr}}$; thus, $\alpha = \alpha_{I}$ for a normal subgroup $I$ of $A^{\mathrm{gr}}$ where $I = 1/I$ the congruence class containing the neutral element $1 \in A$. For $x \in A$ and $a \in I$, we have $\{ a,x\} \mathrel{\alpha} \{ 1,x\}=1$ which shows $\{ a,x\} \in I$; similarly, $\{ x,a\} \in I$. We conclude $I$ is an ideal. 
		
Now assume $I \triangleleft A$. Then $\alpha_{I} \in \Con A^{\mathrm{gr}}$. The argument is done once we show $\alpha_{I}$ is compatible with the bracket operation. Let $(a,b), (c,d) \in \alpha_{I}$. Then there exists $x,y \in I$ such that $a = bx$ and $c = dy$. Then 
\begin{align*}
\{ a,c\} &= \{ bx,dy\} = \{ x,d\}^{b} \cdot \{ b,d\} \cdot \{ x,y\}^{db} \cdot \{ b,y\}^{d} = \{ b,d\} \cdot \{ x,d\}^{\{ d,b\}b} \cdot \{ x,y\}^{db} \cdot \{ b,y\}^{d} \in I
\end{align*}
since $I$ is an ideal. This shows $\{ a,c\} \mathrel{\alpha_{I}} \{ b,d\}$; altogether, $\alpha_{I}$ is a congruence of $I$. 
\end{proof}


\begin{lemma}\label{lem:4}
Let $A \in \mathcal M \mathcal L$. The group center, null ideal and center are ideals of $A$; in particular, if $I \triangleleft A$, then $\mathrm{Null} \, I \triangleleft A$ and $\zeta(I) \triangleleft A$. 
\end{lemma}	
\begin{proof}
We already know $Z(A^{\mathrm{gr}})$ is a normal subgroup of $A$. Take $a \in Z(A^{\mathrm{gr}})$ and $x \in A$. Then for all $y \in A$, $[ y , \{ a, x \} ] = \{ y , [a,x] \} = \{ y,1 \} = 1$ which implies $[a,x]$ commutes with every element of $A$ and so $[a,x] \in Z(A^{\mathrm{gr}})$; therefore, $Z(A^{\mathrm{gr}}) \triangleleft A$. 

Now assume $I \triangleleft A$. Take $a,b \in \mathrm{Null} \, I$. Then for all $x \in I$, we have $1=\{ bb^{-1},x\} = \{ b^{-1},x\}^{b}\{ b,x\} = \{ b^{-1},x\}^{b}$ which implies $\{ b^{-1},x\}=1$. Then $\{ ab^{-1},x\} = \{ b^{-1},x\}^{a}\{ a,x\} = 1$ which implies $ab^{-1} \in \mathrm{Null} \, A$. Take $a \in \mathrm{Null} \, I$, $x \in A$ and $y \in I$. Then $\{ y,a^{x}\} = \{ y^{x^{-1}},a\}^{x} = 1$ since $y^{x^{-1}} \in I$ because $I$ is a normal subgroup of $A$; similarly, $\{ a^{x},y\}=1$. This shows that $a^{x} \in \mathrm{Null} \, I$; thus, $\mathrm{Null} \, I$ is a normal subgroup of $A$. Note $y \in I$, $a \in \mathrm{Null} \, I$ implies $y^{a^{-1}} \in I$. Then $\{ y^{a^{-1}},x\} \in I$ since $I \triangleleft A$. Because $a^{x} \in \mathrm{Null} \, I$, we see that $\{ \{ y^{a^{-1}},x\},a^{x}\}^{-1} = 1$. Also, $a \in \mathrm{Null} \, I$, $y^{a^{-1}} \in I$ implies $\{ a,y^{a^{-1}}\} = 1$. Then $\left\{\{ a,y^{a^{-1}}\},x^{y^{a^{-1}}} \right\}^{-1}=1$. Then 
\begin{align*}
\{ \{ x,a\},y\} = \{ \{ x,a\},(y^{a^{-1}})^{a}\} = \left\{\{ y^{a^{-1}},x\},a^{x} \right\}^{-1} \cdot \left\{\{ a,y^{a^{-1}} \}, x^{y^{a^{-1}}} \right\}^{-1} = 1
\end{align*}
		We also see that $a^{y^{x^{-1}}} \in \mathrm{Null} \, I$ because $\mathrm{Null} \, I$ is a normal subgroup of $A$. Note $I \triangleleft A$ implies $y^{x^{-1}} \in I$ and $\{ x,y^{x^{-1}}\} \in I$. Then 
\begin{align*}
\{ \{ a,x\},y\} = \left\{\{ a,x\}, (y^{x^{-1}})^{x} \right\} = \left\{\{y^{x^{-1}},a \} , x^{a} \right\}^{-1} \cdot \left\{\{ x,y^{x^{-1}}\} , a^{y^{x^{-1}}} \right\}^{-1} = \{ 1,x^{a}\}^{-1} \cdot 1 = 1. 
\end{align*}
		Then $\{ y,\{ a,x\}\} = \{ \{ a,x\},y\}^{-1} = 1$ and $\{ y,\{ x,a\}\} = \{ \{ x,a\},y\}^{-1} = 1$; altogether, $\{ a,x\}, \{ x,a\} \in \mathrm{Null} \, I$. We have shown $\mathrm{Null} \, I \triangleleft A$.

We now consider the algebraic center. We first observe $\zeta(I)$ is a normal subgroup of $A$. Take $a \in \zeta(I)$, $x \in A$ and $b \in I$. Since $a \in \zeta(I) \subseteq \mathrm{Null} \, I$ and $\mathrm{Null} \, I \triangleleft A$ by Lemma~\ref{lem:4}, we have $a^{x} \in \mathrm{Null} \, I$. If conjugation is written by $i_{x}(a) = a^{x}$, then we have $a^{x} \cdot b = i_{x} \circ i_{x^{-1}}(a^{x} \cdot b) = i_{x} (a \cdot i_{x^{-1}}(b) ) = i_{x} ( i_{x^{-1}}(b) \cdot a ) = b \cdot a^{x}$; thus, $a^{x} \in \zeta(I)$. This also implies that $[a,x] = a \cdot x \cdot a^{-1} \cdot x^{-1} \in \zeta(I)$. By Eq~\eqref{(2.8)}, we have $[\{ a,x\}, b ] = \{[ a, x] , b\} = 1$; similarly, $[\{ x,a\}, b] = 1$. This implies $\{ x,a\},\{ a,x\} \in Z(I^{\mathrm{gr}})$. We already have $\{ x,a\},\{ a,x\} \in \mathrm{Null} \, I$ by Lemma~\ref{lem:4}; altogether, $\zeta(I)$ is an ideal of $A$. 
\end{proof}


Note that $\zeta(A)$ is a \emph{characteristic ideal}; that is, $\phi( \zeta(A)) = \zeta(A)$ for any $\phi \in \Aut \, A$.


\begin{lemma}\label{lem:normalizer}
Let $A \in \mathcal M\mathcal L$ and $H \leq A$. Then $N_{A}(H) = \{a \in A : aH = Ha \text{ and } [a,H] \subseteq H \}$. 
\end{lemma}
\begin{proof}
If we let $S = \{a \in A : aH = Ha \text{ and } [a,H] \subseteq H \}$, then $N_{A}(H) \subseteq S$ since $H \triangleleft N_{A}(H)$. Equality follows once we show $S$ is a subalgebra which contains $H$ as an ideal. We already know $\{a \in A : aH = Ha \}$ is a subgroup. If $a,b \in S$, then for $x \in H$ we see that $\{b,x\}^{-1} \in H$ since $b \in S$ and $H$ is a subalgebra. Then $(\{b,x\}^{-1})^{b^{-1}} \in H$ since $b^{-1}Hb = H$; similarly, $\{a,x\}^{b^{-1}} \in H$. Then we have $\{ ab^{-1},x \} = \{b^{-1},x\}\{a,x\}^{b^{-1}} = (\{b,x\}^{-1})^{b^{-1}}\{a,x\}^{b^{-1}} \in H$; thus, $S$ is a subgroup. 

We also see that for $a,b \in S$, $\{a,b\}H = \{a,b\}H\{a,b\}^{-1} \{a,b\} = \big[ \{a,b\},H \big] H \{a,b\} = \big\{ [ a,b ],H \big\} H \{a,b\} \subseteq H \{a,b\}$ since $S$ is a subgroup which implies the group commutator term $[a,b] \in S$. A similar calculation yields $H\{a,b\} \subseteq \{a,b\}H$; therefore, $S$ is closed under the bracket.  
\end{proof}


For ideals $I,J \triangleleft A$, let $\{ I,J\}$ be the set of all finite products of elements of the form $\{ a,b\}$ or $\{ b,a \}$ with $a\in I, b \in J$ and let $[ I, J ]$ denote the group commutator subgroup.


\begin{lemma}
For $I,J \triangleleft A \in \mathcal M\mathcal L$, $\{ I,J\} \triangleleft A$.
\end{lemma}
\begin{proof}
Application of Eq~\eqref{(2.5)} and Eq~\eqref{(2.6)} show $\{ I,J\}$ is a normal subgroup. Take $a \in I$, $b \in J$. Then for any $x \in A$ we have by Eq~\eqref{(2.4)}
\begin{align*}
\{ x, \{ a,b \} \} = \{ (x^{b^{-1}})^{b},\{ a,b \} \} = \big\{ \{ b,x^{b^{-1}} \} , a^{x^{b^{-1}}} \big\} \cdot \{ \{ x^{b^{-1}},a \}, b^{a} \} \in \{I,J\};
\end{align*}
similarly, $\{ x, \{ b,a \} \} \in \{ I,J \}$. Then for any $w = \{ a_{1},b_{1} \} \cdots \{ a_{n},b_{n} \} \in \{ I,J \}$, we can use Eq~\eqref{(2.12)} to write $\{ x,w \}$ as a finite product of elements from $\{ I,J \}$. This shows $\{ I,J \} \triangleleft A$.
\end{proof}


In the following, we describe the term-condition (TC) commutator from universal algebra in terms of the bracket and group commutators. Since multiplicative Lie algebras form a congruence-permutable variety, our reference for the term-condition commutator is Freese and McKenzie \cite{commod}. The following proposition justifies our terminology in Definition~\ref{ideal_def}(8) for an abeian ideal. The term-condition commutator of congruences $\alpha_{I}$ and $\alpha_{J}$ is denoted by $[\alpha_{I}, \alpha_{J}]$.


\begin{proposition}\label{prop:commutator}
Let $I,J$ be ideals of $A \in \mathcal M\mathcal L$. Then $I_{[\alpha_{I}, \alpha_{J}]}$ is the smallest ideal containing $\{I,J\}$ and $[I,J]$.
\end{proposition}
\begin{proof}  
Let $R$ be the smallest ideal containing $\{I,J\}$ and $[I,J]$. The matrices $\begin{bmatrix} 1 & 1 \\ 1 & [a,b] \end{bmatrix} = \begin{bmatrix} \{1,1\} & \{1,b\} \\ \{a,1\} & \{a,b\} \end{bmatrix} \in M(\alpha_{I},\alpha_{J})$ and $\begin{bmatrix} 1 & 1 \\ 1 & [a,b] \end{bmatrix} = \begin{bmatrix} [1,1] & [1,b] \\ [a,1] & [a,b] \end{bmatrix} \in M(\alpha_{I},\alpha_{J})$ show 
$\alpha_{R} \leq [\alpha_{I}, \alpha_{J}]$. In order to show, $[\alpha_{I}, \alpha_{J}] \leq \alpha_{R}$ we will verify the centralizer condition $A \vDash C(\alpha_{I},\alpha_{J}; \alpha_{R})$. Since multiplicative Lie algebras form a congruence-permutable variety, this is equivalent to the centralizer condition $A' \vDash C(\alpha_{I'},\alpha_{J'};0)$ in the quotient algebra $A' = A/R$ for ideals $I' = I/R$ and $J' = J/R$.

Fix a term $t(\bar{x},\bar{y})$, elements $\bar{a} = (a_{1},\ldots,a_{n}),\bar{b}=(b_{1},\ldots,b_{n})$ pair-wise $\alpha_{I'}$-related, elements $\bar{c} = (c_{1},\ldots,c_{n}),\bar{d}=(d_{1},\ldots,d_{n})$ pair-wise $\alpha_{J'}$-related, and assume 
\begin{align}\label{eqn:TCstart}
t(\bar{a},\bar{c}) = t(\bar{a},\bar{d}) .
\end{align} 
If we fix a choice $s: A'/I' \rightarrow A'$ of coset-representatives for $I'$ and $t: A'/J' \rightarrow A'$ coset-representatives for $J'$, then we can write each $a_{i} = a'_{i}s_{i}$, $b_{i} = b_{i}' s_{i}$ and $c_{i} = c_{i}'t_{i}, d_{i} = d_{i}'t_{i}$ for some $a_{i}',b_{i}' \in I'$ and $c_{i}',d_{i}' \in J'$. 
If we allow for the notation $\bar{x} \cdot \bar{z} = (x_{1} \cdot z_{1},\ldots,x_{n}\cdot z_{n})$, then we have for the term $t(\bar{x} \cdot \bar{z}, \bar{y} \cdot \bar{w})$ the evaluation 
\begin{align}\label{eqn:TCnext}
t(\bar{a}'\cdot \bar{s},\bar{c}'\cdot \bar{t}) = t(\bar{a},\bar{c}) = t(\bar{a},\bar{d}) = t(\bar{a}'\cdot \bar{s},\bar{d}'\cdot \bar{t})
\end{align}
We now observe that we can decompose our term into a group term applied to special bracket terms. Let us denote by $\beta(x_{1},\ldots,x_{n})$ a \emph{bracket term} which is a term built by repeated composition of the bracket operation where each variable appears as either $x_{i}$ or $x_{i}^{-1}$. A bracket term of weight 0 is just $\beta(x)=x$ or $\beta(x) = x^{-1}$ and the bracket terms $\beta(x,y)$ of weight 1 are $[x,y]$, $[x^{-1},y]$ and $[x,y^{-1}]$; thus, bracket terms of weight $n$ have $n$ nested bracket operations. By repeatedly using axioms Eq~\eqref{(2.2)}, Eq~\eqref{(2.3)} and Eq~\eqref{(2.5)}, we see that there is a group term $f$ and bracket terms $\beta_{1},\ldots,\beta_{k}$ such that  
\begin{align*}
t(\bar{x} \cdot \bar{z} , \bar{y} \cdot \bar{w}) = f(\beta_{1}(\bar{x}_{1},\bar{z}_{1},\bar{y}_{1},\bar{w}_{1}),\ldots,\beta_{k}(\bar{x}_{k},\bar{z}_{k},\bar{y}_{k},\bar{w}_{k}))
\end{align*}
for subsets of variables $\bar{x}_{i} \subseteq \bar{x},\bar{z}_{i} \subseteq \bar{z}, \bar{y}_{i} \subseteq \bar{y}, \bar{w}_{i} \subseteq \bar{w}$ in which some may be empty, but not all. Since $R$ contains the bracket ideal $\{I,J\}$, we see that either $\beta_{i}(\bar{a}'_{i},\bar{s}_{i},\bar{c}'_{i},\bar{t}_{i}) = 1$ if $\bar{a}'_{i} \neq \emptyset \neq \bar{c}'_{i}$, $\beta_{i}(\bar{a}'_{i},\bar{s}_{i},\bar{c}'_{i},\bar{t}_{i}) \in I'$ if $\bar{a}'_{i} \neq \emptyset, \bar{c}'_{i} = \emptyset$, $\beta_{i}(\bar{a}'_{i},\bar{s}_{i},\bar{c}'_{i},\bar{t}_{i}) \in J'$ if $\bar{a}'_{i} = \emptyset, \bar{c}'_{i} \neq \emptyset$, or is an evaluation on coset-representative if $\bar{a}'_{i} = \emptyset = \bar{c}'_{i}$. This implies the pair of elements $\beta_{i}(\bar{a}'_{i},\bar{s}_{i},\bar{c}'_{i},\bar{t}_{i})$ and $\beta_{i}(\bar{a}'_{i},\bar{s}_{i},\bar{d}'_{i},\bar{t}_{i})$ will be $\alpha_{J'}$-related and the pair $\beta_{i}(\bar{a}'_{i},\bar{s}_{i},\bar{c}'_{i},\bar{t}_{i})$ and $\beta_{i}(\bar{b}'_{i},\bar{s}_{i},\bar{c}'_{i},\bar{t}_{i})$ will be $\alpha_{I'}$-related.

Because the TC-commutator agrees with the usual commutator for groups (see \cite[Ch. 1]{commod}) and $[I',J'] = 0$ in $A'$, we can apply the group centralizer condition to Eq~\ref{eqn:TCnext} 
\begin{align*}
f(\beta_{1}(\bar{a}'_{1},\bar{s}_{1},\bar{c}'_{1},\bar{t}_{1}),\ldots,\beta_{k}(\bar{a}'_{k},\bar{s}_{k},\bar{c}'_{k},\bar{t}_{k})) &= t(\bar{a}'\cdot \bar{s},\bar{c}'\cdot \bar{t}) \\
&= t(\bar{a}'\cdot \bar{s},\bar{d}'\cdot \bar{t}) \\
&= f(\beta_{1}(\bar{a}'_{1},\bar{s}_{1},\bar{d}'_{1},\bar{t}_{1}),\ldots,\beta_{k}(\bar{a}'_{k},\bar{s}_{k},\bar{d}'_{k},\bar{t}_{k}))
\end{align*}
to conclude 
\begin{align*}
t(\bar{b},\bar{c}) = t(\bar{b}'\cdot \bar{s},\bar{c}'\cdot \bar{t}) &= f(\beta_{1}(\bar{b}'_{1},\bar{s}_{1},\bar{c}'_{1},\bar{t}_{1}),\ldots,\beta_{k}(\bar{b}'_{k},\bar{s}_{k},\bar{c}'_{k},\bar{t}_{k})) \\
&= f(\beta_{1}(\bar{b}'_{1},\bar{s}_{1},\bar{d}'_{1},\bar{t}_{1}),\ldots,\beta_{k}(\bar{b}'_{k},\bar{s}_{k},\bar{d}'_{k},\bar{t}_{k})) \\
&= t(\bar{b}'\cdot \bar{s},\bar{d}'\cdot \bar{t}) \\
&= t(\bar{b},\bar{d}) .
\end{align*}
This verifies that $A' \vDash C(\alpha_{I'},\alpha_{J'};0)$ and finishes the demonstration. 
\end{proof}


It follows that an ideal is abelian in the term-condition commutator if and only if the group operation is commutative and the bracket operation is trivial in the ideal. Our definition of the algebraic center in Definition~\ref{ideal_def}(6) agrees with the terminology of the TC-commutator which defines the center $\zeta \in \Con A$ as the largest congruence $\alpha \in \Con A$ such that $[\alpha,1]=0$. According to Lemma~\ref{lem:congruence} and Proposition~\ref{prop:commutator}, we see that $I_{\zeta}$ is the largest ideal contained in the group center and null ideal; that is, $I_{\zeta} = Z(A^{\mathrm{gr}}) \wedge \mathrm{Null} \, A = \zeta(A)$. Let $C(I,J)$ be the smallest ideal containing $[I,J]$ and $\{I,J\}$. By Proposition~\ref{prop:commutator}, $C(I,J)$ is the ideal corresponding to the TC-commutator of the ideals $I$ and $J$; in terms of the TC-commutator, $\alpha_{C(I,J)} = [\alpha_{I},\alpha_{J}]$.

\vspace{0.2cm}

\section{Extensions of multiplicative Lie algebras}\label{Extensions}


\begin{definition}\label{def:extension}
\emph{Datum} is a pair $(Q,I)$ of algebras in $\mathcal M \mathcal L$. An algebra $A \in \mathcal M \mathcal L$ is an \emph{extension} of $I$ by $Q$ if there is an ideal $K \triangleleft A$ such that $K \approx I$ and $A/K \approx Q$.
\end{definition}


Let $\mathrm{EXT} \, \mathcal M \mathcal L$ be the category whose objects are three term exact sequences 
\begin{align}\label{eqn:exactseq}
1 \longrightarrow I \stackrel{i}{\longrightarrow} A \stackrel{\pi}{\longrightarrow} Q \longrightarrow 1
\end{align}
of algebras in $\mathcal M \mathcal L$ and the morphisms are triples $(\alpha,\lambda,\beta)$ of homomorphisms which yield commuting squares 
between exact sequences.


\begin{lemma}
Two extensions are isomorphic in the category $\mathrm{EXT} \, \mathcal M \mathcal L$ if and only if there is a morphism $(\alpha,\lambda,\beta)$ between the two exacts sequences where each component is an isomorphism. 
\end{lemma}


We can define different equivalence relations on the objects by restricting the type of morphisms which appear in the commuting squares. By fixing a particular extension, the elements in the equivalence blocks for that fixed extension form algebraic structures of interest. By composing with isomorphisms, any extension of $I$ by $Q$ according to Definition~\ref{def:extension} yields an exact sequence in the form of Eq~\eqref{eqn:exactseq}. Our interest will be on reconstructing extension which are concretely built from the algebras $I$ and $Q$. An object in $\mathrm{EXT} \, \mathcal M \mathcal L$ realizes datum $(Q,I)$ if it is of the form in Eq~\eqref{eqn:exactseq}.


\begin{definition}\label{def:equivexactseq}
Two extensions $1 \longrightarrow I \stackrel{i}{\longrightarrow} A \stackrel{\pi}{\longrightarrow} Q \longrightarrow 1$ and $1 \longrightarrow I \stackrel{i}{\longrightarrow} M \stackrel{\pi'}{\longrightarrow} Q \longrightarrow 1$ are \emph{equivalent} if there is an isomorphism $(\id_{I},\lambda,\id_{Q})$ in $\mathrm{EXT} \, \mathcal M \mathcal L$ between them.
\end{definition}


The above definition defines an equivalence relation which refines isomorphism types on the objects of the category. Let $\mathrm{EXT}(Q,I)$ denote the set of all equivalence classes of extensions realizing datum $(Q,I)$ and denote by $[A]$ the equivalence class of an extension in the form of Eq~\eqref{eqn:exactseq}.

We examine how extensions realizing datum $(Q,I)$ can be deconstructed in terms of actions terms and 2-cocycles. Consider an extension $1 \longrightarrow I \stackrel{i}{\longrightarrow} A \stackrel{\pi}{\longrightarrow} Q \longrightarrow 1$ and fix a section $l: Q \rightarrow A$ for the surjection $\pi : A \rightarrow Q$ such that $l(1)=1$. We consider the representation of the operations in terms of the right-coset determined by $l$. The group operation takes the familiar form
\begin{align*}
a \cdot l(x) \cdot b \cdot l(y) &= a \cdot \sigma_{x}(b) \cdot T(x,y) \cdot l(xy) 
\end{align*}
where we have defined $\sigma_{x}(a) := l(x) \cdot a \cdot l(x)^{-1} =  a^{l(x)}$ and $T(x,y) := l(x) \cdot l(y) \cdot l(xy)^{-1}$. The bracket operation decomposes as
\begin{align*}
\{ a \cdot l(x), b \cdot l(y)\} &= \{ l(x),b\}^{a} \cdot \{ a,b\} \cdot \{ l(x),l(y)\}^{ba} \cdot \{ a,l(y)\}^{b} \\
&= \{ l(x),b\}^{a} \cdot \{ a,b\} \cdot T_{f}(x,y)^{ba} \cdot l(\{ x,y\})^{ba} \cdot \{ a,l(y)\}^{b} \\
&= \{ l(x),b\}^{a} \cdot \{ a,b\} \cdot T_{f}(x,y)^{ba} \cdot \big[ ba , l(\{ x,y\}) \big] \cdot l(\{ x,y\}) \cdot \{ a,l(y)\}^{b} \\
&= \{ l(x),b\}^{a} \cdot \{ a,b\} \cdot T_{f}(x,y)^{ba} \cdot ba \cdot \sigma_{\{ x,y\}} \left((ba)^{-1} \right)\cdot \sigma_{\{ x,y\}} \left( \{ a,l(y)\}^{b} \right) \cdot l(\{ x,y\}) \\
&= \{ l(x),b\}^{a} \cdot \{ a,b\} \cdot ba \cdot T_{f}(x,y) \cdot \sigma_{\{ x,y\}} \left( a^{-1} \cdot \{ a,l(y)\} \cdot b^{-1} \right) \cdot l(\{ x,y\}) \\
&= \tau_{x}(b)^{a} \cdot \{ a,b\} \cdot ba \cdot T_{f}(x,y) \cdot  \sigma_{\{ x,y\}} \left( a^{-1} \cdot \nu_{y}(a) \cdot b^{-1} \right) \cdot l(\{ x,y\})
\end{align*}
where we have defined $\tau_{x}(a):= \{ l(x), a \}$, $\nu_{x}(a):= \{ a, l(x) \}$ and $T_{f}(x,y) := \{ l(x),l(y)\} \cdot l(\{ x,y\})^{-1}$. Note we have the following properties
\begin{enumerate}	
	
	\item[(a)] $T(x,1) = T(1,x) = 1$,
	
	\item[(b)] $T_{f}(x,x) = T_{f}(x,1) = T_{f}(1,x) = 1$,
	
	\item[(c)] $\tau_{x}(1) = \tau_{1}(a) = \nu_{x}(1) = \nu_{1}(a) = 1$,
	
\end{enumerate}	
directly from the definitions; in addition, in the group operation we have the standard calculation
\begin{enumerate}
		
	\item[(d)] $T(x,y) \cdot T(x,yz) = \sigma_{x} \left( T(y,z) \right) \cdot T(xy,z)$,
		
	\item[(e)] $T(x,y) \cdot \sigma_{xy} \left( a \right) = \sigma_{x}(\sigma_{y}(a)) \cdot T(x,y)$.
		
\end{enumerate}
Identity (d) above comes from enforcing associativity in the product of coset representatives $l(x)l(y)l(z)$ and rewriting in terms of the functions $\sigma$ and $T$. Similar identities can be derived by evaluating Eq~\eqref{(2.1)}-Eq~\eqref{(2.5)} on the coset representatives. This process can be formalized in the following definitions.


\begin{definition}
A 2-\emph{cocyle} appropriate for datum $(Q,I)$ is a sequence $T = \{T,T_{f}, \sigma, \tau, \nu \}$ with functions
\begin{enumerate}
			
	\item $T : Q^{2} \rightarrow I$,
			
	\item $T_{f} : Q^{2} \rightarrow I$
			
	\item $\sigma : Q \rightarrow \mathrm{Aut} \, I$,
			
	\item $\tau : Q \times I \rightarrow I$,
			
	\item $\nu : I \times Q \rightarrow I$,
			
\end{enumerate}
such that $T(x,1) = T(1,x) = 1$, $T_{f}(x,x) = T_{f}(x,1) = T_{f}(1,x) = 1$ and $\tau_{x}(1) = \nu_{x}(1) = 1$.
\end{definition}


We will also refer to $\{ \sigma, \tau, \nu \}$ as \emph{action} terms of the datum.


\begin{definition}
Let $T$ a 2-cocycle appropriate for datum $(Q,I)$. The algebra $I \rtimes_{T} Q$ has universe the direct product set $I \times Q$ with operations
\begin{enumerate}
			
	\item $\left\langle a, x \right\rangle \cdot \left\langle b, y \right\rangle := \langle a \cdot \sigma_{x}(b) \cdot T(x,y) \ , \ x \cdot y \rangle $
			
	\item $\left\langle a, x \right\rangle^{-1} := \langle \sigma_{x^{-1}} \big( T(x,x^{-1})^{-1} \cdot a^{-1} \big) \ , \ x^{-1} \rangle$
			
	\item $1 := \left\langle 1, 1 \right\rangle$ 
			
	\item $\left\{ \left\langle a, x \right\rangle, \left\langle b, y \right\rangle \right\} := \left\langle \tau_{x}(b)^{a} \cdot \{a,b\} \cdot ba \cdot T_{f}(x,y) \cdot \sigma_{\{ x, y\}} \left( a^{-1} \cdot \nu_{y}(a) \cdot b^{-1} \right) \ , \ \{ x,y\} \right\rangle$
			
\end{enumerate}
\end{definition}


From the above definition and the properties of a 2-cocycle, it can be easily seen that $i: I \rightarrow I \rtimes_{T} Q$ is an embedding in which $i(I) = \{ \left\langle a, 1 \right\rangle : a \in I \} \triangleleft I \rtimes_{T} Q$ and projection onto the second-coordinate determines an extension 
\begin{align}\label{eqn:canonicalextension}
1 \longrightarrow I \stackrel{i}{\longrightarrow} I \rtimes_{T} Q \stackrel{p_{2}}{\longrightarrow} Q \longrightarrow 1 
\end{align}
with section $l(x) = \left\langle 1, x \right\rangle$. The above extension is a semidirect product if and only if $T_{f}(x,y) = T(x,y) = 1$ for all $x,y \in Q$. By direct calculation we also see that the action terms of the 2-cocycle are represented by the operation of the algebra in the following manner: 
\begin{align*}
\{ \left\langle 1, x \right\rangle , \left\langle a , 1 \right\rangle \} &= \left\langle \tau_{x}(a) , 1 \right\rangle &  \{ \left\langle a, 1 \right\rangle , \left\langle 1 , x \right\rangle \} &= \left\langle \nu_{x}(a) , 1 \right\rangle & \left\langle a , 1 \right\rangle^{\left\langle 1 , x \right\rangle} &= \left\langle \sigma_{x}(a) , 1 \right\rangle .
\end{align*}


\begin{definition}
Let $T$ a 2-cocycle appropriate for datum $(Q,I)$. An extension $\pi : A \rightarrow Q$ \emph{realizes} the 2-cocycle $T$ if there is an extensions $\pi : A \rightarrow Q$ and an isomorphism $i: I \rightarrow \ker \pi$ and a section $l: Q \rightarrow A$ of $\pi$ with $l(1)=1$ such that 
\begin{enumerate}
			
			\item[(R1)] $i \circ T(x,y) = l(x) \cdot l(y) \cdot l(xy)^{-1}$ 
			
			\item[(R2)] $i \circ T_{f}(x,y) = \{ l(x), l(y)\} \cdot l(\{ x,y\})^{-1}$ 
			
			\item[(R3)] $i ( \sigma_{x}(a) ) = i(a)^{l(x)}$
			
			\item[(R4)] $i ( \tau_{x}(a)) = \{ l(x), i(a)\} $
			
			\item[(R5)] $i ( \nu_{x}(a) ) = \{ i(a), l(x)\}$.
			
\end{enumerate}
\end{definition}


\begin{lemma}\label{lem:realizesdatum}
Let $T$ be a 2-cocycle appropriate for datum $(Q,I)$. Then $I \rtimes_{T} Q$ is an extension which realizes $T$. 
\end{lemma}


We see that the algebra $I \rtimes_{T} Q$ realizes the 2-cocycle, but it may not be a multiplicative Lie algebra. If we enforce associativity of the group operation and the equations Eq~\eqref{(2.1)}-Eq~\eqref{(2.5)} to guarantee that $I \rtimes_{T} Q \in \mathcal M \mathcal L$, then for each identity $t=s$ when evaluated in the algebra $I \rtimes_{T} Q$ produces two equations by equating the coordinates; that is, there will be terms $\omega_{t}$ and $\omega_{s}$ built from the operations of the 2-cocycle $T$ such that 
\begin{align}\label{eqn:800}
\left\langle \omega_{t}, t^{Q} \right\rangle = t^{I \rtimes_{T} Q} = s^{I \rtimes_{T} Q} = \left\langle \omega_{s}, s^{Q} \right\rangle .
\end{align}
Since $Q \in \mathcal M\mathcal L$ we already know that $t^{Q} = s^{Q}$ for any identity $t =s \in \Id \mathcal M \mathcal L$; thus, $I \rtimes_{T} Q \in \mathcal M \mathcal L$ if and only if $\omega_{t} = \omega_{s}$ for all identities $t=s \in \Id \mathcal M \mathcal L$. By equating the first coordinate in Eq~\eqref{eqn:800} for associativity and each of the identities Eq~\eqref{(2.1)}-Eq~\eqref{(2.5)} we arrive at the definition of a compatible 2-cocycle.


\begin{definition}
A 2-cocycle $T$ appropriate for datum $(Q,I)$ is a \emph{compatible} 2-cocycle if the following hold:
\begin{enumerate}
			
	\item $T(x,y) \cdot T(xy,z) = \sigma_{x}(T(y,z)) \cdot T(x,yz)$;
			
	\item $\sigma_{x}(\sigma_{y}(a)) \cdot T(x,y)= T(x,y) \cdot \sigma_{xy}(a)$;
			
	\item $\tau_{x}(a) \cdot \nu_{x}(a) = 1$;
			
	\item $\tau_x(\sigma_y(c)T(y,z))^{ab}\cdot \{ a,b\sigma_y(c)T(y,z)\} \cdot b\sigma_y(c)T(y,z)a \cdot T_f(x,yz) \sigma_{\{ x,yz\}} \Big( a^{-1}\cdot \nu_{yz}(a) \cdot \big( b\sigma_y(c)T(y,z) \big)^{-1} \Big) = \{ a,b\} \cdot ba \cdot T_f(x,y) \sigma_{\{ x,y\}}(a^{-1}\cdot \nu_y(a) \cdot b^{-1} ) \sigma_{\{ x,y\}} \Big( b \sigma_{y} \big(\tau_x(c)^a \{ a,c\}\cdot ca \cdot T_f(x,z)\cdot \sigma_{\{ x,z\}}(a^{-1}\nu_z(a)c^{-1}) \big) \cdot \\ T(y,\{ x,z\}) \cdot T(y\{x,z\}, y^{-1}) \cdot \sigma_{\{ x,z\}^y} \big( T(y,y^{-1})^{-1}b^{-1} \big) \Big) \cdot T(\{ x,y\},{\{ x,z\}^y});$
			
	\item $ \tau_{xy}(c)^{a\sigma_x(b)T(x,y)} \cdot \{ a\sigma_x(b)T(x,y),c\} \cdot ca\sigma_x(b)T(x,y) \cdot T_f(xy,z) \cdot \sigma_{\{ xy,z\}}((a\sigma_x(b)T(x,y))^{-1}\cdot \nu_y(a\sigma_x(b)T(x,y)) \cdot c^{-1} )$ 
			$= a \sigma_x(\tau_y(c)^b \cdot \{ b,c\} \cdot cb \cdot T_f(y,z) \cdot \sigma_{\{ y,z\}}(b^{-1}\nu_z(b)c^{-1} ))\cdot T(x,\{ y,z\})\cdot T(x\{y,z\},x^{-1}) \sigma_{{\{ y,z\}}^x}(T(x,x^{-1})^{-1}\cdot a^{-1}) \cdot \sigma_{\{ y,z\}^x}(\tau_x(c)^a \cdot \{ a,c\} \cdot ca \cdot T_f(x,z) \cdot \sigma_{\{ x,z\}}(a^{-1} \nu_z(a) c^{-1}) \cdot T(\{ y,z\}^x,\{ x,z\});$
	\item $ \tau_{\{x,y\}}(B)^A \cdot \{ A,B\} \cdot BA \cdot T_f(\{ x,y\},z^y) \cdot \sigma_{\{ \{ x,y\},z^y\}}(A^{-1} \cdot \nu_{z^y}(A)\cdot B^{-1}) \cdot 
			\sigma_{\{ \{ x,y\},z^y\}}(\tau_{\{y,z\}}(B')^{A'} \cdot \{ A',B'\} \cdot B'A' \cdot T_f(\{ y,z\},x^z) \cdot \sigma_{\{ \{ y,z\},x^z\}}(A'^{-1} \cdot \nu_{x^z}(A')\cdot B'^{-1})) T(\{ \{ x,y\},z^y\},\{ \{ y,z\},x^z\})  \sigma_{\{ y^x,\{ z,x\}\}}( \tau_{\{z,x\}}(B'')^{A''} \cdot \{ A'',B''\} \cdot B''A'' \cdot T_f(\{ z,x\},y^x) \cdot \sigma_{\{ \{ z,x\},y^x\}}(A''^{-1} \cdot \nu_{^xy}(A'')\cdot B''^{-1}))  T(\{ y^x,\{ z,x\}\}, \{ \{ z,x\},y^x\} ) = 1, 
			$ where $D_A$ defines, respectively $A,A',A''$ by replacing $(a,b,x,y),(b,c,y,z)$ and $(c,a,z,x)$ with $(g,h,r,s):$ $D_A= \tau_r(h)^g\cdot \{ g,h\} \cdot hg \cdot T_f(r,s) \cdot \sigma_{\{ r,s\}}(g^{-1}\cdot \nu_s(g) \cdot h^{-1})$ and $D_B$ defines, respectively $B,B',B''$ by replacing $(b,c,y,z),(c,a,z,x)$ and $(a,b,x,y)$ with $(k,l,u,v):$ $D_B= k \cdot \sigma_v(l)\cdot T(u,v) \cdot T(uv,u^{-1}) \cdot \sigma_{v^u}(T(u,u^{-1})^{-1} \cdot k^{-1})$.
			
	\item $ a \sigma_x(\tau_y(c)^b \cdot \{ b,c\} \cdot cb \cdot T_f(y,z) \cdot \sigma_{\{ y,z\}}(b^{-1}\nu_z(b) c^{-1}) ) \cdot T(x,\{ y,z\}) \cdot \sigma_{x\{ y,z\}}(\sigma_{x^{-1}}(T(x,x^{-1})^{-1} a^{-1})) \cdot \\ T(x\{ y,z\},x^{-1}) = \tau_{y^x}(B)^A \cdot \{ A,B\} \cdot BA \cdot T_f(y^x, z^x) \cdot \sigma_{\{ y^x, z^x\}}(A^{-1}\cdot \nu{z^x}(A) \cdot B^{-1}),$ where $A = a\sigma_x(b) \cdot T(x,y)\cdot \sigma_{xy}(\sigma_{x^{-1}}(T(x,x^{-1})^{-1}\cdot a^{-1}))$ and $B = a\sigma_x(c) \cdot T(x,z)\cdot \sigma_{xz}(\sigma_{x^{-1}}(T(x,x^{-1})^{-1}\cdot c^{-1}))$.
			
\end{enumerate}
\end{definition}


\begin{proposition}\label{prop:extenrep}
Let $(Q,I)$ be datum in $\mathcal M \mathcal L$ and $A \in \mathcal M \mathcal L$. Then $A$ is an extension of $I$ by $Q$ if and only if there is a compatible 2-cocycle $T$ of the datum such that $A \approx I \rtimes_{T} Q$.
\end{proposition}
\begin{proof}
The condition of compatibility for $T$ precisely states that $I \rtimes_{T} Q$ satisfies identities Eq~\eqref{(2.1)}-Eq~\eqref{(2.5)} and so $I \rtimes_{T} Q$ is a multiplicative algebra. Then Lemma~\ref{lem:realizesdatum} states it is an extension of $I$ by $Q$. If $A \in \mathcal M \mathcal L$ is an extension of $I$ by $Q$, then by composing with isomorphisms, we may assume we have an extension $1 \longrightarrow I \stackrel{i}{\longrightarrow} A \stackrel{\pi}{\longrightarrow} Q \longrightarrow 1$. Choose a section $l: Q \rightarrow A$ such that $l(1)=1$. If we make the definitions
\begin{align}\label{eqn:ext1}
T(x,y) &:= i^{-1} \big( l(x) \cdot l(y) \cdot l(x \cdot y)^{-1} \big) &  T_{f}(x,y) &:= i^{-1} \big( \{l(x),l(y)\} \cdot l(\{x,y\})^{-1} \big)
\end{align}
and
\begin{align}\label{eqn:ext2}
\sigma_{x}(a) &:= i^{-1} \big( i(a)^{l(x)} \big) &  \tau_{x}(a) &:= i^{-1} \big( \{l(x),i(a) \} \big) &  \nu_{x}(a) &:= i^{-1} \big( \{i(a), l(x) \} \big),
\end{align}
then the map $\phi(x):= \left\langle  i^{-1}(x \cdot (l \circ \pi(x))^{-1} ) , \pi(x) \right\rangle$ defines an isomorphism $A \approx I \rtimes_{T} Q$.
\end{proof}


The category $\mathcal H^{2}_{\mathcal M \mathcal L}$ has objects the compatible 2-cocycles appropriate to datum in $\mathcal M\mathcal L$ and morphisms are given in the following manner: if $T$ is a 2-cocycle for $(Q,I)$ and $T'$ is a 2-cocycle for $(P,J)$, then a \emph{morphism} is a triple $(\alpha,\beta,h) : T \rightarrow T'$ where $\alpha \in \Hom(I,J)$, $\beta \in \Hom(Q,P)$ and $h \in \mathrm{Map} (Q,J)$ with $h(1)=1$ such that 
\begin{enumerate}
			
	\item[(E1)] $\alpha \circ \sigma_{x}(a) = \left( \sigma_{\beta(x)}' \circ \alpha(a) \right)^{h(x)}$,
			
	\item[(E2)] $\alpha \circ \tau_{x}(a) = \tau_{\beta(x)}'(\alpha(a))^{h(x)} \cdot \{ h(x),\alpha(a)\}$,
			
	\item[(E3)] $\alpha \circ \nu_{x}(a) = \{ \alpha(a),h(x)\} \cdot \nu_{\beta(x)}'(\alpha(a))^{h(x)}$.
	
	\item[(E4)] $ (\alpha \circ T(x,y) ) \cdot h(x \cdot y) = h(x) \cdot \sigma_{\beta(x)}' \left( h(y) \right) \cdot T'(\beta(x),\beta(y))$,
			
	\item[(E5)] $(\alpha \circ T_{f}(x,y) ) \cdot h(\{ x,y\}) = \tau_{\beta(x)}'(h(y))^{h(x)} \cdot \{ h(x),h(y)\} \cdot h(y)h(x) \cdot T_{f}'(\beta(x),\beta(y)) \cdot \sigma_{\{ \beta(x),\beta(y)\}}' \left( h(x)^{-1} \cdot \nu_{\beta(y)}'(h(x)) \cdot h(y)^{-1} \right)$.
			
\end{enumerate}
For the notion of equivalent extensions in Definition~\ref{def:equivexactseq} there is a corresponding notion of equivalence for 2-cocycles. The representation of extensions in Proposition~\ref{prop:extenrep} depends on a choice of section. Different choices of section produce equivalent 2-cocycles.


\begin{definition}\label{def:equiv2cocycle}
Let $T$ and $T'$ be 2-cocycles for the datum $(Q,I)$. Then $T$ and $T'$ are \emph{equivalent} if there is a morphism $(\id_{I},\id_{Q},h): T \rightarrow T'$ in the category $\mathcal H^{2}_{\mathcal M \mathcal L}$.
\end{definition}

We let $\mathcal H^{2}_{\mathcal M \mathcal L}(Q,I)$ denote the restriction of the category $\mathcal H^{2}_{\mathcal M \mathcal L}$ to the compatible 2-cocycles appropriate to datum $(Q,I)$.


\begin{lemma}\label{lem:equiv2cocyc}
If $T$ and $T'$ are 2-cocycles of $(Q,I)$ defined by different liftings of an extension $\pi: A \rightarrow Q$, then $T$ and $T'$ are equivalent.
\end{lemma}
\begin{proof}
Let $l,l':Q \rightarrow I$ be two sections for the extension $\pi: A \rightarrow Q$ with $l(1)=l'(1)=1$. Assume $T$ and $T'$ are respectively determined by lifting $l$ and $l'$ according to the rules Eq~\eqref{eqn:ext1} and Eq~\eqref{eqn:ext2}. Define $h: Q \rightarrow I$ by $l'(x) = h(x) \cdot l(x)$. Then (E1)-(E5) are derived by rewriting the 2-cocycle and action terms for $T$ in terms of T'. 
\end{proof}


This notion of equivalence defines an equivalence relation on the 2-cocycles appropriate to fixed datum. Let $[T]$ denote the equivalence class of the compatible 2-cocycle $T$. The $2^{\mathrm{nd}}$-cohomology of the datum $(Q,I)$ is the set $H^{2}(Q,I)$ of equivalence classes of compatible 2-cocycles.


\begin{theorem}\label{thm:equivbij}
Let $(Q,I)$ be datum and $A$ and $A'$ be extensions of $I$ by $Q$ in $\mathcal M\mathcal L$. Then $A$ and $A'$ are equivalent in $\mathrm{Ext}(Q,I)$ if and only if they realize equivalent compatible 2-cocycles.
\end{theorem}
\begin{proof}
We may assume $A = I \rtimes_{T} Q$ and $A' = I \rtimes_{T'} Q$ for compatible 2-cocycles $T$ and $T'$. Assume $T$ and $T'$ are equivalent which is witnessed by $h: Q \rightarrow I$. Define $\phi: I \rtimes_{T} Q \rightarrow I \rtimes_{T'} Q$ by $\phi(a,x) = \left\langle a \cdot h(x), x \right\rangle$. Then $(\id_{I},\phi,\id_{Q})$ is the isomorphisms which witnesses the equivalence of $A$ and $A'$.

If $A$ and $A'$ are equivalent in $\mathrm{Ext}(Q,I)$, then there is an isomorphism $\phi : A \rightarrow A'$ such that $\pi = \pi' \circ \phi$. If $l: Q \rightarrow A$ is a lifting of $\pi$, then $\phi \circ l: Q \rightarrow A'$ is a lifting of $\pi'$. Then we have that $T$ and $T'$ are compatible 2-cocycles determined by two different liftings and so must be equivalent by Lemma~\ref{lem:equiv2cocyc}.
\end{proof}


\begin{theorem}
The categories $\mathrm{EXT} \, \mathcal M \mathcal L$ and $\mathcal H^{2}(\mathcal M \mathcal L)$ are equivalent. If $(Q,I)$ is datum in $\mathcal M\mathcal L$, then the map $\mathrm{Ext}(Q,I) \ni [A] \longmapsto [T]$ such that $A \approx I \rtimes_{T} Q$ is a bijection between $\mathrm{Ext}(Q,I)$ and $H^{2}(Q,I)$.
\end{theorem}
\begin{proof}	
The second statement is guaranteed by Theorem~\ref{thm:equivbij}. For the first statement, fix a choice of section $l_{A}$ for each object in $\mathrm{EXT} \, \mathcal M \mathcal L$. Take a morphism $(\alpha,\lambda,\beta)$ between exact sequences
Let $T^{A}$ be the 2-cocyce defined by $l_{A}$ and $T^{M}$ defined by $l_{M}$ so that $A \approx I \rtimes_{T^{A}} Q$ and $M \approx J \rtimes_{T^{M}} P$. Define $h: Q \rightarrow J$ by $h(x):= j \big( \lambda(l_{A}(x)) \cdot l_{M}(\beta(x))^{-1} \big)^{-1}$. Then $(\alpha,\beta,h)$ witness (E1)-(E5).

Conversely, from a morphism $(\alpha,\beta,h): T \rightarrow T'$ we construct the corresponding morphism $(\alpha,\lambda,\beta): I \rtimes_{T} Q \rightarrow J \rtimes_{T'} P$ by defining $\lambda(a,x):= \left\langle \alpha(a) \cdot h(x) , \beta(x) \right\rangle$. Conditions (E1)-(E5) are used to show $\lambda$ is a homomorphism.  	
\end{proof}


\vspace{0.2cm}

\subsection{Extensions with abelian kernels}\label{abelianextensions}

We shall now consider the previous machinery in the case of datum $(Q,I)$ in which the kernel algebra $I$ is abelian. In this case, $2^{\mathrm{nd}}$-cohomology naturally forms an abelian group induced by the group operation in the codomain of 2-cocycle terms. Since multiplicative Lie algebras are expansions of groups, they form a weakly associative Mal'cev variety; that is, there is a term $m(x,y,z)$ satisfying the identities
\begin{align*}
x&=m(x,y,y) & x&=m(y,y,x) &  x&=m(m(x,y,z),z,y) .  
\end{align*}
In the present case, we can take the group term $m(x,y,z) := x \cdot y^{-1} \cdot z$. The cohomology for multiplicative Lie extensions with abelian ideals can then be derived by specialization from the more general case of weakly associative difference term varieties exposited in Wires \cite[Sec 4]{wiresI}. This in turn is a specialization of extensions realizing affine datum in arbitrary varieties of universal algebras \cite{wiresI}. In the present section, we derive from the previous section the special form 2-cocycles, 2-coboundaries, derivations and stabilizing automorphism have when the kernel of an extension is abelian.

For an abelian ideal $I \triangleleft A$ with $Q = A/I$ and choice of section $l:Q \rightarrow A$, the representation of the operations simplifies; in this case, 
\begin{align}
a \cdot l(x) \cdot b \cdot l(y) &= \left( a + \sigma_{x}(b) + T(x,y) \right) \cdot l(xy) 
\end{align}
and
\begin{align}
\begin{split}
\{ a \cdot l(x),b \cdot l(y)\} &= \{ l(x),b\}^{a} \cdot \{ a,b\} \cdot T_{f}(x,y)^{ba} \cdot ba \cdot \sigma_{\{ x,y\}} \left( a^{-1} \cdot \{ a,l(y)\} \cdot b^{-1} \right) \cdot l(\{ x,y\}) \\
&= \Big( a + b + \tau_{x}(b) + T_{f}(x,y) + \sigma_{\{ x,y\}} \left( \nu_{y}(a) - a - b \right) \Big) \cdot l(\{ x,y\}).
\end{split}
\end{align}
Note that if $(Q,I)$ is datum with $I$ abelian and $A$ and $A'$ are equivalent extensions realizing the datum, then the associated compatible 2-cocycles have the same action terms. To see this, let $A$ and $A'$ have associated 2-cocycles $T$ and $T'$, respectively. Then we can assume there is a fixed extension $A$ in which $T$ and $T'$ are defined by different liftings $l$ and $l'$, respectively. Then $l(x) = h(x) \cdot l'(x)$ for a map $h: Q \longrightarrow I$ and we have
\begin{align*}
\tau_{x}(a) &= \{ l(x),a\} = \{ h(x) \cdot l'(x),a\} = \{ h(x),a\} + \{ l'(x),a\}^{h(x)} = \{ l'(x),a\} = \tau_{x}'(a) \\
\nu_{x}(a) &= \{ a,l(x)\} = \{ a,h(x) \cdot l'(x)\} = \{ a,l'(x)\}^{h(x)} + \{ a,h(x)\} = \{ a,l'(x)\} = \nu_{x}'(a) .
\end{align*}
Since the equivalence on extensions fixes the action terms $\chi = \{\sigma,\tau,\nu\}$, the equivalence classes of compatible 2-cocycles which realize the fixed pair $(Q,I)$ can themselves be partitioned into different classes according to the action terms which they also realize. We can consider the set $\mathrm{Ext}_{\chi}(Q,I)$ of equivalence classes of extensions which realize $\chi$; that is, we can formally the modify the notion of datum to include the action terms. Since the action terms are fixed, a 2-cocycle of the datum is then considered to only consist of the factor sets and it is then possible to define a group operation on the 2-cocycles induced by the codomains.

For $(Q,I)$ with $I$ abelian, a 2-\emph{cocycle }consists only of the factor sets $T=\{T(x,y), T_{f}(x,y) \}$ corresponding to the group and bracket operation. Given action terms $\chi = \{\sigma, \tau,\nu \}$ for the pair $(Q,I)$, the operations in the algebra $I \rtimes_{\chi,T} Q$ are then given by
\begin{enumerate}
		
	\item $\left\langle a,x \right\rangle \cdot \left\langle b,y \right\rangle = \left\langle a + \sigma_{x}(b) + T(x,y), x \cdot y \right\rangle$
		
	\item $\left\{ \left\langle a,x \right\rangle , \left\langle b,y \right\rangle \right\} =\left\langle a + b + \tau_{x}(b) + T_{f}(x,y) + \sigma_{[x,y]} \left( \nu_{y}(a) - a - b \right) \ , \ \{ x,y\} \right\rangle$.
		
\end{enumerate}
We can then refer to the triple $(Q,I,\chi)$ as \emph{datum} and consider extensions which realize the datum. The representation of extensions with abelian ideals by $I \rtimes_{\chi,T} Q$ yields a corresponding decomposition of term operations which can often by useful in discussing the satisfaction of equations in the extension algebras. Let $\tau$ be the signature of the variety $\mathcal M\mathcal L$ and let $\tau^{\mathrm{gr}}$ be the reduct to the group operation symbols.


\begin{lemma}\label{lem:equationcondition}
Let $T$ be a compatible 2-cocycle for datum $(Q,I,\chi)$ in $\mathcal M\mathcal L$ with $I$ abelian. For every term $t(\vec{x})$ in the signature of $\mathcal M\mathcal L$, there is term $t^{\partial}$ in the signature $\tau^{\mathrm{gr}}$, a term $t^{\chi}$ in the signature $\chi \cup \tau$, and a term $t^{T}$ in the signature $T \cup \chi \cup \tau$ such that
\begin{align}\label{eqn:eqnrep}
t^{I \rtimes_{\chi,T} Q} \left( \left\langle a_{1} , x_{1} \right\rangle, \ldots, \left\langle a_{n} , x_{n} \right\rangle \right) = \left\langle \ t^{\partial}(\vec{a}) + t^{\chi}(\vec{a},\vec{x}) + t^{T}(\vec{x}) \ , \ t^{Q}(\vec{x}) \ \right\rangle .
\end{align} 
The term $t^{T}(\vec{x})$ can be chosen so that $t^{0}(\vec{x}) \equiv 1$ where 0 signifies the trivial 2-cocycle, and $t^{\chi}(\vec{a},\vec{x}) \equiv 1$ if $\chi = \{ \id, 1, 1 \}$.
\end{lemma}
\begin{proof}
Inductively on the generation of terms. On evaluation of $t^{I \rtimes_{\chi,T} Q}$, the term $t^{T}$ collects all the nested operations which utilize the factor sets $\{ T(x,y),T_{f}(x,y) \}$, the term $t^{\chi}$ collects those terms which utilize action terms but no factor sets, and $t^{\partial}$ is the remaining term of pure group operations.  
\end{proof}


If we are given that $I,Q \in \mathcal M\mathcal L$ and $I \rtimes_{\chi,T} Q \in \mathcal M\mathcal L$, then for any $t=s \in \Id \mathcal M\mathcal L$ we have that 
\begin{align}\label{eqn:sat}
t^{\partial}(\vec{a}) + t^{\chi}(\vec{a},\vec{x}) + t^{T}(\vec{x}) = s^{\partial}(\vec{a}) + s^{\chi}(\vec{a},\vec{x}) + s^{T}(\vec{x})
\end{align} 
for all appropriate tuples $\vec{a}$ in $I$ and $\vec{x}$ in $Q$. Taking the identity tuple $\vec{a}=(1,\ldots,1)$ we conclude that $t^{T}(\vec{x}) = s^{T}(\vec{x})$ and subtracting this from Eq~\eqref{eqn:sat} we have $t^{\partial}(\vec{a}) + t^{\chi}(\vec{a},\vec{x}) = s^{\partial}(\vec{a}) + s^{\chi}(\vec{a},\vec{x})$. We say that a 2-cocycle $T = \{ T(x,y), T_{f}(x,y) \}$ is \emph{compatible} if $t^{T}(\vec{x}) = s^{T}(\vec{x})$ for all $t=s \in \Id \mathcal M\mathcal L$. We say the action terms $\chi = \{\sigma, \tau,\nu \}$ are \emph{compatible} if $t^{\partial}(\vec{a}) + t^{\chi}(\vec{a},\vec{x}) = s^{\partial}(\vec{a}) + s^{\chi}(\vec{a},\vec{x})$ for all $t=s \in \Id \mathcal M\mathcal L$. Then as a consequence of Lemma~\ref{lem:equationcondition} we see the following:
\begin{itemize}

	\item the semidirect product $I \rtimes_{\chi} Q \in \mathcal M\mathcal L$ if and only if $\chi$ is compatible; 
	
	\item for fixed compatible $\chi$, taking $T$ as a parameter which runs through all compatible 2-cocycles, the algebras $I \rtimes_{\chi,T} Q$ reconstruct all extensions which realize datum $(Q,I,\chi)$.

\end{itemize}


If we then take $T$ as a parameter which runs through all compatible 2-cocycles, the algebras $I \rtimes_{\chi,T} Q$ reconstruct all extensions which realize datum $(Q,I,\chi)$. The compatible 2-cocycle equations take the form 
\begin{enumerate}
		
	\item $T(x,y) + T(xy,z) = \sigma_{x}(T(y,z)) + T(x,yz)$
		
	\item $\tau_x(\sigma_y(c) + T(y,z)) + \sigma_y(c)+ T(y,z) + T_f(x,yz) + \sigma_{\{ x,yz\}} \Big( \nu_{yz}(a) - a - b - \sigma_y(c) - T(y,z)^{-1} \Big) = T_f(x,y) + \sigma_{\{ x,y\}}( \nu_y(a) -a) + \sigma_{\{ x,y\}} \Big(\sigma_y \big(\tau_x(c) + c + a + T_f(x,z) + \sigma_{\{ x,z\}} ( \nu_z(a) -a - c) \big) \Big) + T(y,\{ x,z\}) + T(y \{x, z\}, y^{-1}) - \sigma_{{\{ x,z\}}^y}( T(y,y^{-1}) + b) + T(\{ x,y\},\{ x,z\}^y)$
		
	\item $\tau_{xy} (c) 
		+ c + T(x,y) + T_f(xy,z) + \sigma_{\{ xy,z\}} (\nu_y (a + \sigma_x(b) + T(x,y)) - a- c - \sigma_x(b) - T(x,y)) = \sigma_x (\tau_y(c) + c + T_f(y,z) + \sigma_{\{ y,z\}} ( \nu_z(b) -b - c)) + T(x,\{ y,z\}) + T(x\{y,z\}, x^{-1}) - \sigma_{{\{ y,z\}}^x}(T(x,x^{-1})) + \sigma_{{\{ y,z\}}^x} (\tau_x(c)+ c + T_f(x,z) + \sigma_{\{ x,z\}} ( \nu_z(a) - a - c)) + T(\{ y,z\}^x, \{ x,z\})$
		
	\item $ \tau_{\{x,y\}}(B) + B + A + T_f(\{ x,y\},z^y) + \sigma_{\{ \{ x,y\},z^y\}}( \nu_{z^y}(A) - A - B) +
		\sigma_{\{ \{ x,y\},z^y\}}(\tau_{\{y,z\}}(B') + B' + A' + T_f(\{ y,z\},x^z) + \sigma_{\{ \{ y,z\},x^z\}}( \nu_{x^z}(A') - A' - B')) + T(\{ \{ x,y\},z^y\},\{ \{ y,z\},x^z\}) +  \sigma_{\{ y^x,\{ z,x\}\}}( \tau_{\{z,x\}}(B'') + B'' + A'' + T_f(\{ z,x\},y^x) + \sigma_{\{ \{ z,x\},y^x\}}( \nu_{y^x}(A'') - A'' - B'')) + T(\{ y^x,\{ z,x\}\}, \{ \{ z,x\},y^x\} ) = 1, 
		$ where $D_A$ defines, respectively $A,A',A''$ by replacing $(a,b,x,y),(b,c,y,z)$ and $(c,a,z,x)$ with $(g,h,r,s):$ $D_A= \tau_r(h) + h + g + T_f(r,s) + \sigma_{\{ r,s\}}( \nu_s(g) - g - h)$ and $D_B$ defines, respectively $B,B',B''$ by replacing $(b,c,y,z),(c,a,z,x)$ and $(a,b,x,y)$ with $(k,l,u,v):$ $D_B= k + \sigma_v(l) + T(u,v) + T(uv,u^{-1}) - \sigma_{v^u}(k + T(u,u^{-1}))$.
		
	\item $a + \sigma_x(\tau_y(c) + c+b + T_f(y,z) + \sigma_{\{ y,z\}}(\nu_z(b) -b - c)) + T(x,\{ y,z\}) + T(x\{ y,z\},x^{-1}) - \sigma_{{\{ y,z\}}^x}(T(x,x^{-1}) + a) = \tau_{y^x}(B) + B + A + T_f(y^x, z^x) + \sigma_{\{ y^x, z^x\}}(\nu_{z^x}(A) -A - B),$ where $A = a + \sigma_x(b) + T(x,y) + T(xy, x^{-1}) - \sigma_{y^x}(T(x,x^{-1}) + a)$ and $B = a + \sigma_x(c) + T(x,z) + T(xz, x^{-1}) - \sigma_{z^x}(T(x,x^{-1}) + a)$.
		
\end{enumerate}
	

By taking only the parts of (E1)-(E5) concerning factor sets and simplifying when $I$ is abelian, we arrive at the definition of 2-coboundaries.

\begin{definition}
Let $(Q,I,\chi)$ be datum with $I$ abelian. Then a sequence $G = \{ G(x,y),G_{f}(x,y) \}$ is a 2-\emph{coboundary} of the datum $(Q,I,\chi)$ if there is a map $h : Q \longrightarrow I$ such that
\begin{enumerate}
			
	\item $G(x,y) = h(x) + \sigma_{x}(h(y)) - h( x \cdot y )$,
			
	\item $G_{f}(x,y) = \tau_{x}(h(y)) + h(x) + h(y) + \sigma_{\{ x,y\}} \Big( \nu_{y}(h(x)) - h(x) - h(y) \Big) - h(\{ x,y \})$.
			
\end{enumerate}
\end{definition}


For compatible 2-cocycles $T$ and $T'$, the equivalence of Definition~\ref{def:equiv2cocycle} then translates to asserting that $T - T'$ is a 2-coboundary. Equivalence respects the addition on 2-cocycles induced by the codomain.

	
\begin{theorem}
Let $(Q,I,\chi)$ be datum in $\mathcal M\mathcal L$ with $I$ abelian. Then $H^{2}(Q,I,\chi)$ is an abelian group and the map
\begin{align*}
\mathrm{Ext}_{\chi}(Q,I) \ni [A] \longmapsto [T] \in H_{2}(Q,I,\chi)
\end{align*}
such that $A \approx I \rtimes_{\chi,T} Q$ is a bijective correspondence.
\end{theorem}


\begin{definition}
Let $(Q,I,\chi)$ be datum in $\mathcal M\mathcal L$ with $I$ abelian. A \emph{derivation} (or 1-cocycle) of the datum is witness for the null 2-coboundary; that is, a map a map $d: Q \longrightarrow I$ such that
	\begin{enumerate}
			
		\item $d( x \cdot y ) = d(x) + \sigma_{x}(d(y))$
			
		\item $d(\{ x,y\}) = \tau_{x}(d(y)) + d(x) + d(y) + \sigma_{\{ x,y\}} \Big( \nu_{y}(d(x)) - d(x) - d(y) \Big)$.
			
	\end{enumerate}
\end{definition}


It is easy to see that the set of derivations of the datum $\mathrm{Der} (Q,I,\chi)$ forms an abelian group under point-wise addition in the codomain. For an extension
\begin{align}
1 \longrightarrow I \stackrel{i}{\longrightarrow} A \stackrel{\pi}{\longrightarrow} Q \longrightarrow 1
\end{align}
a stabilizing automorphism is $\gamma \in \Aut A$ such that $(\id_{I},\gamma,\id_{Q})$ is an automorphism	in $\mathrm{EXT}(Q,I)$. The stabilizing automorphisms form a subgroup isomorphic to the abelian group of derivations for the associated datum.


\begin{theorem}\label{thm:stab}
Let $(Q,I,\chi)$ be datum in $\mathcal M\mathcal L$ with $I$ abelian and $[A] \in \mathrm{Ext}_{\chi}(Q,I)$. Then $\mathrm{Stab}(A) \approx \mathrm{Der} (Q,I,\chi)$.
\end{theorem}
\begin{proof}
The isomorphism between the group of derivations and stabilizing automorphisms follows the familiar group case. Given $\gamma \in \mathrm{Stab}(A)$, the fact that $(\id_{I},\gamma,\id_{Q})$ is an automorphism in $\mathrm{EXT}(Q,I)$ implies $\gamma$ restricts to the identity on $I$. Then there must be a function $d_{\gamma} : Q \rightarrow I$ such that $\left\langle a, 1 \right\rangle^{-1} \cdot \gamma(a,x) = \gamma(1,x) = \left\langle d_{\gamma}(x), x \right\rangle$. Define the correspondence $\gamma \longmapsto d_{\gamma}$ where $\gamma(a,x) = \left\langle a + d_{\gamma}(x), x \right\rangle$. Verification offers no surprises.  
\end{proof}


\begin{remark}
If $(Q,I,\chi)$ is datum in $\mathcal M\mathcal L$ with $I$ abelian, then evaluation of the following identities in $I \rtimes_{\chi} Q$
\begin{align*}
[ \left\langle 1, x \right\rangle \cdot \left\langle 1, y \right\rangle, \left\langle a, 1 \right\rangle ] &= [ \left\langle 1, y \right\rangle, \left\langle a, 1 \right\rangle ]^{\left\langle 1, x \right\rangle} \cdot [ \left\langle 1, x \right\rangle, \left\langle a, 1 \right\rangle ] \\
[ \left\langle a, 1 \right\rangle , \left\langle 1, x \right\rangle \cdot \left\langle 1, y \right\rangle ] &= [ \left\langle a, 1 \right\rangle , \left\langle 1, x \right\rangle ] \cdot [ \left\langle a, 1 \right\rangle, \left\langle 1, y \right\rangle ]^{\left\langle 1, x \right\rangle}
\end{align*}
yields
\begin{enumerate}

	\item $\tau_{xy}(a) = \tau_{x}(a) + \sigma_{x}(\tau_{y}(a))$,
	
	\item $\nu_{xy}(a) = \nu_{x}(a) + \sigma_{x}(\nu_{y}(a))$

\end{enumerate}
Then (1)-(2) show that the functions $x \mapsto \tau_{x}(a)$ and $x \mapsto \nu_{x}(a)$ for fixed $a \in I$ are derivations of the group reduct. By the correspondence in Theorem~\ref{thm:stab}, for each $a \in I$ the functions defined by 
\begin{align*}
	\alpha_{a}(b,x) &:= \left\{\left\langle a , 1 \right\rangle , \left\langle b , x \right\rangle \right\} \cdot \left\langle b , x \right\rangle = \left\langle b + \nu_{x}(a) , x \right\rangle \\
	\delta_{a}(b,x) &:= \left\{\left\langle b , x \right\rangle , \left\langle a , 1 \right\rangle \right\} \cdot \left\langle b , x \right\rangle = \left\langle b + \tau_{x}(a), x \right\rangle
\end{align*}
are stabilizing automorphisms of the group reduct. 
\end{remark}


From Proposition~\ref{prop:commutator}, we see in an extension $I \rtimes_{\chi,T} Q$ that $I$ is central if and only if $\chi$ is \emph{trivial}: that is, $\sigma_{x}(a)=a$, $\tau_{x}(a) = \nu_{x}(a) = 1$ for all $x \in Q$, $a \in I$. The following description of nilpotent algebras follows directly from \cite[Sec 2]{wiresI} since $\mathcal M\mathcal L$ is an example of a variety with a difference term.


\begin{proposition}
Let $A \in \mathcal M\mathcal L$. Then $A$ is n-step nilpotent if and only if there exists abelian algebras $Q_{i}$ for $i=1,\ldots,n$ and compatible 2-cocycles $T_{i}$ appropriate to each $(Q_{i},Q_{i+1})$ such that $A \approx Q_{n} \rtimes_{T_{n-1}} \cdots \rtimes_{T_{1}} Q_{1}$. 
\end{proposition}


The derived series of the ideal $I$ is $I^{0}=I$, $I^{1}=C(I,I)$ and $I^{k+1} = C(I^{k},I^{k})$. The ideal $I$ is n-step solvable if $I^{n} = 0$.

\begin{proposition}
Let $A \in \mathcal M\mathcal L$. Then $A$ is n-step solvable if and only if $A \approx Q_{n} \rtimes_{\chi_{n-1},T_{n-1}} \cdots \rtimes_{\chi_{1},T_{1}} Q_{1}$ where each $Q_{i}$ is abelian.
\end{proposition}
\begin{proof}
If we set $Q_{k} = A^{k-1}/A^{k}$, then there is an action $\chi_{k}$ and 2-cocycle $T_{k}$ such that $A/A^{k+1} \approx A^{k}/A^{k+1} \rtimes_{\chi_{k},T_{k}} A/A^{k} = Q_{k+1} \rtimes_{\chi_{k},T_{k}} A/A^{k}$ where $Q_{k+1}$ is abelian. The product representation follows recursively.

Now assume $A \approx Q_{n} \rtimes_{\chi_{n-1},T_{n-1}} \cdots \rtimes_{\chi_{1},T_{1}} Q_{1}$ where each $Q_{i}$ is abelian. Set $A_{k} = Q_{k} \rtimes_{\chi_{k-1},T_{k-1}} \cdots \rtimes_{\chi_{1},T_{1}} Q_{1}$. The second-projection $p_{2} : A_{k} = Q_{k} \rtimes_{\chi_{k-1},T_{k-1}} A_{k-1} \rightarrow A_{k-1}$ realizes an extension in which $Q_{k} = \ker p_{2}$ is an abelian ideal. Then 
\begin{align*} 
\big( C(A_{2},A_{2}) \vee Q_{2} \big)/Q_{2} = C \big( \big( A_{2} \vee Q_{2}/Q_{2} \big) , \big( A_{2} \vee Q_{2}/Q_{2} \big) \big) = C(Q_{1},Q_{1}) = 0
\end{align*}
because $Q_{1}$ is abelian; therefore, $C(A_{2},A_{2}) \leq Q_{2}$. This implies $A_{2}^{2}=0$; that is, $A^{2}$ is 2-step solvable. The proof is completed by induction since $A_{n} = A$.
\end{proof}


\begin{lemma}\label{lem:52}
Let $T$ be a compatible 2-cocycle for datum $(Q,I,\chi)$ in $\mathcal M\mathcal L$. For every term $t(\vec{x})$ in the signature of $\mathcal M\mathcal L$ and $\vec{a} \in \zeta(I)^{\ar t}$,
\begin{align*}
t^{I \rtimes_{\chi,T} Q} \left( \left\langle \ \vec{a} \ , \ \vec{x} \ \right\rangle \right) = \left\langle \ t^{\partial}(\vec{a}) + t^{\chi}(\vec{a},\vec{x}) \ , \ 1 \ \right\rangle \cdot t^{I \rtimes_{\chi,T} Q} \left( \left\langle \vec{1} \ , \ \vec{x} \right\rangle \right)
\end{align*}
\end{lemma}
\begin{proof}
Inductively on the generation of terms.
\end{proof}

\vspace{0.2cm}

\subsection{Correspondence}\label{sec:correspondence}

In this section, we will prove that if $\chi= (\sigma,\tau,\nu)$ are action terms for datum $(Q,I)$ which are realized by an extension, then the equivalence classes of extensions over $(Q,I)$ which realize $\chi$ is in bijective correspondence with the abelian $2^{\mathrm{nd}}$-cohomology group with coefficients in the algebraic center $\zeta(I)$. We first define a category whose objects are the sequence of action terms. The morphisms in the category determine an equivalence on objects which corresponds to the equivalence on extensions; in particular, an extension yields action terms from the induced 2-cocycle and the equivalence class of an extension corresponds to an equivalence class of action terms in the category.

The category $\mathrm{Act} \mathcal M\mathcal L$ has \emph{objects} which are sequences $\chi=(\sigma,\tau,\nu)$ where
\begin{itemize}
			
	\item $\sigma : Q \rightarrow \Aut I$,
			
	\item $\tau: Q \longrightarrow \mathrm{Map}(I,I)$,
			
	\item $\nu: Q \longrightarrow \mathrm{Map}(I,I)$,
			
\end{itemize}
for algebras $(Q,I) \in \mathcal M\mathcal L$ such that there exists a function $S : Q^{2} \rightarrow I$ satisfying
\begin{enumerate}

	\item[(A1)] $\sigma_{x}(\sigma_{y}(a)) \cdot S(x,y) = S(x,y) \cdot \sigma_{xy}(a)$;

	\item[(A2)] $S(x,y) \cdot \tau_{xy}(a) = \sigma_{x} \big( \tau_{y}(a) \big) \cdot \tau_{x}(a) \cdot \big\{ S(x,y)^{-1} , a^{S(x,y)} \big\} \cdot S(x,y)$;
	
	\item[(A3)] $S(x,y) \cdot \nu_{xy}(a) = \big\{ a^{S(x,y)^{-1}} , S(x,y)^{-1} \big\} \cdot \nu_{x}(a) \cdot \sigma_{x} \big( \nu_{y}(a) \big) \cdot S(x,y)$;

\end{enumerate}
Let $\chi^{1} = (\sigma^{1},\tau^{1},\nu^{1})$ be an action of $(Q,I)$ and $\chi^{2} = (\sigma^{2},\tau^{2},\nu^{2})$ and be an action of $(P,K)$. A \emph{morphism} from $\chi^{1}$ to $\chi^{2}$ is a triple $(\alpha,\beta,h) : \chi^{1} \rightarrow \chi^{2}$ where
	\begin{itemize}
		
		\item $\alpha \in \Hom(Q,P)$,
		
		\item $\beta \in \Hom(I,K)$,
		
		\item $h \in \mathrm{Map} (Q,K)$
		
	\end{itemize}
	such that
	\begin{enumerate}
		
		\item $\beta \circ \sigma_{x}^{1}(a) = \left( \sigma_{\alpha(x)}^{2} \circ \beta(a) \right)^{h(x)}$,
		
		\item $\beta \circ \tau_{x}^{1}(a) = \left( \tau_{\alpha(x)}^{2} \circ \beta(a) \right)^{h(x)} \cdot [h(x),\beta(a)]$
		
		\item $\beta \circ \nu_{x}^{1}(a) = \{\beta(a),h(x)\} \cdot \left( \nu_{\alpha(x)}^{2} \circ \beta(a) \right)^{h(x)}$.
		
	\end{enumerate}

Let $\mathrm{Act}(Q,I)$ be the actions of fixed datum $(Q,I)$. The following is the restriction of the notion of equivalence for nonabelian 2-cocycles adapted to the action terms over fixed datum $(Q,I)$.


\begin{definition}\label{def:4.2}
Let $\chi$ and $\chi'$ be action terms for $(Q,I)$. We say $\chi$ is \emph{equivalent} to $\chi'$ and write $\chi \sim \chi'$ if there is a map $h: Q \rightarrow I$ such that 
\begin{enumerate}
			
	\item $\sigma_{x}(a) = h(x) \cdot \sigma_{x}(a) \cdot h(x)^{-1}$,
			
	\item $\tau_{x}(a) = \tau_{x}'(a)^{h(x)} \cdot \{ h(x),a\}$, 
			
	\item $\nu_{x}(a) = \{ a,h(x)\} \cdot \nu_{x}'(a)^{h(x)}$.
			
\end{enumerate}
\end{definition}


\begin{lemma}
The equivalence of action terms of $(Q,I)$ determines an actual equivalence relation.
\end{lemma}
\begin{proof}
Suppose $h(x)$ witness that $\chi$ is equivalent to $\chi'$ and $s(x)$ witness that $\chi'$ is equivalent to $\chi''$. Then 
\begin{align*}
\sigma_{x}(a) = \sigma_{x}'(a)^{h(x)} = \sigma_{x}''(a)^{h(x)s(x)}
\end{align*}
and
\begin{align*}
\tau_{x}(a) = \tau_{x}'(a)^{h(x)} \cdot \{ h(x),a\} &= \left( \tau_{x}''(a)^{s(x)} \cdot \{ s(x),a\} \right)^{h(x)} \cdot \{ h(x),a\} \\
&= \tau_{x}''(a)^{h(x)s(x)} \cdot \{ s(x),a\}^{h(x)} \cdot \{ h(x),a\} \\
&= \tau_{x}''(a)^{h(x)s(x)} \cdot \{ h(x)s(x),a\}.
\end{align*}
Similarly, for $\nu$; thus, $\chi \sim \chi''$. Likewise, if $\chi \sim \chi'$ witnessed by $h$, then $\sigma_{x}'(a) = \sigma_{x}'(a)^{h(x)^{-1}}$. We also see that $\tau_{x}(a) = \tau_{x}'(a)^{h(x)} \cdot \{ h(x),a\}$ implies 
\begin{align*}
\tau_{x}'(a) = \tau_{x}(a)^{h(x)^{-1}} \cdot \left( \{ h(x),a\}^{-1} \right)^{h(x)^{-1}} &= \tau_{x}(a)^{h(x)^{-1}} \cdot \{ a,h(x)\}^{h(x)^{-1}} \\
&= \tau_{x}(a)^{h(x)^{-1}} \cdot \{ h(x)^{-1},a^{h(x)}\}^{h(x)^{-1}} \\
&= \tau_{x}(a)^{h(x)^{-1}} \cdot \{ h(x)^{-1},a\} 
\end{align*}
Similarly, for $\nu$; thus, $\chi' \sim \chi$.
\end{proof}


The equivalence class of action terms $\chi$ over fixed datum $(Q,I)$ is denoted by $[\chi]$. Given an extension $E: 1 \rightarrow I \rightarrow A \rightarrow Q \rightarrow1$, any choice of section $l$ determines a 2-cocycle $T=\{T^{l},T^{l}_{f}, \sigma^{l}, \tau^{l}, \nu^{l} \}$. Then the action terms of the 2-cocycle give $\Psi(A) = [\chi^{l}] =\left[ \{ \phi^{l}, \tau^{l}, \nu^{l} \} \right]$. By taking equivalence classes of action terms, the map 
\begin{align}
\Psi : \mathrm{EXT}(Q,I) \rightarrow \mathrm{Act}(Q,I)/\sim
\end{align}
is well-defined.


\begin{definition}
The equivalence class of action terms $[\chi]$ is realized by an extension $1 \longrightarrow I \stackrel{i}{\longrightarrow} A \stackrel{\pi}{\longrightarrow} Q \longrightarrow 1$ if there is a section $l: Q \rightarrow A$ such that $[\chi] = [\chi^{l}] = \Psi(E)$.
\end{definition}


\begin{remark}
An attentive reader may wonder why the bracket operation is not included in our list of properties for action terms in (A1)-(A3). While it is possible to do that now, it is (1) not strictly necessary for the results of this manuscript, and (2) we feel it is to premature since it is intimately connected with a proper notion of abstract kernel which would form the basis for the future development of $3^{\mathrm{rd}}$-cohomology and obstructions. At that point the category $\mathrm{Act} \mathcal M\mathcal L$ would be replaced by an appropriate category of abstract kernels.
\end{remark}


Let $\chi = \{ \sigma, \tau, \nu\}$ be action terms of datum $(Q,I)$ and $B \leq I$. We say $B$ is \emph{closed} under $\chi$ if $\tau_{x}(B) \subseteq B, \nu_{x}(B) \subseteq B$ and $\sigma_{x}(B) \subseteq B$ for all $x \in Q$. Then $\chi_{B}$ denotes the action terms whose operations are the restriction to $B$.


\begin{lemma}
Let $\chi$ be action terms of datum $(Q,I)$ which is realized by some extension. Then $\mathrm{Null} \, I$ is closed under $\chi$.
\end{lemma}
\begin{proof}
Take a lifting $\sigma : Q \rightarrow \Aut I$ such that $\phi_{x} = \sigma_{x} \mathrm{Inn} \, I$. Then there exists $T(x,y) : Q^{2} \rightarrow I$ such that $\sigma_{xy}(a)^{T(x,y)} = \sigma_{x} \circ \sigma_{y}(a)$; in particular, $a^{T(x,x^{-1})} = \sigma_{x} \circ \sigma_{x^{-1}}(a)$. Take $y \in I$, $a \in \mathrm{Null} \, I$. Then 
\begin{align*}
\{ y , \sigma_{x}(a)\} = \left\{ \left( y^{T(x,x^{-1})^{-1}} \right)^{T(x,x^{-1})} , \sigma_{x}(a) \right\} = \sigma_{x} \left( \left\{ \sigma_{x^{-1}}\left( y^{T(x,x^{-1})^{-1}} \right) , a \right\} \right) = \sigma_{x}(1)=1
\end{align*}
since $a \in \mathrm{Null} \, I$; similarly, $\{ \sigma_{x}(a),y\} = 1$. We conclude $\mathrm{Null} \, I$ is closed under $\sigma$. Lemma~\ref{lem:4} shows $\mathrm{Null} \, I$ is closed under $\tau$ and $\nu$.
\end{proof}


\begin{lemma}\label{lem:zetaclosed}
Let $\chi$ be action terms of datum $(Q,I)$ which is realized by some extension. Then $\zeta(I)$ is closed under $\chi$.
\end{lemma}
\begin{proof}
Note $Z(I^{\mathrm{gr}})$ is a characteristic subgroup of $I^{\mathrm{gr}}$; thus, we see that $Z(I^{\mathrm{gr}})$ is closed under $\sigma$. Let $\chi$ be realized by $A$ and take $a \in \zeta(I)$ and $x \in A$. For $a \in \zeta(I)$ we have $\tau_{x}(a) = \{ l(x),a\} \in \zeta(I)$ and $\nu_{x}(a) = \{a,l(x)\} \in \zeta(I)$ by Lemma~\ref{lem:4}.; thus, $\zeta(I)$ is closed under the action terms. 
\end{proof}


\begin{lemma}\label{lem:50}
Let $T$ be a 2-cocycle for datum $(Q,I)$ and $S$ a 2-cocycle for datum $(Q,\zeta(I))$. Then for every term $t(\vec{x})$ in the signature of $\mathcal M\mathcal L$ we have
\begin{align*}
t^{I \rtimes_{T \cdot S} Q}( \left\langle \ \vec{a} , \vec{x} \ \right\rangle) = \left\langle \ t^{S}(\vec{x}) , 1 \ \right\rangle \cdot t^{I \rtimes_{T} Q}( \left\langle \ \vec{a} , \vec{x} \ \right\rangle) 
\end{align*}
\end{lemma}
\begin{proof}
Inductively by term generation.
\end{proof}


With the appropriate preparations, correspondence for the group case can be extended to include multiplicative Lie algebras.

\begin{theorem}\label{thm:corresp1}(Correspondence)
Assume $(Q,I)$ is datum in $\mathcal M\mathcal L$. Suppose the action terms $\chi$ are realized by an extension of $I$ by $Q$. Then $H^{2}(Q,\zeta(I),\chi_{\zeta(I)})$ and $\mathrm{Ext}_{\chi}(Q,I)$ are in bijective correspondence. 
\end{theorem}
\begin{proof}
Since $\chi$ is realized by an extension, $\zeta(I)$ is closed under $\chi$ by Lemma~\ref{lem:zetaclosed}. We show there is a sharply transitive action of the group $H^{2}(Q,\zeta(I),\chi_{\zeta(I)})$ on the set $\mathrm{Ext}_{\chi}(Q,I)$. Let $\chi = \{ \phi, \tau,\nu \}$. Take $[A] \in \mathrm{Ext}_{\chi}(Q,I)$ and $[S] \in H^{2}(Q,\zeta(I),\chi_{\zeta(I)})$. We can assume $A = I \rtimes_{T} Q$ for a compatible 2-cocycle $T$ defined by a lifting $l : Q \rightarrow I$ with $l(1)=1$; that is, 
		\begin{itemize}
			
			\item $T(x,y) = l(x)l(y)l(xy)^{-1}$,
			
			\item $T_{f}(x,y) = \{ l(x),l(y)\} l(\{ x,y\})^{-1}$,
			
			\item $\sigma^{l}_{x}(a) = a^{l(x)}$,
			
			\item $\tau_{x}(a) = \{ l(x),a\}$,
			
			\item $\nu_{x}(a) = \{ a, l(x)\}$
			
		\end{itemize}
		where $\phi_{x} = \sigma_{x}^{l} \mathrm{Inn} \, I$. Define $A \ast S$ to be the extension $A \ast S := I \rtimes_{T \ast S} Q$ where the 2-cocyce $T \ast S$ is defined by 
		\begin{itemize}
			
			\item $(T \ast S)(x,y) := T(x,y) \cdot S(x,y)$,
			
			\item $(T \ast S)_{f}(x,y) := T_{f}(x,y) \cdot S_{f}(x,y)$.
			
		\end{itemize}
		and the action terms are the same as $T$. We first check this action is well-defined.

		Let us assume $[A] = [A']$ in $\mathrm{Ext}_{\chi}(Q,I)$. Let $A \approx I \rtimes_{T} Q$ and $A' \approx I \rtimes_{T'} Q$. Then there exists $h: Q \rightarrow I$ which witnesses $T \sim T'$; that is, $h(1)=1$ and 
		\begin{enumerate}
			
			\item $\sigma_{x}(a) = h(x) \cdot \sigma'_{x}(a) \cdot h(x)^{-1}$, 
			
			\item $T(x,y) \cdot h(x,y) = h(x) \cdot \sigma'_{x}(h(y)) \cdot T'(x,y)$,
			
			\item $T_{f}(x,y) \cdot h(\{ x,y\}) = \tau_{x}(h(y))^{h(x)} \cdot \{ h(x),h(y)\} \cdot h(y) h(x) \cdot T_{f}'(x,y) \cdot \sigma'_{\{ x,y\}} \big( h(x)^{-1} \cdot \nu_{y}(h(x)) \cdot h(y)^{-1} \big)$,
			
			\item $\tau_{x}(a) = \tau_{x}(a)^{h(x)} \cdot \{ h(x),a\}$,
			
			\item $\nu_{x}(a) = \{ h(x),a\} \cdot \nu_{x}(a)^{h(x)}$.
			
		\end{enumerate}
		If we take $[S] \in H^{2}(Q,\zeta(I),\chi_{\zeta(I)})$, then because $\im S \subseteq \zeta(I)$ we can add $S(x,y)$ to both sides of (2) and (3) above and use commutativity to shift to the right position. Then we can conclude $T \ast S \sim T' \ast S$ and so, $[A \ast S] = [A' \ast S]$.

		Now suppose $[S] = [S']$ in $H^{2}(Q,\zeta(I),\chi_{\zeta(I)})$. Then there is $h : Q \rightarrow \zeta(I)$ such that 
		\begin{itemize}
			
			\item $S(x,y) - S'(x,y) = h(x) + \sigma_{x}(h(y)) - h(xy)$, 
			
			\item $S_{f}(x,y) - S_{f}'(x,y) = \tau_{x}(h(y)) + h(y) + h(x) + \sigma_{\{ x,y\}} \big( \nu_{y}(h(x)) - h(x) - h(y) \big) - h(\{ x,y\})$.
			
		\end{itemize}
		Since $\im h \subseteq \zeta(I)$, the equations
		\begin{align*}
			\sigma_{x}(a) &= h(x) \cdot \sigma_{x}(a) \cdot h(x)^{-1}, &\tau_{x}(a) = \tau_{x}(a)^{h(x)} \cdot \{ h(x),a\},& &\nu_{x}(a) = \{ h(x),a\} \cdot \nu_{x}(a)^{h(x)}
		\end{align*}
		are trivially true. Then for $[A] \in \mathrm{Ext}_{\chi}(Q,I)$ with $A \approx I \rtimes_{T} Q$ we have
		\begin{align*}
			T(x,y) \cdot S(x,y) &= T(x,y) \cdot \big( h(x) + \sigma_{x}(h(y)) - h(xy) + S'(x,y) \big) \\
			&= h(x) \cdot \sigma_{x}(h(y)) \cdot T(x,y) \cdot S'(x,y) \cdot h(x,y)^{-1}
		\end{align*}
		since $\im h, \im S \subseteq \zeta(I)$ and $\im T \subseteq I$. For the same reasons we see that
		\begin{align*}
			T_{f}(x,y) &\cdot S_{f}(x,y) \\
			&= T_{f}(x,y) \cdot \big( \tau_{x}(h(y)) + h(y) + h(x) + \sigma_{\{ x,y\}} \big( \nu_{y}(h(x)) - h(x) - h(y) \big) - h(\{ x,y\}) + S_{f}(x,y) \big) \\
			&= \tau_{x}(h(y))^{h(x)} \cdot \{ h(x),h(y)\} \cdot h(y) h(x) \cdot \big( T_{f}(x,y) \cdot S_{f}'(x,y) \big) \\
			&\quad \cdot \sigma_{\{ x,y\}} \big( h(x)^{-1} \cdot \nu_{y}(h(x)) \cdot h(y)^{-1} \big) \cdot h(\{ x,y\})^{-1}.
		\end{align*}
		We have that $[A \ast S] = [A \ast S']$; altogether, the action is well-defined.

		We now show $T \ast S$ is a compatible 2-cocycle for $(Q,I)$. Since $T$ is a compatible 2-cocycle, we have $I \rtimes_{T} Q \in \mathcal M\mathcal L$ which implies 
		\begin{align}\label{eqn:name1}
			t^{I \rtimes_{T} Q} \left( \left\langle \ \vec{a} , \vec{x} \ \right\rangle \right) = s^{I \rtimes_{T} Q} \left( \left\langle \ \vec{a} , \vec{x} \ \right\rangle \right)
		\end{align}
		for all identities $t=s \in \Id \, \mathcal M\mathcal L$; in particular, 
		\begin{align}\label{eqn:name2}
			\left\langle \ t^{\partial}(\vec{a}) \ , \ \vec{1} \ \right\rangle = t^{I \rtimes_{T} Q} \left( \left\langle \ \vec{a} , \vec{1} \ \right\rangle \right) = s^{I \rtimes_{T} Q} \left( \left\langle \ \vec{a} , \vec{1} \ \right\rangle \right) = \left\langle \ s^{\partial}(\vec{a}) \ , \ \vec{1} \ \right\rangle .
		\end{align}
		If we restrict to a tuple $\vec{a}$ in $\zeta(I)$, then by Lemma~\ref{lem:52} we have
		\begin{align*}
			\left\langle \ t^{\partial}(\vec{a}) + t^{\chi}(\vec{a},\vec{x}) \ , \ 1 \ \right\rangle \cdot t^{I \rtimes_{T} Q} \left( \left\langle \vec{1} \ , \ \vec{x} \right\rangle \right) &= t^{I \rtimes_{T} Q} \left( \left\langle \ \vec{a} \ , \ \vec{x} \ \right\rangle \right) \\
			&= s^{I \rtimes_{T} Q} \left( \left\langle \ \vec{a} \ , \ \vec{x} \ \right\rangle \right) \\
			&= \left\langle \ s^{\partial}(\vec{a}) + s^{\chi}(\vec{a},\vec{x}) \ , \ 1 \ \right\rangle \cdot s^{I \rtimes_{T} Q} \left( \left\langle \vec{1} \ , \ \vec{x} \right\rangle \right).
		\end{align*} 
		Then using Eq~\ref{eqn:name1} to cancel the factors on the right we conclude that $t^{\partial}(\vec{a}) + t^{\chi}(\vec{a},\vec{x}) = s^{\partial}(\vec{a}) + s^{\chi}(\vec{a},\vec{x})$; thus, Eq~\ref{eqn:name2} yields $t^{\chi}(\vec{a},\vec{x}) = s^{\chi}(\vec{a},\vec{x})$. This argument shows that if the action terms $\chi$ are realized by an extension of datum $(Q,I)$ in $\mathcal M\mathcal L$, then the restriction $\chi_{\zeta(I)}$ is compatible with $\mathcal M\mathcal L$. Then for each $[S] \in H^{2}(Q,\zeta(I),\chi_{\zeta(I)})$, $S$ is compatible with $\mathcal M\mathcal L$. It then follows from Lemma~\ref{lem:50} that
		\begin{align*}
			t^{I \rtimes_{T \cdot S} Q}( \left\langle \ \vec{a} , \vec{x} \ \right\rangle) = \left\langle \ t^{S}(\vec{x}) , 1 \ \right\rangle \cdot t^{I \rtimes_{T} Q}( \left\langle \ \vec{a} , \vec{x} \ \right\rangle) = \left\langle \ s^{S}(\vec{x}) , 1 \ \right\rangle \cdot s^{I \rtimes_{T} Q}( \left\langle \ \vec{a} , \vec{x} \ \right\rangle) = s^{I \rtimes_{T \cdot S} Q}( \left\langle \ \vec{a} , \vec{x} \ \right\rangle) 
		\end{align*} 
		for all $t=s \in \Id \, \mathcal M\mathcal L$; that is, $T \cdot S$ is a 2-cocycle compatible with $\mathcal M\mathcal L$.

		Now assume $[A] \ast [S] = [A \ast S] = [A]$ in $\mathrm{Ext}_{\chi}(Q,I)$. There is a map $h: Q \rightarrow I$ such that $h(1)=1$ and
		\begin{enumerate}
			
			\item $\sigma_{x}(a) = h(x) \cdot \sigma_{x}(a) \cdot h(x)^{-1}$, 
			
			\item $T(x,y) \cdot S(x,y) \cdot h(x,y) = h(x) \cdot \sigma_{x}(h(y)) \cdot T(x,y)$,
			
			\item $T_{f}(x,y) \cdot S_{f}(x,y) \cdot h(\{ x,y\}) = \tau_{x}(h(y))^{h(x)} \cdot \{ h(x),h(y)\} \cdot h(y) h(x) \cdot T_{f}(x,y) \cdot \sigma_{\{ x,y\}} \big( h(x)^{-1} \cdot \nu_{y}(h(x)) \cdot h(y)^{-1} \big)$,
			
			\item $\tau_{x}(a) = \tau_{x}(a)^{h(x)} \cdot \{ h(x),a\}$,
			
			\item $\nu_{x}(a) = \{ h(x),a\} \cdot \nu_{x}(a)^{h(x)}$.
			
		\end{enumerate}
		Again, we see that $\im h \subseteq \zeta(I)$. Then (2) yields 
		\begin{align}\label{eqn:500}
			S(x,y) = h(x) \cdot \sigma_{x}(h(y)) \cdot h(xy)^{-1} .
		\end{align}
		Also, (3) simplifies to 
		\begin{align*}
			T_{f}(x,y) \cdot S_{f}(x,y) \cdot h(\{ x,y\}) = \tau_{x}(h(y)) \cdot h(y) h(x) \cdot T_{f}(x,y) \cdot \sigma_{\{ x,y\}} \big( h(x)^{-1} \cdot \nu_{y}(h(x)) \cdot h(y)^{-1} \big) .
		\end{align*} 
		which yields
		\begin{align}\label{eqn:501}
			S_{f}(x,y) = \tau_{x}(h(y)) \cdot h(y) h(x) \cdot \sigma_{\{ x,y\}} \big( h(x)^{-1} \cdot \nu_{y}(h(x)) \cdot h(y)^{-1} \big) \cdot h(\{ x,y\})^{-1}.
		\end{align}
		Then Eq~\ref{eqn:500} and Eq~\ref{eqn:501} shows $[S] = 0$ in $H^{2}(Q,\zeta(I),\chi_{\zeta(I)})$.

At this point, we have shown the action is sharp. We now establish transitivity. Take $[A],[A'] \in \mathrm{Ext}_{\chi}(Q,I)$ and show there is $[S] \in H^{2}(Q,\zeta(I),\chi_{\zeta(I)})$ such that $[A] \ast [S] = [A \ast S] = [A']$. We have $A \approx I \rtimes_{T} Q$ and $A' \approx I \rtimes_{T'} Q$ where $T$ is defined by a lifting $l: Q \rightarrow I$ and $T'$ is defined by a lifting $s : Q \rightarrow I$ which realize $\chi$; that is, $T = T^{l}$, $T' = T^{s}$ and $\left[ \{ \sigma^{l}, \tau^{l}, \nu^{l} \} \right] = [\chi] = \left[ \{ \sigma^{s}, \tau^{s}, \nu^{s} \} \right]$. There is $h: Q \rightarrow I$ such that $h(1)=1$ and
		\begin{enumerate}
			
			\item $\sigma_{x}^{l} (a) = i_{h(x)} \circ \sigma_{x}^{s} (a)$;
			
			\item $\tau_{x}^{l}(a) = \tau_{x}^{s}(a)^{h(x)} \cdot \{ h(x),a\}$;
			
			\item $\nu_{x}^{l}(a) = \{ a,h(x)\} \cdot \nu_{x}^{s}(a)^{h(x)}$.
			
		\end{enumerate}
So we have $l(x) \cdot a \cdot l(x)^{-1} = h(x) \cdot s(x) \cdot a \cdot s(x)^{-1} \cdot h(x)^{-1}$ for all $a \in I$. Define $t(x) = h(x) \cdot s(x)$ which is a section of the extension $A'$; thus, $[T'] = [T^{t}]$ and $\sigma^{t}_{x}(a) = i_{h(x)} \circ \sigma_{x}^{s} (a) = \sigma^{l}_{x}(a)$. Note $\tau_{x}^{l}(a) = [s(x),a]^{h(x)} \cdot [h(x),a] = [h(x)s(x),a] = \tau_{x}^{t}(a)$; similarly, $\nu_{x}^{l}(a) = \nu_{x}^{t}(a)$.

		Define
		\begin{align*}
			S(x,y) &:= T^{l}(x,y)^{-1} \cdot T^{t}(x,y) \\
			S_{f}(x,y) &:= T_{f}^{l}(x,y)^{-1} \cdot T_{f}^{t}(x,y) .
		\end{align*}
		We will show that $S, S_{f} \subseteq \zeta(I)$. We calculate that
		\begin{align}\label{eqn:1001}
			\begin{split}
				i_{T^{l}(x,y)}(a) &= l(x)l(y)l(xy)^{-1} \cdot a \cdot l(xy) l(y)^{-1} l(x)^{-1} \\
				&= l(x)l(y)l((xy)^{-1}) \cdot a^{T^{l}(xy,(xy)^{-1})^{-1}} \cdot l((xy)^{-1}) l(y)^{-1} l(x)^{-1} \\
				&= t(x)t(y)t((xy)^{-1}) \cdot a^{f^{l}(xy,(xy)^{-1})^{-1}} \cdot t((xy)^{-1}) t(y)^{-1} t(x)^{-1} \\
				&= t(x)t(y)t(xy)^{-1} \cdot a^{T^{t}(xy,(xy)^{-1})T^{l}(xy,(xy)^{-1})^{-1}} \cdot t(xy) t(y)^{-1} t(x)^{-1} \\
				&= i_{T^{t}(x,y)} \big( a^{T^{t}(xy,(xy)^{-1})T^{l}(xy,(xy)^{-1})^{-1}} \big) \\
				&= i_{T^{t}(x,y)} \big( i_{T^{t}(xy,(xy)^{-1})T^{l}(xy,(xy)^{-1})^{-1}}(a) \big) .
			\end{split}
		\end{align}
		This implies 
		\begin{align*}
			i_{T^{l}(x,x^{-1})}(a) = i_{T^{t}(x,x^{-1})} \left( a^{T^{t}(1,1)T^{l}(1,1)^{-1}} \right) = i_{T^{t}(x,x^{-1})}(a);
		\end{align*}
		equivalently, $i_{T^{t}(xy,(xy)^{-1})T^{l}(xy,(xy)^{-1})^{-1}}(a) = a$. Then evaluating back into Eq~\eqref{eqn:1001} yields
		\begin{align}\label{eqn:1023}
			i_{T^{l}(x,y)}(a) = i_{T^{t}(x,y)}(a);
		\end{align}
		therefore, $S(x,y) = T^{l}(x,y)^{-1} \cdot T^{t}(x,y) \in Z(I^{\mathrm{gr}})$. 
		
		We also observe from Eq~\eqref{eqn:1023} the useful equation
		\begin{align}
			\begin{split}
				a^{t(x)^{-1}} = t(x^{-1}) T^{t}(x,x^{-1})^{-1} a T^{t}(x,x^{-1})^{-1} t(x^{-1})^{-1} &= \sigma^{t}_{x^{-1}} \circ i_{T^{t}(x,x^{-1})^{-1}}(a) \\
				&= \sigma^{l}_{x^{-1}} \circ i_{T^{l}(x,x^{-1})^{-1}}(a) \\
				&= a^{l(x)^{-1}} .
			\end{split}
		\end{align}
		Now we observe that
		\begin{align}\label{eqn:2345}
			\begin{split}
				\{ a,T^{t}(x,x^{-1})\} = \{ a,t(x)t(x^{-1})\} &= \{ a,t(x)\} \cdot \{ a,t(x^{-1})\}^{t(x)} \\
				&= \nu_{x}^{t}(a) \cdot \sigma_{x}^{t}(\nu_{x^{-1}}^{t}(a)) \\
				&= \nu_{x}^{l}(a) \cdot \sigma_{x}^{l}(\nu_{x^{-1}}^{l}(a)) \\
				&= \{ a,T^{l}(x,x^{-1})\} ;
			\end{split}
		\end{align}
		consequently, $\{ T^{t}(x,x^{-1}),a\} = \{ T^{l}(x,x^{-1}),a\}$. It then follows from Eq~\eqref{eqn:2345} and Eq~\eqref{eqn:1023} that
		\begin{align*}
			\{ T^{t}(x,x^{-1})^{-1},a\} = \{ a,T^{t}(x,x^{-1})\}^{T^{t}(x,x^{-1})^{-1}} &= i_{T^{t}(x,x^{-1})^{-1}} \left( \{ a,T^{t}(x,x^{-1})\} \right) \\
			&= i_{T^{l}(x,x^{-1})^{-1}} \left(\{ a,T^{l}(x,x^{-1})\} \right) \\
			&= \{ T^{l}(x,x^{-1})^{-1},a\};
		\end{align*}
		consequently, $\{ T^{t}(x,x^{-1})^{-1},a\} = \{ T^{l}(x,x^{-1})^{-1},a\}$. Then be expanding the product we have
		\begin{align}\label{eqn:2000}
			\begin{split}
				\{ a,T^{t}(x,y)\} &= \{ a,t(x)t(y)t(xy)^{-1}\} \\ 
				&= \{ a,t(x)\} \cdot \{ a,t(y)\}^{t(x)} \cdot \{ a,t(xy)^{-1}\}^{t(x)t(y)} \\
				&= \{ a,t(x)\} \cdot \{ a,t(y)\}^{t(x)} \cdot \{ a,t(xy)^{-1}\}^{t(x)t(y)} \\
				&= \{ a,t(x)\} \cdot \{ a,t(y)\}^{t(x)} \cdot \left(\{ a,t((xy)^{-1})\} \cdot \{ a,T^{t}(xy,(xy)^{-1})^{-1}\}^{t((xy)^{-1})} \right)^{t(x)t(y)}  \\
				&= \nu^{t}_{x} \cdot \sigma^{t}_{x}(\nu^{t}_{y}) \cdot \sigma^{t}_{x} \left( \sigma^{t}_{y} \left( \nu^{t}_{(xy)^{-1}}(a) \cdot \sigma^{t}_{(xy)^{-1}} \big(\{ a,T^{t}(xy,(xy)^{-1})^{-1}\} \big) \right) \right) \\
				&= \nu^{l}_{x} \cdot \sigma^{l}_{x}(\nu^{l}_{y}) \cdot \sigma^{l}_{x} \left( \sigma^{l}_{y} \left( \nu^{l}_{(xy)^{-1}}(a) \cdot \sigma^{l}_{(xy)^{-1}} \big(\{ a,T^{l}(xy,(xy)^{-1})^{-1}\} \big) \right) \right) \\
				&= \{ a,T^{l}(x,y)\}.
			\end{split}
		\end{align}
Then we can use Eq~\eqref{eqn:2000} with Eq~\eqref{(2.13)} to conlcude $\{ a,S(x,y)\} = \{ a,T^{l}(x,y)^{-1} \cdot T^{t}(x,y)\} =1$ and so $\{ S(x,y),a\}=1$; thus, $S(x,y) \in \zeta(I)$.

		The task now is to show $S_{f}(x,y) \in \zeta(I)$. Let us first observe that
		\begin{align}
			\begin{split}
				\{ t(x)^{-1},a\} = \{ t(x^{-1}) \cdot T^{t}(x,x^{-1})^{-1}, a\} &= \{ T^{t}(x,x^{-1})^{-1}, a\}^{t(x^{-1})} \cdot \{ t(x^{-1}),a\} \\
				&= \sigma^{t}_{x^{-1}} \left(\{ T^{t}(x,x^{-1})^{-1}, a\} \right) \cdot \tau^{t}_{x^{-1}}(a) \\
				&= \sigma^{l}_{x^{-1}} \left(\{ T^{l}(x,x^{-1})^{-1}, a\} \right) \cdot \tau^{l}_{x^{-1}}(a) \\
				&= \{ l(x)^{-1},a\}.
			\end{split}
		\end{align}
		Then we can consider the calculation
		\begin{align*}
			\{ t(x),t(y)\} \cdot a \cdot \{ t(x),t(y)\}^{-1} &= [\{t(x),t(y)\}, a ] \cdot a \\
			&= \left\{ [t(x),t(y)],a \right\} \cdot a \\
			&= \{ t(y)t(x)^{-1}t(y)^{-1},a\}^{t(x)} \cdot \{ t(x),a\} \cdot a \\
			&= \{ t(x)^{-1}t(y)^{-1},a\}^{t(x)t(y)} \cdot \{ t(y),a\}^{t(x)} \cdot \{ t(x),a\} \cdot a \\
			&= \{ t(y)^{-1},a\}^{t(x)t(y)t(x)^{-1}} \cdot \{ t(x)^{-1}, a\}^{t(x)t(y)} \cdot \{ t(y),a\}^{t(x)} \cdot \{ t(x),a\} \cdot a \\
			&= \sigma^{t}_{x} \circ \sigma^{t}_{y} \left(\{ t(y)^{-1},a\}^{t(x)^{-1}} \right) \cdot \sigma^{t}_{x} \circ \sigma^{t}_{y} \left(\{ t(x)^{-1},a\} \right) \cdot \sigma^{t}_{x} \left( \tau^{t}_{y}(a) \right) \cdot \tau^{t}_{x}(a) \cdot a \\
			&= \sigma^{l}_{x} \circ \sigma^{l}_{y} \left(\{ l(y)^{-1},a\}^{l(x)^{-1}} \right) \cdot \sigma^{l}_{x} \circ \sigma^{l}_{y} \left(\{ l(x)^{-1},a\} \right) \cdot \sigma^{l}_{x} \left( \tau^{l}_{y}(a) \right) \cdot \tau^{l}_{x}(a) \cdot a \\
			&= \{ l(x),l(y) \} \cdot a \cdot \{ l(x),l(y)\}^{-1} .
		\end{align*}
		Then we can use this to see that
		\begin{align*}
			i_{T^{t}_{f}(x,y)}(a) &= \{ t(x),t(y)\} \cdot t(\{ x,y\})^{-1} \cdot a \cdot t(\{ x,y\}) \cdot \{ t(x),t(y)\}^{-1} \\ 
			&= \{ t(x),t(y)\} \cdot \sigma^{t}_{\{ x,y\}^{-1}} \left( a^{T^{t}(\{ x,y\},\{ x,y\}^{-1})^{-1}} \right) \cdot \{ t(x),t(y)\}^{-1} \\
			&= \{ l(x),l(y)\} \cdot \sigma^{l}_{\{ x,y\}^{-1}} \left( a^{T^{l}(\{ x,y\},\{ x,y\}^{-1})^{-1}} \right) \cdot \{ l(x),l(y)\}^{-1} \\
			&= i_{T^{l}_{f}(x,y)}(a) ;
		\end{align*}
		therefore, $S_{f}(x,y) = T^{l}_{f}(x,y)^{-1} \cdot T^{t}_{f}(x,y) \in Z(I^{\mathrm{gr}})$. Putting the previous observations together we have
		\begin{align*}
			\{ a,\{ t(x),t(x)\}\} &= \{ (a^{t(y)^{-1}})^{t(y)},\{ t(x),t(x)\}\} \\
			&= \left\{ \{ a^{t(y)^{-1}} , t(x) \}, t(y)^{t(x)} \right\} \cdot \left\{ \{ t(y) , a^{t(y)^{-1}} \}, t(x)^{a^{t(y)^{-1}}} \right\} \\
			&= \sigma^{t}_{x} \left( \nu^{t}_{y} \left( \nu^{t}_{x} \big( a^{t(y)^{-1}} \big)^{t(x)^{-1}} \right) \right) \cdot \nu^{t}_{x} \left( \tau^{t}_{x} \left( a^{t(y)^{-1}} \right)^{(a^{-1})^{t(y)^{-1}}} \right)^{a^{t(y)^{-1}}} \\
			&= \sigma^{l}_{x} \left( \nu^{l}_{y} \left( \nu^{l}_{x} \big( a^{l(y)^{-1}} \big)^{l(x)^{-1}} \right) \right) \cdot \nu^{l}_{x} \left( \tau^{l}_{x} \left( a^{l(y)^{-1}} \right)^{(a^{-1})^{l(y)^{-1}}} \right)^{a^{l(y)^{-1}}}  \\
			&= \{ a,\{ l(x),l(x)\}\} 
		\end{align*}
		Then we finally observe that
		\begin{align*}
			\{ a,S^{t}_{f}(x,y)\} = \{ a,\{ t(x),t(x)\} \cdot t(\{ x,y\})^{-1}\} &= \{ a,\{ t(x),t(x)\}\} \cdot \{ a,t(\{ x,y\})^{-1}\}^{\{ t(x),t(x)\}} \\
			&= \{ a,\{ l(x),l(x)\}\} \cdot \{ a,l(\{ x,y\})^{-1}]^{\{ l(x),l(x)\}} \\
			&= \{ a,\{ l(x),l(x)\} \cdot l(\{ x,y\})^{-1}\} = \{ a,S^{l}_{f}(x,y)\} \\
		\end{align*}
		which implies $\{ a, S_{f}(x,y)\} = \{a, S^{l}_{f}(x,y)^{-1} \cdot S^{l}_{f}(x,y)\} = 1$. This yields $S_{f}(x,y) \in \zeta(I)$; altogether, we have shown that $S$ is a 2-cocycle of datum $(Q,\zeta(I))$. Then $[S] \ast [A] = [I \rtimes_{T \cdot S} Q] = [A']$.
	\end{proof}

\vspace{0.2cm}

\section{Kernel-Preserving Automorphisms of Extensions}\label{section:4}

In this section, we prove for an arbitrary extension of multiplicative Lie algebras an analogue of Wells's exact sequence from group theory \cite{wells}. This is a true generalization since the original theorem is recovered by restriction to the subvariety $\mathcal M\mathcal L^{\mathrm{gr}}$.

\begin{theorem}\label{thm:wells}
Assume the extension $1 \rightarrow I \longrightarrow A \stackrel{\pi}{\longrightarrow} Q \rightarrow 1$ of multiplicative Lie algebras realizes the action terms $\chi$. There is an exact sequence 
\begin{align*}
1 \longrightarrow \mathrm{Der}(Q,\zeta(I),\chi_{\zeta}) \longrightarrow \Aut_{I} A \longrightarrow C(Q,I,\chi) \stackrel{W_{T}}{\longrightarrow} H^{2}(Q,\zeta(I),\chi_{\zeta}).
\end{align*}
\end{theorem}

There is a group action $\Aut I \times \Aut Q \curvearrowright \mathrm{Act}(Q,I)$ given by $\chi^{(\omega,\kappa)} = ( \sigma^{(\omega,\kappa)}, \tau^{(\omega,\kappa)}, \nu^{(\omega,\kappa)} )$ where
\begin{align*}
\sigma^{(\omega,\kappa)}(x) &= \omega \circ \sigma( \kappa^{-1}(x) ) \circ \omega^{-1} \\
\tau^{(\omega,\kappa)}(x) &= \omega \circ \tau( \kappa^{-1}(x) ) \circ \omega^{-1} \\
\nu^{(\omega,\kappa)}(x) &= \omega \circ \nu( \kappa^{-1}(x) ) \circ \omega^{-1} .
\end{align*}
There is also a group action $\Aut I \times \Aut Q \curvearrowright \mathcal H^{2}_{\mathcal M\mathcal L}(Q,I)$ given by $T^{(\omega,\kappa)} = \left\{ T^{(\omega,\kappa)}(x,y) , T_{f}^{(\omega,\kappa)}(x,y) , \chi^{(\omega,\kappa)} \right\}$ where
\begin{align*}
T^{(\omega,\kappa)}(x,y) &:= \omega \circ T ( \kappa^{-1}(x),\kappa^{-1}(y) ) \\
T_{f}^{(\omega,\kappa)}(x,y) &:= \omega \circ T_{f} ( \kappa^{-1}(x),\kappa^{-1}(y) ) .
\end{align*}
It can be shown that $T^{(\omega,\kappa)}$ is compatible whenever $T$ is compatible.

The \emph{compatible automorphisms} are the set of pairs $C(Q,I,\chi) = \{ (\omega,\kappa) \in \Aut I \times \Aut Q : [\chi] = [\chi^{(\omega,\kappa)}] \} $. By the Correspondence Theorem~\ref{thm:corresp1}, if $[T^{(\omega,\kappa)}]$ represent an extension $\mathrm{Ext}_{[\chi]}(Q,I)$ then there is a unique element $W_{T}(\omega,\kappa) \in H^{2}(Q,\zeta(I),\chi_{\zeta})$ such that 
\begin{align*}
[I \rtimes_{T^{(\omega,\kappa)}} Q ] \cdot W_{T}(\omega,\kappa) = [I \rtimes_{T} Q ] .
\end{align*}
Given $(\omega,\kappa) \in C(Q,I,\chi)$, we have that $\chi^{(\omega,\kappa)} \sim \chi$ is witnessed by a map $h:Q \rightarrow I$. If $A \approx I \rtimes_{T} Q$ where $T$ is defined by a section $l : Q \rightarrow A$ and $\chi = \chi^{l}$, then set $t(x) = h(x) \cdot l(x)$ for $x \in Q$. Then $T^{t} \sim T^{l} = T$. But since $(\omega,\kappa) \in C(Q,I,\chi)$ we see that $T^{(\omega,\kappa)} = \left\{ T^{(\omega,\kappa)}(x,y) , T_{f}^{(\omega,\kappa)}(x,y) , \chi^{(\omega,\kappa)} \right\} = \left\{ T^{(\omega,\kappa)}(x,y) , T_{f}^{(\omega,\kappa)}(x,y), \chi^{t} \right\}$ and both $T^{t}$ and $T^{(\omega,\kappa)}$ have the same action terms. By the argument in the Correspondence theorem we see that $\im (T^{t} - T^{(\omega,\kappa)}) \subseteq \zeta(I)$ and $[ T - T^{(\omega,\kappa)} ] = [ T^{t} - T^{(\omega,\kappa)} ]$. It follows that the \emph{Wells map} $W_{T}: C(Q,I,\chi) \rightarrow H^{2}(Q,\zeta(I),\chi_{\zeta})$ can be defined by 
	\begin{align}
		W_{T}(\omega,\kappa) := [ T - T^{(\omega,\kappa)} ].
	\end{align}

	For $d \in \mathrm{Der}(Q,\zeta(I),\chi_{\zeta})$ define $\phi_{d}(a,x) = \left\langle a\cdot d(x) , x \right\rangle$. The proof that $\mathrm{Stab} A \approx \mathrm{Der}(Q, I,\chi)$ for an extension $A$ realizing datum $(Q,I,\chi)$ with $I$ abelian also proves $\phi_{d} \in \Aut_{I} A$ since $\im d \subseteq \zeta(I)$. Then there is a homomorphism $i: \mathrm{Der}(Q,\zeta(I),\chi_{\zeta}) \rightarrow \Aut_{I} A$ given by $i(d):= \phi_{d}$

	Fix the extension $1 \rightarrow I \longrightarrow A \stackrel{\pi}{\longrightarrow} Q \rightarrow 1$ of multiplicative Lie algebras which realizes the abstract kernel $[\chi]$. Fix a section $t: Q \rightarrow A$ of $\pi$. For $\phi \in \Aut_{I} A$, define $\phi_{t}(x) := \pi \circ \phi \circ t(x)$. Then define $\psi(\phi) := (\phi|_{I},\phi_{t})$. We must show $\phi_{t}$ does not depend on the choice of section $t$, $\psi$ is a homomorphism and $(\phi|_{I},\phi_{t}) \in C(Q,I,\chi)$.

	\begin{lemma}
		Let $\pi : A \to Q$ be an extension realizing the datum $(Q,I,\chi)$ and $t: Q \to M$ be a lifting of $\pi$. Then $\psi: Aut_IA \to C(Q,I,\chi)$ defined by $\psi(\phi)= (\phi_I,\phi_t)$ is a homomorphism.
	\end{lemma}
	\begin{proof}
		First we show that $\phi_t$ is independent of the choice of section $t$. Let $t,t':Q\to A$ be liftings of $\pi$. Then $\pi(t(x))=\pi(t'(x))=x$, $t(x)t'(x)^{-1}\in I$. Thus $\pi(\phi(t(x)t'(x)^{-1}))=1$. Hence $\phi_t(x)=\phi_{t'}(x)$ for all $x\in Q$. To show that $\psi$ is a homomorphism. Let $\phi,\phi' \in Aut_IA$ such that $\psi(\phi)=(\phi|_I,\phi_t)$ and $\psi(\phi)=(\phi'|_I,\phi'_t)$. Consider $\pi(t(\phi_t(x))\phi(t(x))^{-1})=\phi_t(x)\phi_t(x)^{-1}=1$, there exist $h: Q \to I$ such that $h(x)\phi(t(x))=t(\phi_t(x))$. Consider $\pi(\phi(\phi'(t(x))))=\pi(\phi(h'(x)^{-1}t(\phi'_t(x))))=\pi(\phi(t(\phi'_t(x))))=\pi(h(\phi'_t(x))^{-1}t(\phi_t(\phi'_t(x)))) = \phi_t\circ \phi'_t(x) $. Thus 
		\[
		\pi(\phi\circ \phi')= (\phi\circ \phi'|_I, \pi\circ \phi \circ \phi' \circ t)= (\phi\circ \phi'|_I, \phi_t\circ \phi'_t). 
		\]
		Now we prove that $(\phi|_I, \phi_t) \in C(Q, I , \chi)$. Consider
		\begin{align*}
			\sigma_{\phi_t(x)}(a) = t(\phi_t(x))\cdot a \cdot t(\phi_t(x))^{-1} = h(x)\cdot \phi(t(x)) a \phi(t(x))^{-1} \cdot h(x)^{-1} & = h(x) \cdot \phi(t(x) \phi^{-1}(a) t(x)^{-1}) \cdot h(x)^{-1}\\
			& = h(x)\cdot \phi|_I(\sigma_x(\phi_I^{-1}(a))) \cdot h(x)^{-1} \\
			& = h(x) \cdot \sigma_{\phi_t(x)}^{(\phi|_I,\phi_t)}(a) \cdot h(x)^{-1}
		\end{align*}

		\begin{align*}
			\tau_{\phi_t(x)}(a) = \{ t(\phi_t(x)),a\} = \{ h(x)\phi (t(x)),a\} & = \{ \phi(t(x)),a\}^{h(x)} \cdot \{ h(x),a\} \\
			& = \phi(\{ t(x), \phi^{-1}(a)\})^{h(x)}\cdot \{ h(x),a\} \\
			& = \phi|_I(\tau_x (\phi|_I^{-1}(a)))^{h(x)}\cdot \{ h(x),a\} \\
			& = \tau_{\phi_t(x)}^{(\phi|_I,\phi_t)}(a)^{h(x)} \cdot \{ h(x),a\}.
		\end{align*}
		Similarly, we can show that $\nu_{\phi_t(x)}(a) = \{ a,h(x)\}\cdot \nu_{\phi_t(x)}^{(\phi|_I,\phi_t)}(a)^{h(x)}$. Hence $(\phi|_I, \phi_t) \in C(Q, I , \chi)$.
	\end{proof}
	

\begin{proof}(proof of Theorem~\ref{thm:wells})
	We first show $i$ is injective. Assume $i(d) = \id \in \Aut_{I} A$. Then $\left\langle a, x \right\rangle = i(d)(a,x) = \left\langle a \cdot d(x), x \right\rangle$ implies $d(x) = 1$ for all $x \in Q$; thus, $d \equiv 0$.

	We now show $\im i = \ker \psi$. Take $d \in \mathrm{Der}(Q,\zeta(I),\chi_{\zeta})$. Then $\phi_{d}(a,1) = \left\langle a \cdot d(1) , 1 \right\rangle = \left\langle a,1 \right\rangle = \id_{I}$. We also see that for $x \in Q$, $\pi \circ \phi_{d} \circ t(x) = \pi \circ \phi_{d} \left( 1,x \right) = \pi (1 \cdot d(x), x ) = x$; altogether, $\im i \subseteq \ker \psi$.

	Now assume $(\phi|_{I},\phi_{t}) = \psi(\phi) = \id$. This yields $\left\langle a, 1 \right\rangle = \phi(a,1)$ and $x = \phi_{t}(x) = \pi \circ \phi \circ t(x) = \pi \circ \phi(1,x)$. This implies there is a function $d : Q \rightarrow I$ such that $\phi(1,x) = \left\langle d(x), x \right\rangle$. Then 
	\begin{align*}
		\phi(a,x) = \phi(a,1) \cdot \phi(1,x) = \left\langle a,1 \right\rangle \cdot \left\langle d(x), x \right\rangle = \left\langle a \cdot d(x) , x \right\rangle.
	\end{align*}
	Since $\phi \in \Aut_{I} A$, we have
	\begin{align*}
		\left\langle a \cdot \sigma_{x}(b) \cdot T(x,y) \cdot d(xy) , xy \right\rangle = \phi \left( \left\langle a, x \right\rangle \cdot \left\langle b, y \right\rangle \right) &= \phi(a,x) \cdot \phi(b,y) \\	
		&= \left\langle a \cdot d(x), x \right\rangle \cdot \left\langle b \cdot d(y), y \right\rangle \\
		&= \left\langle a \cdot d(x) \cdot \sigma_{x}(b) \cdot \sigma_{x}(d(y)) \cdot T(x,y) , xy \right\rangle
	\end{align*}
	which yields 
	\begin{align*}
		\sigma_{x}(b) \cdot T(x,y) \cdot d(xy) = d(x) \cdot \sigma_{x}(b) \cdot \sigma_{x}(d(y)) \cdot T(x,y);
	\end{align*}
	thus, $d$ is a derivation of the group reduct. By taking $a=1,y=1$ we see from the above that $\sigma_{x}(b) \cdot d(x) = d(x) \cdot \sigma_{x}(b)$ which implies $\im d \subseteq Z(I^{\mathrm{gr}})$ because $\sigma_{x}$ is a bijection of $I$. Similarly, since $\phi$ respects the bracket operation we have 
	\begin{align*}
		\tau_{x}(b) &\cdot \{ a,b\} \cdot ba \cdot T_{f}(x,y) \cdot \sigma_{\{ x,y\}}(a^{-1} \cdot \nu_{y}(a) \cdot b^{-1}) \cdot d(\{ x,y\}) \\
		&= \tau_{x}(b) \cdot \tau_{x}(d(y)) \cdot \{ a \cdot d(x),b \cdot d(y)\} \cdot bd(y) \cdot a d(x) \cdot T_{f}(x,y) \cdot \\
		&\quad \quad \quad \quad \sigma_{\{ x,y\}}(d(x)^{-1}a^{-1} \cdot \nu_{y}(a) \cdot \nu_{y}(d(x)) \cdot d(y)^{-1}b^{-1} ) .
	\end{align*}
	Again, if we evaluate $a=1,y=1$ we conclude from the above that $1 = \{ d(x),b\} \cdot b \cdot d(x) \cdot b^{-1} d(x)^{-1} = \{ d(x),b\}$ since we already saw that $\im d \subseteq Z(I^{\mathrm{gr}})$. We see that $\im d \subseteq \zeta(I)$. We have established $d \in \mathrm{Der}(Q,\zeta(I),\chi_{\zeta})$. This shows $\phi \in \im i$.

	We now show $\im \psi = \ker W_{T}$. Take $\phi \in \Aut_{I} A$. Define $t'(x) := \phi \circ t \circ \phi_{t}^{-1}(x)$. Then $(a,t \circ \pi(a)) \in \alpha_{I}$ for all $a \in A$ implies $( t(x) , t'(x) ) = (\phi \circ \phi^{-1} \circ t(x) , \phi \circ t \circ \pi \circ \phi^{-1} \circ t(x) ) \in \alpha_{I}$; thus, $t$ and $t'$ are both section of $\pi : A \rightarrow Q$. Then $T=T^{t}$ and $T^{t'}$ are equivalent compatible 2-cocycles. We also see that
	\begin{align*}
		T^{(\phi|_{I},\phi_{t})}(x,y) = \phi \circ T(\phi_{t}^{-1}(x),\phi_{t}^{-1}(y)) &= \phi \left( t(\phi_{t}^{-1}(x)) \cdot t(\phi_{t}^{-1}(y)) \cdot t(\phi_{t}^{-1}(x) \cdot \phi_{t}^{-1}(y))^{-1} \right) \\
		&= \phi \circ t \circ \phi_{t}^{-1}(x) \cdot \phi \circ t \circ \phi_{t}^{-1}(y) \cdot \phi \circ t \circ \phi_{t}^{-1}(x \cdot y) \\
		&= t'(x) \cdot t'(y) \cdot t'(x \cdot y)^{-1} \\
		&= T^{t'}(x,y) .
	\end{align*}
	Similarly, we can show that $T_{f}^{(\phi|_{I},\phi_{t})}(x,y) = T_{f}^{t'}(x,y)$. We see that $[T] = [T^{(\phi|_{I},\phi_{t})}]$ and so $W_{T} \circ \psi(\phi) = W_{T}(\phi|_{I},\phi_{t}) = [T - T^{(\phi|_{I},\phi_{t})}] = 0$. This establishes that $\im \psi \subseteq \ker W_{T}$.

	Now assume $(\omega,\kappa) \in C(Q,I,\chi)$ such that $W_{T}(\omega,\kappa) = 0$, i.e. $[T-T^{(\omega,\kappa)}]=0$. Then there exist a map $h: Q \to \zeta(I)$ such that 
	$\omega(T(\kappa^{-1}(x),\kappa^{-1}(y)))\cdot h(x\cdot y) = h(x)\cdot \sigma_x(h(y)) \cdot T(x,y) $
	and $\omega(T_f(\kappa^{-1}(x),\kappa^{-1}(y)))\cdot h([x,y])= \tau_x(h(y)) \cdot h(y) \cdot h(x)\cdot T_f(x,y)\cdot \sigma_{[x,y]}(h(x)^{-1}\cdot \nu_y(h(x))\cdot h(y)^{-1}) $. 
	Define $\phi : A \to A $ by $\phi(a,x)=\langle
	\omega(a)h(x), \kappa(x)\rangle$. It is easy to see that $\phi|_I=\omega$ and $\kappa= \pi \circ \phi \circ t$, so $\phi|_t=\kappa$. Thus $\psi(\phi)=(\omega, \kappa)$. Now we prove that $\phi \in \Aut_IA$.
	\begin{align*}
		\phi(\langle a, x \rangle \cdot \langle b , y \rangle ) &= \phi (\langle a \cdot \sigma_x(b) \cdot T(x,y) , x \cdot y \rangle) \\
		&= \langle \omega(a) \cdot \omega(\sigma_x(b)) \cdot \omega (T(x,y)) h (x \cdot y), \kappa(x \cdot y) \rangle \\
		&= \langle \omega(a) \cdot \omega(\sigma_x(b)) \cdot h(x) \cdot \sigma_x(h(y)) \cdot T(\kappa(x),\kappa(y)), \kappa(x) \cdot \kappa(y) \rangle \\
		&= \langle \omega(a) \cdot h(x) \cdot \sigma_{\kappa(x)}\omega(b) \cdot \sigma_{\kappa(x)}(h(y)) \cdot T(\kappa(x),\kappa(y)), \kappa(x) \cdot \kappa(y) \rangle\\
		&= \langle \omega(a) h(x) , \kappa(x) \rangle \cdot \langle \omega(b) h(y) , \kappa(y) \rangle = \phi(\langle a, x \rangle ) \cdot \phi(\langle b, y \rangle ).
	\end{align*}
	Also, 
	\begin{align*}
		\phi(\{ \langle a,x \rangle , \langle b, y \rangle \}) &= \phi(\langle \tau_x(b)^a \{ a,b\} \cdot ba \cdot T_f(x,y) \cdot \sigma_{\{ x,y\}}(a^{-1}\cdot \nu_y(a) \cdot b^{-1}),\{ x,y\} \rangle)\\
		&= \langle \omega(\tau_x(b)^a)\cdot \omega(\{ a,b\} \cdot ba) \cdot \omega (T_f(x,y)) \omega(\sigma_{\{ x,y\}}(a^{-1} \cdot \nu_y(a) \cdot b^{-1}))\cdot h\{ x,y\}, \kappa(\{ x,y\}) \rangle \\
		&= \langle (\tau_{\kappa(x)}(\omega(b))^{\omega(a)})\cdot \omega(\{ a,b\}\cdot ba) \cdot \tau_{\kappa(x)}(h(y))\cdot h(y) \cdot h(x) \cdot T_f(\kappa(x),\kappa(y))\cdot  \\
		& \phantom{=}~~ \sigma_{\{ \kappa(x),\kappa(y)\}}(h(x)^{-1} \cdot \nu_{\kappa(y)}(h(x))\cdot h(y)^{-1}) \cdot 
		\sigma_{\{ \kappa(x),\kappa(y)\}}(\omega (a^{-1}\cdot \nu_y(a)\cdot b^{-1})) , \{ \kappa(x),\kappa(y)\} \rangle \\
		&= \langle \tau_{\kappa(x)}(\omega(b)h(y))^{\omega(a)h(x)}\{ \omega(a)h(x),\omega(b)h(y)\} \cdot \omega(b)h(y) \cdot \omega(a)h(x)\cdot T_f(\kappa(x),\kappa(y)) \cdot \\
		& \phantom{=}~~ \sigma_{\{ \kappa(x),\kappa(y)\}}(h(x)^{-1}\omega(a)^{-1} \cdot \nu_{\kappa(y)}(\omega(a)h(x))\cdot h(y)^{-1}\omega(b)^{-1}), \{ \kappa(x),\kappa(y)\} \rangle \\
		&= \{ \langle \omega(a) h(x), \kappa(x) \rangle , \langle \omega(b) h(y), \kappa(y) \rangle\}= \{ \phi\langle a, x \rangle, \phi \langle b , y \rangle \}
	\end{align*}
	Thus $\ker W_{T} \subseteq \im \psi$.

\end{proof}

\begin{proposition}
	Let $A=A_{0} \triangleright A_{1} \triangleright \cdots \triangleright A_{n} = 1$ be an ideal series of $A \in \mathcal M\mathcal L$. Then the stabilizer of the ideal series is a nilpotent group of class n-1.
\end{proposition}
\begin{proof}
	Since each ideal is a normal subgroup satisfying an additional property, any stabilizer of the ideal series is also a stabilizer of the corresponding normal series in the group reduct $A^{\mathrm{gr}}$. This means the stabilizer of the ideal series in $A$ is a subgroup of the stabilizer of the series in $A^{\mathrm{gr}}$. The result follows from the group result. 
\end{proof}

\vspace{0.2cm}

\section{Hochschild-Serre sequence}\label{section:5}

In this section we present and prove for multiplicative Lie algebras our version of the cohomological inflation/deflation exact sequence for groups \cite[Thm 2]{hochserre}.

\begin{theorem}\label{thm:hoschserre} 
Let $1 \rightarrow I \longrightarrow M \stackrel{\pi}{\longrightarrow} Q \rightarrow 1$ be an exact sequence of multiplicative Lie algebras and $\chi = \{\sigma, \tau,\nu\}$ action terms for datum $(M,A)$ with $A$ abelian which is realized by some extension. If $H^{0}(M,A,\chi)=0$, then we have an exact sequence
	\begin{align*}
		1 \rightarrow H^{1}(Q,A^{I},\hat{\chi}) \longrightarrow H^{1}(M,A,\chi) \longrightarrow H^{1}(I,A,\chi_{I})^{\square} \stackrel{\partial}{\longrightarrow} H^{2}(Q,A^{I},\hat{\chi}) \longrightarrow H^{2}(M,A,\chi).
	\end{align*}
\end{theorem}

We begin by defining the $1^{\mathrm{st}}$ and zeroth-cohomology groups. The following definitions privilege the ``conjugation'' action in the action terms over the ``bracket'' actions which are relegated to secondary importance.


\begin{lemma}
Let $\chi$ be compatible action terms for datum $(Q,I)$ with $I$ abelian. For each $a \in I$, $d_{a}(x) := a - \sigma_{x}(a)$ is a multiplicative Lie algebra derivation of the datum. 
\end{lemma}
\begin{proof}
Note we have a correspondence between derivations and stabilizing automorphism of the semidirect product $I \rtimes_{\chi} Q$ given by $\phi(a,x) = \left\langle a + d_{\phi}(x), x \right\rangle$. Recall the principal group derivation can be implemented by $\left\langle a , x \right\rangle^{\left\langle b , 1 \right\rangle} = \left\langle a + d_{b}(x) , x \right\rangle$. The evaluations of the identities 
\begin{align*}
\Big( \left\langle 1 , x \right\rangle \cdot \left\langle 1 , y \right\rangle \Big)^{\left\langle b , 1 \right\rangle} &= \left\langle 1 , x \right\rangle^{\left\langle b , 1 \right\rangle} \cdot \left\langle 1 , y \right\rangle^{\left\langle b , 1 \right\rangle}& \Big\{ \left\langle 1 , x \right\rangle, \left\langle 1 , y \right\rangle \Big\}^{\left\langle b , 1 \right\rangle} &= \Big\{ \left\langle 1 , x \right\rangle^{\left\langle b , 1 \right\rangle} , \left\langle 1 , y \right\rangle^{\left\langle b , 1 \right\rangle} \Big\}
\end{align*}
in $I \rtimes_{\chi} Q$ show that $d_{b}(x)$ is a multiplicative lie algebra derivation of the datum $(Q,I,\chi)$.  
\end{proof}


A \emph{principal derivation} for datum $(Q,A,\chi)$ with $A$ abelian is a function of the form $d_{a}(x) = a - \sigma_{x}(a)$. The set of principal derivations $\mathrm{PDer}(Q,A,\chi)$ forms a subgroup of derivations.


\begin{definition}
Let $(Q,A,\chi)$ datum in $\mathcal M\mathcal L$ with $A$ abelian. The $1^{\mathrm{st}}$-\emph{cohomology} of the datum $(Q,A,\chi)$ is
\begin{align}
H^{1}(Q,A,\chi) := \mathrm{Der}(Q,A,\chi)/\mathrm{PDer}(Q,A,\chi) . 
\end{align}
The \emph{zeroth-cohomology} of datum $(Q,A,\chi)$ is
\begin{align}
H^{0}(Q,A,\chi) := \{a \in A: \sigma_{x}(a)=a, \forall x \in Q \}.
\end{align}
\end{definition}


Fix the extension $1 \rightarrow I \stackrel{i}{\longrightarrow} M \stackrel{\pi}{\longrightarrow} Q \rightarrow 1$. We are given action terms $\chi = (\sigma,\tau,\nu)$ for $(M,A)$ in which $A$ is abelian. We can represent the extension $M \approx I \rtimes_{T} Q$ for a 2-cocycle $T = \{ T(x,y), T_{f}(x,y) , \sigma^{s}, \tau^{s}, \nu_{s} \}$ defined by a section $s: Q \rightarrow M$ with $s(1)=1$.

Define $A^{I} = \{a \in A : \sigma_{x}(a) = a, \tau_{x}(a) = \nu_{x}(a) = 1, x \in I \}$ and note it is a subalgebra of $A$ since $A$ has a trivial bracket. There are induced action terms $\hat{\chi} = \{ \hat{\sigma}, \hat{\tau}, \hat{\nu} \}$ for datum $(Q,A^{I})$ given by 
\begin{align}
\hat{\sigma}_{q}(a) &:= \sigma_{m}(a)& \hat{\tau}_{q}(a) := \tau_{m}(a)& &\hat{\nu}_{q}(a) := \nu_{m}(a)
\end{align}
for $a \in A^{I}$ and any $m \in M$ such that $\pi(m)=q$. Since $a \in A^{I}$, the action terms are well-defined. By restriction to $I \triangleleft M$, we also have action terms for $\chi_{I} = \{ \sigma^{I}, \tau^{I}, \nu^{I}, \}$ for $(I,A)$.

The \emph{restriction} map on cohomology $\check{r} : H^{1}(M,A,\chi) \rightarrow H^{1}(I,A,\chi_{I})$ is defined by $\check{r}([d]):=[d|_{I}]$. This is well-defined since the restriction of principal derivation on $M$ using coefficients in $A$ is again a principal derivation on $I$. It is easy to see that restriction is a homomorphism on cohomology.

Given a 2-cocycle $S=\{ S(x,y), S_{f}(x,y) \}$ appropriate for $(Q,A^{I})$, then $S \circ \pi = \{ (S\circ \pi)(x,y), (S\circ \pi)_{f}(x,y) \}$ is a 2-cocycle appropriate for $(M,A)$ where 
\begin{align*}
(S \circ \pi)(x,y) &:= S(\pi(x),\pi(y))  &\quad \text{ and }& \quad &(S\circ \pi)_{f}(x,y) := S_{f}(\pi(x),\pi(y)) .
\end{align*} 
The \emph{inflation} maps $\check{\sigma} : H^{1}(Q,A^{I},\hat{\chi}) \rightarrow H^{1}(M,A,\chi)$ and $\check{\sigma} : H^{2}(Q,A^{I},\hat{\chi}) \rightarrow H^{2}(M,A,\chi)$ are defined by pre-composition 
\begin{align}
\hat{\sigma}([d]) &:= [d \circ \pi]  &\quad \text{ and }& \quad &\hat{\sigma}([S]) := [S \circ \pi] 
\end{align} 
for $[d] \in H^{1}(Q,A^{I},\hat{\chi})$ and $[S] \in H^{2}(Q,A^{I},\hat{\chi})$. It is straightforward to see that the inflation of a principal derivation is again a principal derivation and that if $G$ is a 2-coboundary for $(Q,A^{I})$ witnessed by $h:Q \rightarrow A^{I}$, then $G \circ \pi$ is a 2-coboundary of $(M,A)$ witnessed by $h \circ \pi : M \rightarrow A$. We also see that pre-composition respects the point-wise group operation induced in the codomain so that inflation is a group homomorphism. It remains to observe that $S \circ \pi$ is compatible whenever $S$ is compatible which follows from the fact that for any term $t(\vec{x})$ and evaluation $\epsilon(\vec{x}) \in M^{\ar t}$ we have $t^{S \circ \pi}(\epsilon(\vec{x})) = t^{S}(\pi \circ \epsilon(\vec{x}))$.

The following is a useful characterization of the inflation map $\check{\sigma} : H^{2}(Q, A^{I},\hat{\chi}) \rightarrow H^{2}(M,A,\chi)$ in terms of the pull-back construction for multiplicative Lie algebras.


\begin{lemma}\label{lem:201}
	Let $\pi : M \rightarrow Q$ be an extension realizing datum $(Q,I,\chi)$ with $A$ an abelian multiplicative Lie algebra. Assume $\rho: E \rightarrow Q$ is an extension represented by the 2-cocycle $T$ where $E \approx A^{I} \rtimes_{\hat{\chi},T} Q$. Then the inflation along $\pi$ is given by $\check{\sigma}([T]) = [S]$ where $P \approx A \rtimes_{\chi, S} M$ and $P$
	is the pull-back of $\pi$ and $\rho$.
\end{lemma}
\begin{proof}
	We can write $M \approx I \rtimes_{T^{s}} Q$ for 2-cocycle $T^{s}$. Define $S := T \circ \pi$. Define a map $\phi: A \rtimes_{S} M \rightarrow P$ by $\phi(a,(b,x)) := \left\langle (a,x), (b,x) \right\rangle \in P \subseteq E \times M$. Bijectivity is immediate. We verify $\phi$ is a homomorphism. Observe that
	\begin{align*}
		\phi \left( \left\langle a , (b,x) \right\rangle \cdot \left\langle c , (d,y) \right\rangle \right) &= \phi \left( a + \sigma_{(b,x)}(c) + S((b,x),(d,y)) , (b,x) \cdot (d,y) \right) \\
		&= \phi \left( a + \hat{\sigma}_{x}(c) + T(x,y) , ( b \cdot \sigma^{s}_{x}(d) \cdot T^{s}(x,y) , x \cdot y ) \right) \\
		&= \left\langle ( a + \hat{\sigma}_{x}(c) + T(x,y) , x \cdot y) , ( b \cdot \sigma^{s}_{x}(d) \cdot T^{s}(x,y) , x \cdot y) \right\rangle \\
		&= \left\langle (a,x) \cdot (c,y) , (b,x) \cdot (d,y) \right\rangle \\
		&= \left\langle (a,x) , (b,x) \right\rangle \cdot \left\langle (c,y) , (d,y) \right\rangle \\
		&= \phi(a,(b,x)) \cdot \phi(c,(d,y)) .
	\end{align*}
	For the bracket operation we see that
	\begin{align*}
		\phi &\left(\Big\{ \left\langle a, (b,x)\right\rangle, \left\langle c, (d,y)\right\rangle \Big\} \right) \\
		&= \phi \left( a + c + \tau_{(b,x)}(c) + S_{f}((b,x),(d,y)) + \sigma_{[(b,x),(d,y)]}( \nu_{(d,y)}(a)-a-c), \{(b,x),(d,y)\} \right) \\
		&= \phi\left(a+c+\hat{\tau}_{x}(c)+T_{f}(x,y)+\hat{\sigma}_{\{x,y\}}(\hat{\nu}_{y}(a)-a-c), (b + d + \tau^{s}_{x}(d) + T^{s}_{f}(x,y) + \sigma_{\{x, y\}}( \nu^{s}_{y}(b) - b - d), \{x, y\}) \right) \\
		&= \left\langle(a + c + \hat{\tau}_{x}(c) + T_{f}(x,y) + \hat{\sigma}_{\{x,y\}}(\hat{\nu}_{y}(a) - a - c), \{x,y\}), (b + d + \tau^{s}_{x}(d) + T^{s}_{f}(x,y) + \sigma_{\{x,y\}}(\nu^{s}_{y}(b) - b - d),\{x,y\})\right\rangle \\
		&= \left\langle \big\{(a,x),(c,y) \big\}, \big\{(b,x), (d,y) \big\} \right\rangle \\
		&= \big\{\left\langle (a,x), (b,x) \right\rangle, \left\langle (c,y), (d,y) \right\rangle \big\} \\
		&= \big\{\phi(a,(b,x)), \phi(c,(d,y)) \big\}.
	\end{align*}
	So $\phi$ witnesses that $A \rtimes_{S} M$ is isomorphic to the pullback of $\pi$ and $\rho$.
\end{proof}


The corresponding pushout construction for multiplicative Lie algebras yields a characterization of the map $\rho^{\ast} : H^{2}(Q, A, \chi) \rightarrow H^{2}(Q, B,\gamma)$ induced by an injective homomorphism $\rho: A \rightarrow B$.


\begin{lemma}\label{lem:202}
	Let $A$ and $B$ be abelian multiplicative Lie algebras with $\chi$ action terms for $(Q,A)$ and $\gamma$ action terms for $(Q,B)$. Let $\rho: A \rightarrow B$ be a homomorphism which respects the action terms $\chi$ and $\gamma$. If $1 \rightarrow A \stackrel{j}{\rightarrow} E \stackrel{\pi}{\rightarrow} Q \rightarrow 1$ is an extension represented by the 2-cocycle $T$ and $\chi$ with $E \approx A \rtimes_{\chi,T} Q$, then $\rho^{\ast}([T]) = [R]$ where $P \approx B \rtimes_{\gamma,R} Q$ and $P$ is the pushout of $j$ and $\rho$.
\end{lemma}
\begin{proof}
Let $\chi$ be action terms for $(Q,A)$ induced by the extension $E$. There are action terms $\chi^{\ast}$ for datum $(E,B)$ given by $\chi^{\ast} = \pi \circ \chi$. Define $R := \rho \circ T$. The pushout of $j$ and $\rho$ is given by $P:= ( B \rtimes_{\chi^{\ast}} E )/ S$. We define a map $\phi : B \rtimes_{R} Q \rightarrow P$ by $\phi(a,x):= \left\langle a, (1,x) \right\rangle S$. To see that $\phi$ is a homomorphism we verify for the group operation that
	\begin{align*}
		\phi \left(\left\langle a, x \right\rangle \cdot \left\langle b, y \right\rangle \right) &= \phi \left( a + \sigma^{\gamma}_{x}(b) + R(x,y), x\cdot y \right) \\
		&= \phi \left( a + \sigma^{\gamma}_{x}(b) + \rho \circ T(x,y) , x \cdot y \right) \\
		&= \big< a + \sigma^{\gamma}_{x}(b) + \rho \circ T(x,y) , (1, x \cdot y ) \big> S \\
		&= \big< a + \sigma^{\gamma}_{x}(b) , (T(x,y), x \cdot y ) \big> S \\
		&= \big< a + \sigma^{\ast}_{(1,x)}(b) , (1,x) \cdot (1,y) \big> S \\
		&= \big< a, (1,x) \big> S \cdot \big< b, (1,y) \big> S \\
		&= \phi(a,x) \cdot \phi(b,y) . 
	\end{align*}
	For the bracket operation we see that 
	\begin{align*}
		\phi \left(\big\{\left\langle a, x \right\rangle , \left\langle b, y \right\rangle \big\}\right) &= \phi \left(a + b + \tau^{\gamma}_{x}(b)  + R_{f}(x,y) + \sigma^{\gamma}_{[x,y]}(\nu^{\gamma}_{y}(a) - a - b), \{x,y\} \right) \\
		&= \phi \left(a + b + \tau^{\gamma}_{x}(b) + \rho \circ T_{f}(x,y) + \sigma^{\gamma}_{\{x,y\}}(\nu^{\gamma}_{y}(a) - a - b ), \{x,y\} \right)\\
		&= \Big<a + b + \tau^{\gamma}_{x}(b) + \rho \circ T_{f}(x,y) + \sigma^{\gamma}_{\{x,y\}}(\nu^{\gamma}_{y}(a) - a - b ), (1,\{x, y\}) \Big> S\\
		&= \Big< a + b + \tau^{\gamma}_{x}(b) + \sigma^{\gamma}_{\{x,y\}}(\nu^{\gamma}_{y}(a) - a - b), \big(T_{f}(x,y), \{x, y\}\big) \Big> S\\
		&= \Big< a + b + \tau^{\ast}_{(1,x)}(b) + \sigma^{\ast}_{\{(1,x),(1,y)\}}(\nu^{\ast}_{(1,y)}(a) - a - b), \big\{ (1,x),(1,y) \big\} \Big> S\\
		&= \big\{ \left\langle a, (1,x) \right\rangle S , \left\langle b, (1,y) \right\rangle S \big\}\\
		&= \big\{ \phi(a,x) , \phi(b,y)\big\}. 
	\end{align*}
	Observe that $\left\langle a, (b,x) \right\rangle S = \left\langle a + \rho(b), (1,x) \right\rangle S = \phi(a + \rho(b),x)$. Also $\phi(a,x) = \left\langle a, (1,x) \right\rangle \in S$ implies $x = 1$ and $a= - \rho(c)$ and $c=1$. Then $a = 1$; altogether, $\phi$ is an isomorphism.
\end{proof}


The last step in preparation is to define the transgression homomorphism $\partial : H^{1}(I,A,\chi_{I})^{Q} \rightarrow H^{2}(Q,A^{I},\hat{\chi})$. The approach is directly inspired by the exposition of the transgression map in the case of groups given in Dekimpe, Hartl and Wauters \cite{dekimpe}. Consider the semidirect product 
\begin{align*}
1 \longrightarrow A \longrightarrow A \rtimes_{\chi} M \stackrel{\rho}{\longrightarrow} M \longrightarrow 1.
\end{align*}
For any map $d: I \rightarrow A$, the map $t_{d}(x) := \left\langle d(x), x \right\rangle$ determines a partial section $\rho \circ t_{d} = \id_{I}$ on $I$. We see that $d$ is a derivation if and only $t_{d}$ is a homomorphic section on $I$. If we start with a derivation $d:I \rightarrow A$, then we can pass to the normalizer $N_{A \rtimes_{\chi} M}(\im t_{d})$ of the subgroup $\im t_{d}$ and restrict to an extension 
\begin{align}\label{eq:transgressionidea}
1 \longrightarrow A \cap N_{A \rtimes_{\chi} M}(\im t_{d})  \longrightarrow N_{A \rtimes_{\chi} M}(\im t_{d})/\im t_{d} \longrightarrow \rho(N_{A \rtimes_{\chi} M}(\im t_{d}))/I \longrightarrow 1.
\end{align} 
We can then pair the 2-cocycle associated with the above extension to the derivation $d$ - this gives us a mapping from 1-cocycles to 2-cocycles.  We shall restrict ourselves to derivations such that $\rho(N_{A \rtimes_{\chi} M}(\im t_{d})) = M$ and $A^{I} = A \cap N_{A \rtimes_{\chi} M}(\im t_{d})$; in this manner, the associated paired 2-cocycle will represent a cohomology class in $H^{2}(Q,A^{I},\hat{\chi})$ which will be the image of the transgression map. Let us begin by giving a characterization of the normalizer associated to a section determined by a derivation $d : I \rightarrow A$.


\begin{lemma}\label{lem:1000}
Let $A \rtimes_{\chi,S} M$ be an extension and $I \triangleleft M$. Let $d: I \rightarrow A$ be a derivation of datum $(I,A,\chi_{I})$ and assume $H = \{ \left\langle d(n) , n \right\rangle : n \in I \}$ is a subalgebra. Then $\left\langle a, m \right\rangle \in N_{A \rtimes_{\chi,S} M}(H)$ if and only if for all $n \in I$,
	\begin{itemize}
		
		\item $\sigma_{m} \circ d(m^{-1} \cdot n \cdot m) - d(n) + S(m,m^{-1} \cdot n \cdot m) - S(n,m) = \sigma_{n}(a) - a$,
		
		\item $d(\{m,n\}) - \tau_{m} \circ d(n) + \sigma_{\{m,n\}} \circ d(n) - d(n) - S_{f}(m,n) = a + \sigma_{\{m,n\}} \circ \big( \nu_{n}(a) - a \big)$.
		
	\end{itemize}  
\end{lemma}   
\begin{proof}
Let $\left\langle a,m \right\rangle \in N_{A \rtimes M}(H)$. Then by Lemma~\ref{lem:normalizer} we have $\left\langle a,m \right\rangle H = H \left\langle a,m \right\rangle$ and $\left[ \left\langle a,m \right\rangle , H \right] \subseteq H$. Then for all $n \in I$ there exists $i,j \in I$ such that 
\begin{align}\label{eq:567}
\left\langle a,m \right\rangle \cdot \left\langle d(i), i \right\rangle  &= \left\langle d(n), n \right\rangle \cdot \left\langle a,m \right\rangle \\\label{eq:568}
\left\{ \left\langle a,m \right\rangle , \left\langle d(n), n \right\rangle \right\} &= \left\langle d(j), j \right\rangle .
\end{align} 
The left-hand and right-hand side of Eq~\eqref{eq:567} yields $a + \sigma_{m}(d(i)) + S(m,i) = d(n) + \sigma_{n}(a) + S(n,m)$ and $m \cdot i = n \cdot m$ which gives $\sigma_{m}(d(m^{-1} \cdot n \cdot m)) - d(n) + S(m,m^{-1} \cdot n \cdot m) - S(n,m) = \sigma_{n}(a) - a$. We see that Eq~\eqref{eq:568} yields $d(j) = \tau_{m}(d(n)) + a + d(n) + S_{f}(m,n) + \sigma_{\{ m,n \}}( \nu_{n}(a) - d(n) - a)$ and $\{m,n\} = j$. Then rewriting we have $d(\{ m, n \}) - \tau_{m}(d(n)) - d(n) - S_{f}(m,n) + \sigma_{\{ m,n \}}(d(n)) = a + \sigma_{\{ m,n \}}(\nu_{n}(a) - a)$. Reversing the calculation establishes the converse.
\end{proof}


The assumption that $H$ in Lemma~\ref{lem:1000} is a subalgebra occurs in two important cases: (1) when we have a semidirect product $A \rtimes M$, or (2) in an extension $A \rtimes_{\chi,S} M$ where the ideal $A \rtimes_{\chi,S} I$ is itself a semidirect product so that $[S|_{I}] = 0$. We will be interested in applying Lemma~\ref{lem:1000} in the case (2) where the factor sets will be zero.


\begin{definition}\label{def:square}
A derivation $d : I \rightarrow A$ satisfies the $\square$-condition if there exists $\eta : M \rightarrow A$ with $\eta(1)=1$ such that for all $m \in M$, $n \in I$,
\begin{itemize}
	
	\item[(S1)] $\sigma_{m} \circ d(m^{-1} \cdot a \cdot m) - d(a) = \sigma_{n} \circ \eta(m) - \eta(m)$, 
	
	\item[(S2)] $d(\{ m,n\}) - \tau_{m} \circ d(n) + \sigma_{\{ m,n\}} \circ d(n) - d(n) = \eta(m) + \sigma_{\{ m,n\}} \circ \big( \nu_{n}(\eta(m)) - \eta(m) \big)$.
	
\end{itemize}
\end{definition}


We then define $H^{1}(I,A,\chi_{I})^{\square} := \{ [d] \in H^{1}(I,A,\chi_{I}) : d \text{ satisfies the } \square-\text{condition} \}$.


\begin{lemma}
The restriction of derivations to those which satisfy the $\square$-condition respects cohomology. The image of the restriction homomorphism $\check{r}$ on $H^{1}(M,A,\chi)$ is contained in $H^{1}(I,A,\chi_{I})^{\square}$.
\end{lemma}
\begin{proof}
First statement is secured if we show the principal derivations all satisfy the $\square$-condition. Take $a \in A$ and set $d_{a}(x) = a - \sigma^{I}_{x}(a)$ for $x \in I$. Then for all $m \in M, n \in I$ we see that 
\begin{align*}
\sigma_{m} \circ d_{a}(m^{-1} \cdot n \cdot m) - d_{a}(n) &=  \sigma_{m} (a - \sigma^{I}_{m^{-1} \cdot n \cdot m}(a)) - (a - \sigma^{I}_{n}(a)) \\
&= \sigma_{m}(a) - \sigma_{n \cdot m}(a) - a  + \sigma^{I}_{n}(a) \\
&=  \sigma_{n}^{I}(a - \sigma_{m}(a)) - (a - \sigma_{m}(a)) .
\end{align*}
For the bracket operation observe that we can evaluate the identity Eq~\eqref{(2.5)} in $A \rtimes_{\chi} M$ to see that 
\begin{align*}
\left\langle a - \sigma^{I}_{\{ m,n \}}(a) , [m,n] \right\rangle  &= \left\{ \left\langle 1,m \right\rangle , \left\langle 1,n \right\rangle \right\}^{\left\langle a,1 \right\rangle} \\
&= \left\{ \left\langle 1,m \right\rangle^{\left\langle a,1 \right\rangle} , \left\langle 1,n \right\rangle^{\left\langle a,1 \right\rangle} \right\} \\
&= \left\{ \left\langle a - \sigma_{m}(a), m \right\rangle , \left\langle a - \sigma^{I}_{n}(a), n \right\rangle  \right\} \\
&= \Big< \tau_{m}(a - \sigma^{I}_{n}(a)) + a - \sigma_{m}(a) + a - \sigma^{I}_{n}(a)  \\
&\quad + \sigma^{I}_{\{ m,n \} } \big( \nu^{I}_{n}(a - \sigma_{m}(a)) - a + \sigma_{m}(a) - a + \sigma^{I}_{n}(a) \big) , [m,n] \Big> 
\end{align*}
which yields by the first coordinate
\begin{align*}
d_{a}(\{m,n\}) - \tau_{m}(d_{a}(n)) - d_{a}(n) + \sigma^{I}_{\{ m,n \} }(d_{a}(n)) = a - \sigma_{m}(a) + \sigma^{I}_{\{ m,n \} } \big( \nu^{I}_{n}(a - \sigma_{m}(a)) - a + \sigma_{m}(a) \big) .
\end{align*}
The conditions in Definition~\ref{def:square} are met if we take $\eta(m) := a - \sigma_{m}(a)$. In a similar manner, if we take a derivation $d : M \rightarrow A$, then the restriction $d|_{I}$ satisfies conditions (S1) and (S2) where we can take $\eta(m) = d(m)$. This shows the image of the restriction homomorphism is contained in $H^{1}(I,A,\chi_{I})^{\square}$. 
\end{proof}


We can now derive an explicit form for the transgression map $\partial : H^{1}(I,A,\chi^{I})^{\square} \rightarrow H^{2}(Q,A^{I},\hat{\chi})$. Fix $[d] \in H^{1}(I,A,\chi^{I})^{\square}$ and the extension  
\begin{align*}
1 \longrightarrow A \longrightarrow A \rtimes_{\chi} M \stackrel{\rho}{\longrightarrow} M \longrightarrow 1 .
\end{align*}
We set $t : I \rightarrow A \rtimes_{\chi} M$ given by $t(n) = \left\langle d(n) , n\right\rangle$ and write $\im t = H = \{ \left\langle d(n) , a \right\rangle : n \in I \}$ is a subalgebra of $A \rtimes_{\chi} M$. Since $d$ satisfies (S1) and (S2), we see that $\rho(N_{A \rtimes_{\chi} M}(H)) = M$ and $A \cap N_{A \rtimes_{\chi} M}(H) = A^{I}$ by the characterization in Lemma~\ref{lem:1000}. Then the restriction to the normalizer yields the extension
\begin{align}\label{eqn:super}
1 \longrightarrow A^{I} \longrightarrow  N_{A \rtimes_{\chi} M}(H) \stackrel{\rho}{\longrightarrow} M \longrightarrow 1 .
\end{align}
We can then fix a section $\tilde{t}: M \rightarrow N_{A \rtimes_{\chi} M}(H)$ given by $\tilde{t}(m) := \left\langle \tilde{\eta}(m), m \right\rangle$ where $\tilde{\eta}(m)$ is given by Definition~\ref{def:square}. Note we can choose $\tilde{\eta}$ such that $\tilde{\eta}(n) = d(n)$ for $n \in I$; consequently, $\tilde{t}(n) = \left\langle d(n) , n \right\rangle$ for all $n \in I$. Then passing to the quotient we derive the extension 
\begin{align}\label{eqn:millions}
 1 \longrightarrow A^{I} \longrightarrow N_{A \rtimes_{\chi} M}(H)/H \longrightarrow Q \longrightarrow 1 
\end{align}
with a section $\bar{t}: Q \rightarrow N_{A \rtimes_{\chi} M}(H)/H$ given by $\bar{t}(x) := \tilde{t}(l(x)) H = \left\langle \eta(x) , l(x) \right\rangle H$ where we have written $\eta(x) := \tilde{\eta}(l(x))$. Then a 2-cocyle is calculated by
\begin{align*}
	\bar{t}(x) \cdot \bar{t}(y) &= \tilde{t}(x) \cdot \tilde{t}(y) H \\
	&= T^{\tilde{t}}(l(x),l(y)) \cdot \tilde{t}(l(x) \cdot l(y)) H \\
	&= T^{\tilde{t}}(l(x),l(y)) \cdot \tilde{t}( T^{s}(x,y) \cdot l(x \cdot y)) H \\
	&= T^{\tilde{t}}(l(x),l(y)) \cdot T^{t}(T^{s}(x,y),l(x \cdot y))^{-1} \cdot \tilde{t}(T^{s}(x,y)) \cdot \tilde{t}(l(x \cdot y)) H \\
	&= T^{\tilde{t}}(l(x),l(y)) \cdot T^{t}(T^{s}(x,y),l(x \cdot y))^{-1} \cdot \tilde{t}(l(x \cdot y)) H 
\end{align*}
since $\tilde{t}|_{I} = t$; thus, the group part of the 2-cocycle is given by 
\begin{align*}
	T^{\bar{t}}(x,y) &= T^{\tilde{t}}(l(x),l(y)) \cdot T^{t}(T^{s}(x,y),l(x \cdot y))^{-1} \\
	&= \tilde{t}(l(x)) \cdot \tilde{t}(l(y)) \cdot \Big( \tilde{t}(T^{s}(x,y)) \cdot \tilde{t}(l(x \cdot y)) \Big)^{-1} \\
	&= \left\langle \eta(x) , l(x) \right\rangle \cdot \left\langle \eta(y) , l(y) \right\rangle \cdot \Big( \left\langle d(T^{s}(x,y)), T^{s}(x,y) \right\rangle \cdot \left\langle \eta( x \cdot y) , l( x \cdot y) \right\rangle \Big)^{-1} \\
	&= \left\langle \eta(x) + \sigma_{l(x)} \circ \eta(y), l(x) \cdot l(y) \right\rangle \cdot \left\langle d(T^{s}(x,y)) + \sigma_{T^{s}(x,y)} \circ \eta( x \cdot y) , T^{s}(x,y) l( x \cdot y ) \right\rangle^{-1} \\
	&= \left\langle \eta(x) + \sigma_{l(x)} \circ \eta(y), l(x) \cdot l(y) \right\rangle \cdot \left\langle - \sigma_{(l(x) \cdot l(y))^{-1}} \big(d(T^{s}(x,y)) + \sigma_{T^{s}(x,y)} \circ \eta( x \cdot y) \big) , (l(x) \cdot l(y))^{-1} \right\rangle \\
	&= \left\langle \eta(x) + \sigma_{l(x)} \circ \eta(y) - d(T^{s}(x,y)) - \sigma_{T^{s}(x,y)} \circ \eta(xy) , 1 \right\rangle .
\end{align*}
For the Lie bracket we can calculate
\begin{align*}
	\big\{ \bar{t}(x), \bar{s}(y)\big\} &= \big\{ \tilde{t}(x), \tilde{t}(y) \big\} H \\
	&= T^{\tilde{t}}_{f}(l(x),l(y)) \cdot \tilde{t}(\{ l(x), l(y)\})H \\
	&= T^{\tilde{t}}_{f}(l(x),l(y)) \cdot \tilde{t}( T^{s}_{f}(x,y) \cdot l(\{ x, y\})H \\
	&= T^{\tilde{t}}_{f}(l(x),l(y)) \cdot T^{t}(T^{s}_{f}(x,y),l(\{ x, y\}))^{-1} \cdot \tilde{t}(T^{s}(x,y)) \cdot \tilde{t}(l(\{ x , y\})) H \\
	&= T^{\tilde{t}}_{f}(l(x),l(y)) \cdot T^{t}(T^{s}_{f}(x,y),l(\{ x , y\}))^{-1} \cdot \tilde{t}(l(\{ x, y\})H .
\end{align*}
Then 
\begin{align*}
	T^{\bar{t}}(x,y) &= T^{\tilde{t}}_{f}(l(x),l(y)) \cdot T^{t}(T^{s}_{f}(x,y),l(\{ x, y\}))^{-1} \\
	&= \big\{\tilde{t}(l(x)),\tilde{t}(l(y)) \big\} \cdot \Big(\tilde{t}(T^{s}(x,y)) \cdot \tilde{t}(l(\{ x, y\})) \Big)^{-1} \\
	&= \big\{ \left\langle \eta(x), l(x)\right\rangle, \left\langle \eta(y) , l(y) \right\rangle \big\} \cdot \Big( \left\langle d(T^{s}_{f}(x,y)), T^{s}_{f}(x,y)\right\rangle \cdot \left\langle \eta(\{x, y\}), l(\{ x, y\}) \right\rangle \Big)^{-1} \\
	&= \left\langle \eta(x) + \eta(y) + \tau_{l(x)} \circ \eta(y) + \sigma_{\{ l(x),l(y)\}} \circ \big( \nu_{l(y)} \circ \eta(x) - \eta(x) - \eta(y) \big), \{ l(x), l(y)\} \right\rangle \\
	&\quad \cdot \left\langle d(T^{s}_{f}(x,y)) + \sigma_{T^{s}_{f}(x,y)} \circ \eta(\{ x, y\}),T^{s}_{f}(x,y)l(\{ x, y\}) \right\rangle^{-1} \\
	&= \left\langle \eta(x) + \eta(y) + \tau_{l(x)} \circ \eta(y) + \sigma_{\{ l(x),l(y)\}} \circ \big(\nu_{l(y)} \circ \eta(x) - \eta(x) - \eta(y) \big),\{ l(x),l(y)\} \right\rangle \\
	&\quad \cdot \left\langle - \sigma_{\{ l(x), l(y)\}^{-1}} \big(d(T_f^{s}(x,y)) + \sigma_{T_f^{s}(x,y)} \circ \eta(\{ x, y\} )\big), \{l(x), l(y)\}^{-1} \right\rangle \\
	&= \big<\ \eta(x) + \eta(y) + \tau_{l(x)} \circ \eta(y) + \sigma_{\{l(x),l(y)\}} \circ \big( \nu_{l(y)} \circ \eta(x) - \eta(x) - \eta(y) \big)- d(T^{s}_{f}(x,y)) \\
	&\quad - \sigma_{T^{s}_{f}(x,y)} \circ \eta(\{ x,y\}) , 1 \ \big> .
\end{align*}
The transgression is then determined by  
\begin{align}\label{eqn:trans5}
	\partial([d])(x,y) &:= \Big[ \eta(x) + \sigma_{l(x)} \circ \eta(y) - \sigma_{T^{s}(x,y)} \circ \eta(xy) - d(T^{s}(x,y)) \Big] \\
\begin{split}	\label{eqn:trans6}
	\partial([d])_{f}(x,y) &:= \Big[ \eta(x) + \eta(y) + \tau_{l(x)} \circ \eta(y) + \sigma_{\{ l(x),l(y) \}} \circ \big( \nu_{l(y)} \circ \eta(x) - \eta(x) - \eta(y) \big) \\
	&\quad - \sigma_{T^{s}_{f}(x,y)} \circ \eta(\{ x,y\}) - d(T^{s}_{f}(x,y)) \Big]
\end{split}
\end{align}
for a choice of $\tilde{\eta}: M \rightarrow A$ such that $\tilde{\eta} (n) = d(n)$ for $n \in I$. To see that the transgression is well-defined, note that for a different choice of $\tilde{\eta}' : M \rightarrow A$ which extends $d$, both $\tilde{\eta}$ and $\tilde{\eta}'$ define different respective sections $\tilde{t}$ and $\tilde{t}'$ of the extension in Eq~\eqref{eqn:super}; in turn, they induce different sections $\bar{t}$ and $\bar{t}'$ of the extension in Eq~\eqref{eqn:millions}. This implies the constructed 2-cocycles will be equivalent and determine the same cohomology class. To see that transgression is a homomorphism note that if $d$ and $d'$ are derivations which satisfy the $\square$-condition witnessed by $\eta$ and $\eta'$, respectively, then $d + d'$ satisfies the $\square$-condition witnessed by $\eta + \eta'$. The result now follows from the explicit form of the transgression in Eq~\eqref{eqn:trans5} and Eq~\eqref{eqn:trans6} since the action terms for abelian datum are homomorphisms of the group operation.

The following lemma gives a useful description of the three exact sequences used to construct the transgression map. The description is in terms of the pull-back and push-out constructions and the proof is a straightforward application of the mapping properties which offers no surprises.


\begin{lemma}\label{lem:exactnesshelp}
In the constructed diagram below
$N_{A \rtimes_{\chi} M}(H)$ is the pull-back of the maps $\tilde{\rho}$ and $\pi$ and $A \rtimes_{\chi} M$ is the pushout of the inclusion maps $i : A^{I} \hookrightarrow N_{A \rtimes_{\chi} M}(H)$ and $j: A^{I} \hookrightarrow A$.
\end{lemma}


We are know ready to complete the proof of the theorem. The hypothesis $H^{1}(M,A,\chi) = 0$ asserts that if $\sigma_{m}(a)=a$ for all $m \in M$, then $a=1$.


\begin{proof}(of Theorem \ref{thm:hoschserre})

\underline{$\ker \check{\sigma} = 0$} : Assume $[d] \in H^{1}(Q,A^{I},\hat{\chi})$ such that $\check{\sigma}([d]) = 0$. Then $d \circ \pi: M \rightarrow A^{I}$ is a principle derivation for $(M,A,\chi)$. So there exist $a \in A$ such that $d(x) = d \circ \pi(b,x) = a - \sigma_{\left\langle b,x \right\rangle}(a)$ for $\left\langle b,x \right\rangle \in M = I \rtimes_{T} Q$. Then $a - \sigma_{\left\langle 1,x \right\rangle}(a) = d \circ \pi(1,x) = d(x) = d \circ \pi(b,x) = a - \sigma_{(b,x)}(a)$ implies $\sigma_{\left\langle 1,x \right\rangle}(a) = \sigma_{\left\langle b,x \right\rangle}(a)$. We also see that $1 = d(1) = d \circ \pi(b,1) = a - \sigma_{\left\langle b,1 \right\rangle}(a)$ implies $\sigma_{\left\langle b, 1 \right\rangle}(a)=a$ for all $b \in I$. Since $\chi$ is compatible, it is realized by the semidirect product $A \rtimes_{\chi} M$; that is, there is a section $t$ such that $\chi = \chi^{t}$. Then using the fact that $A$ is abelian, we can directly calculate for all $m,x \in M,a \in A$ that 
\begin{align*}
\tau_{m}(d_{a}(x)) = \{ t(m) , a - \sigma_{x}(a) \} = \big\{ t(m) , [ a, t(x) ] \big\} &= \big[ t(m) , \{ a, t(x) \} \big] \\
&= - \big[ \{ a, t(x) \} , t(m) \big] \\
&= - d_{\nu_{x}(a)}(m)
\end{align*}
Then we see that $1 = \tau_{\left\langle c,y \right\rangle}(d \circ \pi(b,1) ) = \tau_{\left\langle c,y \right\rangle} \big( d_{a}( \left\langle b , 1 \right\rangle ) \big) = \nu_{\left\langle b,1 \right\rangle}(a) - \sigma_{\left\langle c,y \right\rangle}(\nu_{\left\langle b,1 \right\rangle}(a))$ which implies $\sigma_{\left\langle c,y \right\rangle}(\nu_{\left\langle b,1 \right\rangle}(a)) = \nu_{\left\langle b,1 \right\rangle}(a)$ for all $c \in I$; similarly, $\sigma_{\left\langle c,y \right\rangle}(\tau_{\left\langle b,1 \right\rangle}(a)) = \tau_{\left\langle b,1 \right\rangle}(a)$ for all $\left\langle c,y \right\rangle \in M$. Then the hypothesis $H^{0}(M,A,\chi)=0$ implies $a \in A^{I}$. This shows $d(x) = d \circ \pi(b,x) = a - \sigma_{\left\langle b,x \right\rangle}(a) = a - \sigma_{\left\langle 1,x \right\rangle}(a) = a - \hat{\sigma}_{x}(a)$ with $a \in A^{I}$; thus, $d$ is a principle derivation for datum $(Q,A^{I})$.

	\underline{$\im \check{\sigma} \subseteq \ker \check{r}$} : For $[d] \in H^{1}(Q,A^{I},\hat{\chi})$ we see that $d \circ \pi(b,1) = d(1)=1$ which yields $(d \circ \pi)|_{I} = 0$ and so $\check{r} \circ \check{\sigma}([d]) = \check{r}([d \circ \pi]) = 0$.

	\underline{$\im \check{\sigma} = \ker \check{r}$} : Take $[d'] \in H^{1}(M,A,\chi)$ such that $\check{r}([d']) = 0$. There exists $a \in A$ such that $d'(b,1) = a - \sigma_{\left\langle b, 1\right\rangle}(a)$ for $b \in I$. By taking $d(x) := d'(x) - d_{a}(x)$ we have $[d] = [d']$ in $H^{1}(M,A,\chi)$ and $d(b,1) = 1$ for all $b \in I$. Note for $x \in Q, b \in I$ we have $d( \left\langle 1,x \right\rangle \cdot \left\langle b,1 \right\rangle ) = d(1,x) + \sigma_{\left\langle 1,x \right\rangle}(d(b,1))$. Then 
	\begin{align*}
		d(b,x) = d( \left\langle b,1 \right\rangle \cdot \left\langle 1,x \right\rangle ) = d(\left\langle 1,x\right\rangle \cdot \left\langle \sigma^{s}_{x^{-1}}(b),1\right\rangle ) = d(1,x).
	\end{align*} 
	Altogether, we have 
	\begin{align}\label{eqn:1002}
		d(1,x) = d(b,x) = d(\left\langle b,1 \right\rangle \cdot \left\langle 1,x \right\rangle) = d(b,1) + \sigma_{\left\langle b,1 \right\rangle}(d(1,x)) = \sigma_{\left\langle b,1 \right\rangle}(d(1,x)).
	\end{align}
	From Eq~\eqref{eqn:1002} we see that for all $m \in M$, $b \in I$,
	\begin{align*}
		\sigma_{m} \circ \nu_{b}(d(1,x)) - \nu_{b}(d(1,x)) = \tau_{m} \big( d(1,x) - \sigma_{b}(d(1,x)) \big) = 0
	\end{align*}
	which implies $ \nu_{b}(d(1,x)) = 1$; similarly, $\tau_{b}(d(1,x)) = 1$. This yields $d(1,x) \in A^{I}$. Set $\gamma(x) := d(1,x)$. Then 
	\begin{align*}
		\gamma(x \cdot y) = d(1,x \cdot y) = d(T^{s}(x,y), x \cdot y) &= d\left( \left\langle 1,x \right\rangle \cdot \left\langle 1,y \right\rangle \right) \\
		&= d(1,x) + \sigma_{\left\langle 1, x \right\rangle}(d(1,y)) \\
		&= \gamma(x) + \hat{\sigma}_{x}(\gamma(y))
	\end{align*}
	and 
	\begin{align*}
		\gamma(\{x,y\}) = d(1,\{x,y\}) &= d(T_{f}(x,y),\{x,y\}) \\
		&= d\Big(\big\{\left\langle 1,x \right\rangle,\left\langle 1,y \right\rangle \big\} \Big) \displaybreak[0]\\
		&= d(1,x) + d(1,y) + \tau_{\left\langle 1, x \right\rangle}(d(1,y)) + \sigma_{\left\{ \left\langle 1,x \right\rangle, \left\langle 1,y \right\rangle \right\}} \big(\nu_{\left\langle 1,y \right\rangle}(d(1,x)) - d(1,x) - d(1,y) \big) \displaybreak[0]\\
		&= \gamma(x) + \gamma(y) + \hat{\tau}_{x}(\gamma(y)) + \sigma_{\left\langle T_{f}(x,y), \{x,y\} \right\rangle} \big( \hat{\nu}_{y}(\gamma(x)) - \gamma(x) - \gamma(y) \big) \displaybreak[0]\\
		&= \gamma(x) + \gamma(y) + \hat{\tau}_{x}(\gamma(y)) + \sigma_{\left\langle 1, \{x,y\} \right\rangle} \big( \hat{\nu}_{y}(\gamma(x)) - \gamma(x) - \gamma(y) \big) \displaybreak[0]\\
		&= \gamma(x) + \gamma(y) + \hat{\tau}_{x}(\gamma(y)) + \hat{\sigma}_{\{x,y\}} \big(\hat{\nu}_{y}(\gamma(x)) - \gamma(x) - \gamma(y) \big) .
	\end{align*}
	So we see that $d \in \mathrm{Der}(Q,A^{I},\hat{\chi})$. Then $\gamma \circ \pi(b,x) = \gamma(x) = d(1,x) = d(b,x)$ which gives $\hat{\sigma}([\gamma]) = [d]$.

	\underline{$\im \check{r} \subseteq \ker \partial$} : Take $[d] \in H^{1}(M,A,\chi)$. Then we observed from the development of the transgression that we can take $\eta : Q \rightarrow A$ to be $\eta(x):= d(1,x)$. Then 
	\begin{align*}
		\partial \circ \check{r}([d])(x,y) &= \partial ([d|_{I}])(x,y) \displaybreak[0]\\
		&= \Big\{ d(1,x) + \sigma_{\left\langle 1,x \right\rangle}(d(1,y)) - \sigma_{\left\langle T^{s}(x,y), 1 \right\rangle}(d(1,xy)) - d(T^{s}(x,y),1)\Big\} \displaybreak[0]\\
		&= \Big\{ d(\left\langle 1, x \right\rangle \cdot \left\langle 1, y \right\rangle ) - d(\left\langle T^{s}(x,y), 1 \right\rangle \cdot \left\langle 1, xy \right\rangle ) \Big\} \displaybreak[0]\\
		&= \Big\{ d(T^{s}(x,y), xy) - d( \left\langle T^{s}(x,y), 1 \right\rangle \cdot \left\langle 1, xy \right\rangle ) \Big\} \displaybreak[0]\\
		&= [ 0 ]
	\end{align*}
	and
	\begin{align*}
		\partial \circ \check{r}([d])_{f}(x,y) &= \partial ([d|_{I}])_{f}(x,y) \displaybreak[0]\\
		&= \Big\{ d(1,x) + d(1,y) + \tau_{}(d(1,y)) + \sigma_{\{ \left\langle 1,x \right\rangle, \left\langle 1,y \right\rangle \}}\big( \nu_{\left\langle 1,y \right\rangle}(d(1,x)) - d(1,x) - d(1,y) \big) \displaybreak[0]\\
		&\quad - \sigma_{\left\langle T^{s}_{f}(x,y),1 \right\rangle}(d(1,\{x,y\})) - d(T^{s}_{f}(x,y),1) \Big\} \displaybreak[0]\\
		&= \Big\{ d\big( \big\{ \left\langle 1,x \right\rangle, \left\langle 1,y \right\rangle \big\} \big) - d \big( \left\langle T^{s}_{f}(x,y),1 \right\rangle \cdot \left\langle 1,\{x,y\} \right\rangle \big)\Big\} \displaybreak[0]\\
		&= \Big\{ d(T^{s}_{f}(x,y),\{x,y\}) - d \big( \left\langle T^{s}_{f}(x,y),1 \right\rangle \cdot \left\langle 1,\{x,y\} \right\rangle \big) \Big\} \displaybreak[0]\\
		&= [0] .
	\end{align*}
	\underline{$\im \check{r} = \ker \partial$} : Let $[d] \in \ker \partial$; that is, there exists $h : Q \rightarrow A^{I}$ such that 
	\begin{align*}
		\partial([d])(x,y) &= h(x) + \hat{\sigma}_{x}(h(y)) - h(x \cdot y) \\
		\partial([d])_{f}(x,y) &= h(x) + h(y) + \hat{\tau}_{x}(h(y)) + \hat{\sigma}_{\{x,y\}} \Big( \hat{\nu}_{y}(h(x)) - h(x) - h(y) \Big) - h\big(\{x,y\}\big)
	\end{align*}
	Since $[d] \in H^{1}(I,A,\chi_{A})^{Q}$, there exists $\eta : Q \rightarrow A$ such that 
	\begin{align*}
		\sigma_{\left\langle 1, x \right\rangle} \circ d( \left\langle 1,x \right\rangle^{-1} \cdot \left\langle a,1 \right\rangle \cdot \left\langle 1,x \right\rangle ) - d(a) &= \sigma_{\left\langle a , 1 \right\rangle}(\eta(x)) - \eta(x) \\
		d(\tau^{s}_{x}(a)) - \tau_{\left\langle 1, x \right\rangle} \circ d(a) + \sigma_{\{\left\langle 1, x \right\rangle,\left\langle a, 1 \right\rangle\}} \circ d(a) - d(a) &= \eta(x) + \sigma_{\{\left\langle 1, x \right\rangle , \left\langle a, 1 \right\rangle\}} \circ \big( \nu_{\left\langle a, 1 \right\rangle}(\eta(x)) - \eta(x) \big) 
	\end{align*}
	for all $a \in I, x \in Q$; in addition, $\sigma_{\left\langle a , 1 \right\rangle}(\eta(x)) - \eta(x) = \sigma_{\left\langle 1 , x \right\rangle}(\eta(x^{-1})) - \sigma_{\left\langle a , x \right\rangle}(\eta(x^{-1}))$. Since $\im \partial \subseteq A^{I}$, we see that 
	\begin{align*}
		h(x) &+ \hat{\sigma}_{x}(h(y)) - h( x \cdot y) \displaybreak[0]\\
		&= \partial([d])(x,y) \\
		&= \sigma_{\left\langle a\sigma^{s}_{x}(b) , 1 \right\rangle}(\partial([d])(x,y)) \displaybreak[0]\\
		&= \sigma_{\left\langle a\sigma^{s}_{x}(b) , 1 \right\rangle}(\eta(x)) + \sigma_{\left\langle a\sigma^{s}_{x}(b) , x \right\rangle}(\eta(y)) - \sigma_{\left\langle a\sigma^{s}_{x}(b)T^{s}(x,y), 1\right\rangle}(\eta( x \cdot y )) - \sigma_{\left\langle a\sigma^{s}_{x}(b) , 1 \right\rangle}(d(T^{s}(x,y)))
	\end{align*}
	Now define a map $\Psi: M \rightarrow A$ by $\Psi(a,x):= d(a) + \sigma_{\left\langle a, 1 \right\rangle}(\eta(x)) - h(x)$. It is easy to see that the restriction of $\Psi$ to $I$ is $d$. We show that $\Psi \in \mathrm{Der}(M,A,\chi)$. For the group operation, let us first observe that
	\begin{align*}
		\sigma_{\left\langle 1, x \right\rangle^{-1}} \big( d( \sigma^{s}_{x}(b) ) \big) - d(b) = \sigma_{\left\langle b, 1 \right\rangle} (\eta(x^{-1})) - \eta(x^{-1}) = \sigma_{\left\langle 1, x \right\rangle^{-1}}(\eta(x)) - \sigma_{\left\langle b, 1 \right\rangle \cdot \left\langle 1, x \right\rangle^{-1}}(\eta(x))
	\end{align*}
	implies
	\begin{align}\label{eqn:7000}
		\sigma_{\left\langle a, 1 \right\rangle} \big( \sigma^{s}_{x}(b) \big) - \sigma_{\left\langle a, x \right\rangle}(d(b)) = \sigma_{\left\langle a, 1 \right\rangle}(\eta(x)) - \sigma_{\left\langle a \cdot \sigma^{s}_{x}(b), 1 \right\rangle}(\eta(x)).
	\end{align}
	Then using Eq~\eqref{eqn:7000} we see that 
	\begin{align*}
		\Psi(a,x) &+ \sigma_{\left\langle a,x \right\rangle }\left( \Psi(b,y) \right) \displaybreak[0]\\
		&= d(a) + \sigma_{\left\langle a, 1 \right\rangle}(\eta(x)) - h(x) + \sigma_{\left\langle a,x \right\rangle}\left( d(b) + \sigma_{\left\langle b, 1 \right\rangle}(\eta(y)) - h(y) \right) \displaybreak[0]\\
		&= d(a) + \sigma_{\left\langle a,x \right\rangle}(d(b)) + \sigma_{\left\langle a, 1 \right\rangle}(\eta(x)) + \sigma_{\left\langle a,\sigma^{s}_{x}(b) , x \right\rangle}(\eta(y)) - h(x) - \sigma_{\left\langle a,x \right\rangle}(h(y)) \displaybreak[0]\\
		&= d(a) + \sigma_{\left\langle a, 1 \right\rangle} \big( d( \sigma^{s}_{x}(b) ) \big) - \sigma_{\left\langle a, 1 \right\rangle}(\eta(x)) + \sigma_{\left\langle a \cdot \sigma^{s}_{x}(b) , 1 \right\rangle} (\eta(x)) + \sigma_{\left\langle a, 1 \right\rangle}(\eta(x)) \displaybreak[0]\\
		&\quad + \sigma_{\left\langle a,\sigma^{s}_{x}(b) , x \right\rangle}(\eta(y)) - h(x) - \sigma_{\left\langle a,x \right\rangle}(h(y)) \displaybreak[0]\\
		&= d(a) + \sigma_{\left\langle a, 1 \right\rangle} \big( d( \sigma^{s}_{x}(b) ) \big) + \sigma_{\left\langle a \cdot \sigma^{s}_{x}(b) , 1 \right\rangle} (\eta(x)) + \sigma_{\left\langle a,\sigma^{s}_{x}(b) , x \right\rangle}(\eta(y)) \displaybreak[0]\\
		&\quad - h(x) - \sigma_{\left\langle a,x \right\rangle}(h(y))  \displaybreak[0]\\
		&= d(a) + \sigma_{\left\langle a, 1 \right\rangle} \big( d( \sigma^{s}_{x}(b) ) \big) + \sigma_{\left\langle a \cdot \sigma^{s}_{x}(b) \cdot T^{s}(x,y), 1\right\rangle}(\eta( x \cdot y )) + \sigma_{\left\langle a \cdot \sigma^{s}_{x}(b) , 1 \right\rangle}(d(T^{s}(x,y))) - h( x \cdot y ) \displaybreak[0]\\
		&= d(a \cdot \sigma^{s}_{x}(b) \cdot T^{s}(x,y)) + \sigma_{\left\langle a \cdot \sigma^{s}_{x}(b) \cdot T^{s}_{x,y} , 1 \right\rangle}(\eta( x \cdot y )) - h( x \cdot y) \displaybreak[0]\\
		&= \Psi \left( a \cdot \sigma^{s}_{x}(b) \cdot T^{s}_{x,y} , x \cdot y \right) \displaybreak[0]\\
		&= \Psi \left( \left\langle a,x \right\rangle \cdot\left\langle b,y \right\rangle \right) .
	\end{align*}
	We need to consider the bracket observation. We first derive an equation among the values of the derivation $d$ and the operation $\eta$ by evaluating a multiplicative lie algebra identity in the semi-direct product $A \rtimes_{\chi} M$. Calculating in the algebra $A \rtimes_{\chi} M$, we see that
		\begin{align*}
		\big \{&\left\langle \sigma_{a}(\eta(x)), a \cdot x \right\rangle , \left\langle \sigma_{b}(\eta(y)), b \cdot y \right\rangle \big\} \displaybreak[0]\\
		&= \big\{\left\langle 1 , a \right\rangle \cdot \left\langle \eta(x) , x \right\rangle , \left\langle 1 , b \right\rangle \cdot \left\langle \eta(y) , y \right\rangle \big\} \displaybreak[0]\\
		&= \big\{\left\langle \eta(x), x \right\rangle, \left\langle 1 , b \right\rangle \big\}^{\left\langle 1 , a \right\rangle} \cdot \big\{ \left\langle 1 , a \right\rangle , \left\langle 1 , b \right\rangle \big\} \cdot \big\{ \left\langle \eta(x) , x \right\rangle , \left\langle \eta(y) , y \right\rangle \big\}^{\left\langle 1 , ba \right\rangle} \cdot \big\{ \left\langle 1 , a \right\rangle , \left\langle \eta(y) , y \right\rangle \big\}^{\left\langle 1 , b \right\rangle} \displaybreak[0]\\
		&= \left\langle \eta(x) + \sigma_{\{x,b\}}(\nu_{b}(\eta(x)) - \eta(x) ) ,\{x,b\} \right\rangle^{\left\langle 1 , a \right\rangle} \cdot \left\langle 1 , \{a,b\} \right\rangle \displaybreak[0]\\
		&\quad \cdot \left\langle \eta(x) + \eta(y) + \tau_{x}(\eta(y)) + \sigma_{\{x,y\}}(\nu_{y}(\eta(x)) - \eta(x) - \eta(y)) , \{x,y\} \right\rangle^{\left\langle 1 , ba \right\rangle} \cdot \left\langle \eta(u) + \tau_{a}(\eta(y)) - \sigma_{\{a,y\}}(\eta(y)) , \{a,y\} \right\rangle^{\left\langle 1 , b \right\rangle} \displaybreak[0]\\
		&= \left\langle \sigma_{a}(\eta(x)) + \sigma_{a\{x,b\}}(\nu_{b}(\eta(x)) - \eta(x) ) , \{x,b\}^{a} \right\rangle \cdot \left\langle 1 , \{a,b\} \right\rangle \displaybreak[0]\\
		&\quad \cdot \left\langle \sigma_{ba} \big(\eta(x) + \eta(y) + \tau_{x}(\eta(y))\big) + \sigma_{ba\{x,y\}}(\nu_{y}(\eta(x)) - \eta(x) - \eta(y)) , \{x,y\}^{ba} \right\rangle \displaybreak[0]\\
		&\quad \cdot \left\langle \sigma_{b} \big( \eta(u) + \tau_{a}(\eta(y)) \big) - \sigma_{b\{a,y\}}(\eta(y)) , \{a,y\}^{b} \right\rangle \displaybreak[0]\\
		&= \Big< \sigma_{s}(\eta(x)) + \sigma_{a\{x,b\}} \big( \nu_{b}(\eta(x)) - \eta(x) \big) + \sigma_{\{x,b\}^{a}\{a,b\}ba} \big( \eta(x) + \eta(y) + \tau_{x}(\eta(y)) \big) + \sigma_{\{x,b\}^{a}\{a,b\}ba\{x,y\}} \big( \nu_{y}(\eta(x)) - \eta(x) - \eta(y) \big)  \displaybreak[0]\\
		&\quad + \sigma_{\{x,b\}^{a}\{a,b\}\{x,y\}^{ba}b} \big( \eta(y) + \tau_{a}(\eta(y)) \big) - \sigma_{\{x,b\}^{a}\{a,b\}\{x,y\}^{ba}b\{a,y\}} \big( \eta(y) \big) \ ,\ \{x,b\}^{a} \cdot \{a,b\} \cdot \{x,y\}^{ba} \cdot \{a,y\}^{b} \Big> 
	\end{align*}
	but evaluating directly we have 
	\begin{align*}
		\big\{ &\left\langle \sigma_{a}(\eta(x)), a \cdot x \right\rangle , \left\langle \sigma_{b}(\eta(y)), b \cdot y \right\rangle \big\} \\
		&= \left\langle \sigma_{a}(\eta(x)) + \sigma_{b}(\eta(y)) + \tau_{ax}(\sigma_{b}(\eta(y))) + \sigma_{\{ax,yb\}}\big(\nu_{by}(\sigma_{a}(\eta(x))) - \sigma_{a}(\eta(x)) - \sigma_{b}(\eta(y)) \big), \{ax,yb\} \right\rangle .
	\end{align*}
	This gives the identity
	\begin{align}
		\begin{split}
			\sigma_{a}(\eta(x)) & + \sigma_{b}(\eta(y)) + \tau_{ax}(\sigma_{b}(\eta(y))) + \sigma_{\{ax,yb\}}\big(\nu_{by}(\sigma_{a}(\eta(x))) - \sigma_{a}(\eta(x)) - \sigma_{b}(\eta(y)) \big) \\
			&= \sigma_{s}(\eta(x)) + \sigma_{a\{x,b\}} \big(\nu_{b}(\eta(x)) - \eta(x) \big) + \sigma_{\{x,b\}^{a}\{a,b\}ba} \big(\eta(x) + \eta(y) + \tau_{x}(\eta(y)) \big) \\
			&\quad + \sigma_{\{x,b\}^{a}\{a,b\}ba\{x,y\}} \big(\nu_{y}(\eta(x)) - \eta(x) - \eta(y) \big) + \sigma_{\{x,b\}^{a}\{a,b\}\{x,y\}^{ba}b} \big(\eta(y) + \tau_{a}(\eta(y)) \big) \\
			&\quad - \sigma_{\{x,b\}^{a}\{a,b\}\{x,y\}^{ba}b\{a,y\}} \big( \eta(y) \big)
		\end{split}
	\end{align}
	Similarly, if we evaluate the equation 
	\begin{align*}
		\big\{\left\langle d(a) , a \cdot x \right\rangle , \left\langle d(b) , b \cdot y \right\rangle\big\} = \big\{ \left\langle d(a) , a \cdot x \right\rangle , \left\langle d(b) , b \right\rangle \big\} \cdot \big\{ &\left\langle d(a) , a \cdot x \right\rangle , \left\langle 1 , y \right\rangle \big\}^{\left\langle d(b) , b \right\rangle}
	\end{align*}
	in $A \rtimes_{\chi} M$ in two ways we derive the identity
	\begin{align}
		\sigma_{\{x,b\}^{a}\{a,b\}} \big(\nu_{b}(d(a)) - d(a) + \sigma_{b}(d(a)) \big) = \sigma_{\{ax,by\}} \big(\nu_{b}(d(a)) - d(a) + \sigma_{b}(d(a)) \big) .
	\end{align}
	Evaluating
	\begin{align*}
		\big\{\left\langle d(a), a \cdot x \right\rangle, \left\langle d(b) , b \right\rangle \big\} = \big\{ \left\langle 1 , x \right\rangle , \left\langle d(b) , b \right\rangle \big\}^{\left\langle d(a), a \right\rangle} \cdot \big\{ &\left\langle d(a) , a \right\rangle , \left\langle d(b), b \right\rangle \big\}
	\end{align*}
	in the semidirect product $A \rtimes_{\chi} M$ in two ways yields the identity
	\begin{align}
		\sigma_{\{x,b\}^{a}} \big(\tau_{a}(d(b)) + d(b)) \big) = \tau_{a}(d(b)) + d(b).
	\end{align}
	We can also record the following evaluation of the derivation $d$:
	\begin{align}
		\begin{split}
			d(a^{-1} &\cdot \nu^{s}_{y}(a) \cdot b^{-1}) \\ 
			&= \sigma_{a^{-1}} \big(\eta(y) + \tau_{a}(\eta(y)) \big) + \sigma_{a^{-1}\{a,y\}} \big(\nu_{y}(d(a)) - d(a) - \eta(y) \big) - \sigma_{a^{-1} \nu^{s}_{y}(a) b^{-1}} \big(d(b)\big).
		\end{split}
	\end{align}
    Now we can use these identities to evaluate
	\begin{align*}
		\Psi(a,x) & + \Psi(b,y) + \tau_{\left\langle a, x \right\rangle}(\Psi(b,y)) + \sigma_{\{\left\langle a, x \right\rangle,\left\langle b, y \right\rangle\}}\Big(\nu_{\left\langle b, y \right\rangle}(\Psi(a,x)) - \Psi(a,x) - \Psi(b,y)\Big) \\
		&= d(a) + d(b) + \sigma_{\left\langle a, 1 \right\rangle}(\eta(x)) + \sigma_{\left\langle b, 1 \right\rangle}(\eta(y)) + \tau_{\left\langle a, x \right\rangle}(d(b)) + \tau_{\left\langle a, x \right\rangle}(\sigma_{\left\langle b, 1 \right\rangle}(\eta(y))) \displaybreak[0]\\
		&\quad + \sigma_{\{\left\langle a, x \right\rangle, \left\langle b, y \right\rangle\}} \Big(\nu_{\left\langle b, y \right\rangle}(d(a)) + \nu_{\left\langle b, y \right\rangle}(\sigma_{\left\langle a, 1 \right\rangle}(\eta(x))) - d(a) - d(b) - \sigma_{\left\langle a, 1 \right\rangle}(\eta(x)) - \sigma_{\left\langle b, 1 \right\rangle} (\eta(y)) \Big) \displaybreak[0]\\
		&\quad - h(x) - h(y) - \tau_{\left\langle 1, x \right\rangle}(h(y)) - \sigma_{} \big( \nu_{\left\langle 1, y \right\rangle}(h(x)) - h(x) - h(y) \big) \displaybreak[0]\\
		&= d(a) + d(b) + \sigma_{\left\langle a, 1 \right\rangle}(\eta(x)) + \sigma_{\left\langle b, 1 \right\rangle}(\eta(y)) + \tau_{\left\langle a, x \right\rangle}(d(b)) + \tau_{\left\langle a, x \right\rangle}(\sigma_{\left\langle b, 1 \right\rangle}(\eta(y))) \displaybreak[0]\\
		&\quad + \sigma_{\{\left\langle a, x \right\rangle, \left\langle b, y \right\rangle\}} \Big( \nu_{\left\langle b, y \right\rangle}(d(a)) - d(a) - d(b) \Big) + \sigma_{\{\left\langle a, x \right\rangle, \left\langle b, y \right\rangle\}} \Big(\nu_{\left\langle b, y \right\rangle}(\sigma_{\left\langle a, 1 \right\rangle}(\eta(x))) - \sigma_{\left\langle a, 1 \right\rangle}(\eta(x)) - \sigma_{\left\langle b, 1 \right\rangle}(\eta(y)) \Big) \displaybreak[0]\\
		&\quad - \sigma_{\left\langle \tau^{s}_{x}(b)^{a} \cdot \{a,b\} \cdot ba, 1 \right\rangle} \Big( \eta(x) - \eta(y) - \tau_{x}(\eta(y)) - \sigma_{\left\langle T^{s}_{f}(x,y), \{x,y\} \right\rangle} \big( \nu_{\left\langle 1, y \right\rangle}(\eta(x)) - \eta(x) - \eta(y) \big)\Big) \displaybreak[0]\\
		&\quad + \sigma_{\left\langle \tau^{s}_{x}(b)^{a} \cdot \{a,b\} \cdot ba \cdot T^{s}_{f}(x,y), 1 \right\rangle}(\eta(\{x,y\})) + \sigma_{\left\langle \tau^{s}_{x}(b)^{a} \cdot \{a,b\} \cdot ba , 1 \right\rangle}(d(T^{s}_{f}(x,y))) - h(\{x,y\}) \displaybreak[0]\\
		&= d(a) + d(b) + \tau_{\left\langle a, 1 \right\rangle}(d(b)) + \sigma_{\left\langle a, 1 \right\rangle} \big( \tau_{\left\langle 1, x \right\rangle}(d(b)) \big) + \sigma_{\{\left\langle a, x \right\rangle, \left\langle b, y \right\rangle\}} \Big( \nu_{\left\langle b, y \right\rangle}(d(a)) - d(a) - d(b) \Big) \displaybreak[0]\\
		&\quad + \sigma_{\left\langle a, 1 \right\rangle}(\eta(x)) + \sigma_{a\{x,b\}} \big( \nu_{b}(\eta(x)) - \eta(x) \big) + \sigma_{\{x,b\}^{a}\{a,b\}\{x,y\}^{ba}b} \big( \eta(y) + \tau_{a}(\eta(y)) \big) \displaybreak[0]\\
		&\quad - \sigma_{\{x,b\}^{a}\{a,b\}\{x,y\}^{ba}b\{a,y\}} \big( \eta(y) \big) + \sigma_{\left\langle \tau^{s}_{x}(b)^{a} \cdot \{a,b\} \cdot ba \cdot T^{s}_{f}(x,y) , 1 \right\rangle}(\eta(\{x,y\})) + \sigma_{\left\langle \tau^{s}_{x}(b)^{a} \cdot \{a,b\} \cdot ba , 1 \right\rangle}(d(T^{s}_{f}(x,y))) - h(\{x,y\}) \displaybreak[0]\\
		&= d \big( \tau^{s}_{x}(b)^{a} \big) + d(b) + \tau_{\left\langle a, 1 \right\rangle}(d(b)) + \sigma_{\left\langle a \cdot \tau^{s}_{x}(b) , 1 \right\rangle}(d(b)) + \sigma_{\left\langle \tau_{x}(b)^{a} , 1 \right\rangle}(d(a)) - \sigma_{\left\langle a, 1 \right\rangle}(d(b)) \displaybreak[0]\\
		&\quad + \sigma_{\{\left\langle a, x\right\rangle, \left\langle b, y \right\rangle\}} \Big(\sigma_{\left\langle b, 1 \right\rangle}(\nu_{\left\langle 1, y \right\rangle}(d(a))) + \nu_{\left\langle b, 1 \right\rangle}(d(a)) - d(a) - d(b) \Big)  \displaybreak[0]\\
		&\quad + \sigma_{\left\langle a, 1 \right\rangle}(\eta(x)) + \sigma_{a\{x,b\}} \big( \nu_{b}(\eta(x)) - \eta(x) \big) + \sigma_{\{x,b\}^{a}\{a,b\}\{x,y\}^{ba}b} \big( \eta(y) + \tau_{a}(\eta(y)) \big) \displaybreak[0]\\
		&\quad - \sigma_{\{x,b\}^{a}\{a,b\}\{x,y\}^{ba}b\{a,y\}} \big( \eta(y) \big) + \sigma_{\left\langle \tau^{s}_{x}(b)^{a} \cdot \{a,b\} \cdot ba \cdot T^{s}_{f}(x,y) , 1 \right\rangle}(\eta(\{x,y\})) + \sigma_{\left\langle \tau^{s}_{x}(b)^{a} \cdot \{a,b\} \cdot ba , 1 \right\rangle}(d(T^{s}_{f}(x,y))) - h(\{x,y\}) \displaybreak[0]\\
		&= d \big( \tau^{s}_{x}(b)^{a} \big) + \sigma_{\left\langle \tau^{s}_{x}(b)^{a} \cdot \{a,b\} \cdot ba \cdot T^{s}_{f}(x,y) , \{x,y\} \right\rangle} \big(d( a^{-1} \cdot \nu^{s}_{y}(a) \cdot b^{-1}) \big) \displaybreak[0]\\
		&\quad + d(b) + \tau_{\left\langle a, 1 \right\rangle}(d(b)) + \sigma_{\left\langle a \cdot \tau^{s}_{x}(b), 1 \right\rangle}(d(b)) + \sigma_{\left\langle \tau_{x}(b)^{a}, 1 \right\rangle}(d(a)) - \sigma_{\left\langle a, 1 \right\rangle}(d(b))  \displaybreak[0]\\
		&\quad + \sigma_{\{\left\langle a, x \right\rangle, \left\langle b, y \right\rangle\}} \Big( \nu_{\left\langle b, 1 \right\rangle}(d(a)) - d(a) \Big) + \sigma_{\left\langle \tau^{s}_{x}(b) \cdot \{a,b\} \cdot ba \cdot T^{s}_{f}(x,y), \{x,y\} \right\rangle \cdot \left\langle a^{-1} \cdot \nu^{s}_{y}(a) , 1 \right\rangle} \Big( d(a) \Big)  \displaybreak[0]\\
		&\quad + \sigma_{\left\langle \tau^{s}_{x}(b)^{a} \cdot \{a,b\} \cdot ba \cdot T^{s}_{f}(x,y) , 1 \right\rangle}(\eta(\{x,y\})) + \sigma_{\left\langle \tau^{s}_{x}(b)^{a} \cdot \{a,b\} \cdot ba , 1 \right\rangle} \big( d(T^{s}_{f}(x,y)) \big) - h(\{x,y\}) \displaybreak[0]\\
		&= d \big(\tau^{s}_{x}(b)^{a} \big) + \sigma_{\left\langle \tau^{s}_{x}(b)^{a} \cdot \{a,b\} \cdot ba \cdot T^{s}_{f}(x,y) , \{x,y\} \right\rangle} \big( d( a^{-1} \cdot \nu^{s}_{y}(a) \cdot b^{-1} ) \big) \displaybreak[0]\\
		&\quad + d(b) + \tau_{\left\langle a, 1 \right\rangle}(d(b)) + \sigma_{\left\langle a \cdot \tau^{s}_{x}(b), 1 \right\rangle}(d(b)) + \sigma_{\left\langle \tau_{x}(b)^{a}, 1 \right\rangle}(d(a)) - \sigma_{\left\langle a, 1 \right\rangle}(d(b)) \displaybreak[0]\\
		&\quad + \sigma_{\{\left\langle a, x \right\rangle, \left\langle b, y \right\rangle\}} \Big(\nu_{\left\langle b, 1\right\rangle}(d(a)) - d(a) + \sigma_{\left\langle b, 1\right\rangle}(d(a)) \Big)  \displaybreak[0]\\
		&\quad + \sigma_{\left\langle \tau^{s}_{x}(b)^{a} \cdot \{a,b\} \cdot ba \cdot T^{s}_{f}(x,y) , 1 \right\rangle}(\eta(\{x,y\})) + \sigma_{\left\langle \tau^{s}_{x}(b)^{a} \cdot \{a,b\} \cdot ba, 1 \right\rangle} \big( d(T^{s}_{f}(x,y)) \big) - h(\{x,y\}) \displaybreak[0]\\
		&= d \big(\tau^{s}_{x}(b)^{a} \big) + \sigma_{\left\langle \tau^{s}_{x}(b)^{a} \cdot \{a,b\} \cdot ba \cdot T^{s}_{f}(x,y), \{x,y\} \right\rangle} \big( d( a^{-1} \cdot \nu^{s}_{y}(a) \cdot b^{-1} ) \big) \displaybreak[0]\\
		&\quad + d(b) + \tau_{\left\langle a, 1 \right\rangle}(d(b)) + \sigma_{\left\langle a \cdot \tau^{s}_{x}(b), 1 \right\rangle}(d(b)) + \sigma_{\left\langle \tau_{x}(b)^{a}, 1 \right\rangle}(d(a)) - \sigma_{\left\langle a, 1 \right\rangle}(d(b)) \displaybreak[0]\\
		&\quad + \sigma_{\{\left\langle \{x,b\}^{a} \cdot \{a,b\}, 1 \right\rangle} \Big( \nu_{\left\langle b, 1 \right\rangle}(d(a)) - d(a) + \sigma_{\left\langle b , 1 \right\rangle}(d(a)) \Big) \displaybreak[0]\\
		&\quad + \sigma_{\left\langle \tau^{s}_{x}(b)^{a} \cdot \{a,b\} \cdot ba \cdot T^{s}_{f}(x,y) , 1 \right\rangle}(\eta(\{x,y\})) + \sigma_{\left\langle \tau^{s}_{x}(b)^{a} \cdot \{a,b\} \cdot ba , 1 \right\rangle} \big( d(T^{s}_{f}(x,y)) \big) - h(\{x,y\}) \displaybreak[0]\\
		=& d \big( \tau^{s}_{x}(b)^{a} \big) + \sigma_{\left\langle \tau^{s}_{x}(b)^{a} \cdot \{a,b\} \cdot ba \cdot T^{s}_{f}(x,y), \{x,y\} \right\rangle} \big(d(a^{-1} \cdot \nu^{s}_{y}(a) \cdot b^{-1})\big) \\
		&\quad + d(b) + \tau_{\left\langle a, 1 \right\rangle}(d(b)) + \sigma_{\left\langle a \cdot \tau^{s}_{x}(b), 1 \right\rangle}(d(b)) + \sigma_{\left\langle \tau_{x}(b)^{a}, 1 \right\rangle}(d(a)) - \sigma_{\left\langle a, 1 \right\rangle}(d(b)) \displaybreak[0]\\
		&\quad + \sigma_{\{\left\langle \{x,b\}^{a} \cdot \{a,b\}, 1 \right\rangle} \Big(\nu_{\left\langle b, 1 \right\rangle}(d(a)) - d(a) - d(b) \Big) + \sigma_{\left\langle \tau^{s}_{x}(b)^{a} \cdot \{a,b\}, 1 \right\rangle} \big(d(b) + \sigma_{\left\langle b, 1\right\rangle}(d(a)) \big)  \displaybreak[0]\\
		&\quad + \sigma_{\left\langle \tau^{s}_{x}(b)^{a} \cdot \{a,b\} \cdot ba \cdot T^{s}_{f}(x,y), 1 \right\rangle}(\eta(\{x,y\})) + \sigma_{\left\langle \tau^{s}_{x}(b)^{a} \cdot \{a,b\} \cdot ba , 1 \right\rangle} \big(d(T^{s}_{f}(x,y)) \big) - h(\{x,y\}) \displaybreak[0]\\
		=& d\big(\tau^{s}_{x}(b)^{a} \big) + \sigma_{\left\langle \tau^{s}_{x}(b)^{a} \cdot \{a,b\} \cdot ba \cdot T^{s}_{f}(x,y), \{x,y\} \right\rangle} \big(d(a^{-1} \cdot \nu^{s}_{y}(a) \cdot b^{-1}) \big)\displaybreak[0]\\
		&\quad + d(b) + \tau_{\left\langle a, 1 \right\rangle}(d(b)) + \sigma_{\left\langle a \cdot \tau^{s}_{x}(b), 1 \right\rangle}(d(b)) + \sigma_{\left\langle \tau_{x}(b)^{a}, 1 \right\rangle}(d(a)) - \sigma_{\left\langle a, 1 \right\rangle}(d(b)) \displaybreak[0]\\
		&\quad + \sigma_{[\left\langle \{x,b\}^{a}, 1 \right\rangle} \Big( d(\{a,b\}) - d(a) - d(b) - \tau_{\left\langle a, 1 \right\rangle}(d(b)) \Big) + \sigma_{\left\langle \tau^{s}_{x}(b)^{a} \cdot \{a,b\}, 1 \right\rangle} \big(d(b) + \sigma_{\left\langle b, 1\right\rangle}(d(a)) \big)  \displaybreak[0]\\
		&\quad + \sigma_{\left\langle \tau^{s}_{x}(b)^{a} \cdot \{a,b\} \cdot ba \cdot T^{s}_{f}(x,y), 1 \right\rangle}(\eta(\{x,y\})) + \sigma_{\left\langle \tau^{s}_{x}(b)^{a} \cdot \{a,b\} \cdot ba , 1 \right\rangle} \big(d(T^{s}_{f}(x,y)) \big) - h(\{x,y\}) \displaybreak[0]\\
		&= d \big(\tau^{s}_{x}(b)^{a} \big) + \sigma_{\left\langle \tau^{s}_{x}(b)^{a} \cdot \{a,b\} \cdot ba \cdot T^{s}_{f}(x,y), \{x,y\} \right\rangle} \big(d(a^{-1} \cdot \nu^{s}_{y}(a) \cdot b^{-1}) \big) \displaybreak[0]\\
		&\quad + \sigma_{\left\langle \{x,b\}^{a}, 1 \right\rangle} \big( d(\{a,b\}) \big) + \sigma_{\left\langle \{x,b\}^{a} \cdot \{a,b\}, 1 \right\rangle} \big( d(ba) \big) + d(b) + \tau_{\left\langle a, 1 \right\rangle}(d(b)) + \sigma_{\left\langle \{x,b\}^{a}, 1 \right\rangle} \big(\sigma_{\left\langle a, 1 \right\rangle}(d(b))) \big) \displaybreak[0]\\
		&\quad - \sigma_{\left\langle a, 1 \right\rangle}(d(b)) - \sigma_{\left\langle \{x,b\}^{a}, 1 \right\rangle} \big(d(b) + \tau_{\left\langle a, 1 \right\rangle}(d(b)) \big) \displaybreak[0]\\
		&\quad + \sigma_{\left\langle \tau^{s}_{x}(b)^{a} \cdot \{a,b\} \cdot ba \cdot T^{s}_{f}(x,y), 1 \right\rangle}(\eta(\{x,y\})) + \sigma_{\left\langle \tau^{s}_{x}(b)^{a} \cdot \{a,b\} \cdot ba, 1 \right\rangle} \big(d(T^{s}_{f}(x,y)) \big) - h(\{x,y\}) \displaybreak[0]\\
		&= d \big(\tau^{s}_{x}(b)^{a} \big) + \sigma_{\left\langle \tau^{s}_{x}(b)^{a} \cdot \{a,b\} \cdot ba \cdot T^{s}_{f}(x,y), \{x,y\} \right\rangle} \big(d(a^{-1} \cdot \nu^{s}_{y}(a) \cdot b^{-1}) \big) \displaybreak[0]\\
		&\quad + \sigma_{\left\langle \{x,b\}^{a}, 1 \right\rangle} \big( d(\{a,b\}) \big) + \sigma_{\left\langle \{x,b\}^{a} \cdot \{a,b\}, 1 \right\rangle} \big(d(ba) \big) + \sigma_{\left\langle \{x,b\}^{a}, 1 \right\rangle} \big(\sigma_{\left\langle a, 1 \right\rangle}(d(b)) \big) - \sigma_{\left\langle a, 1 \right\rangle}(d(b)) \displaybreak[0]\\
		&\quad + \sigma_{\left\langle \tau^{s}_{x}(b)^{a} \cdot \{a,b\} \cdot ba \cdot T^{s}_{f}(x,y) , 1 \right\rangle}(\eta(\{x,y\})) + \sigma_{\left\langle \tau^{s}_{x}(b)^{a} \cdot \{a,b\} \cdot ba , 1 \right\rangle} \big(d(T^{s}_{f}(x,y)) \big) - h(\{x,y\})  \displaybreak[0]\\
		&= d \big(\tau^{s}_{x}(b)^{a} \big) + \sigma_{\left\langle \tau^{s}_{x}(b)^{a} \cdot \{a,b\} \cdot ba \cdot T^{s}_{f}(x,y), \{x,y\} \right\rangle} \big(d(a^{-1} \cdot \nu^{s}_{y}(a) \cdot b^{-1}) \big) \displaybreak[0]\\
		&\quad + \sigma_{\left\langle \{x,b\}^{a}, 1 \right\rangle} \big( d(\{a,b\}) \big) + \sigma_{\left\langle \{x,b\}^{a} \cdot \{a,b\}, 1 \right\rangle} \big( d(ba) \big) \displaybreak[0]\\
		&\quad + \sigma_{\left\langle \tau^{s}_{x}(b)^{a} \cdot \{a,b\} \cdot ba \cdot T^{s}_{f}(x,y), 1 \right\rangle}(\eta(\{x,y\})) + \sigma_{\left\langle \tau^{s}_{x}(b)^{a} \cdot \{a,b\} \cdot ba, 1 \right\rangle} \big(d(T^{s}_{f}(x,y)) \big) - h(\{x,y\}) \displaybreak[0]
	\end{align*}

	But we also have
	\begin{align*}
		\Psi \left(\left\{ \left\langle a,x \right\rangle, \left\langle b,y \right\rangle \right\} \right) &= \Psi \left(\tau^{s}_{x}(b)^{a} \cdot \{a,b\} \cdot ba \cdot T^{s}_{f}(x,y) \cdot \sigma^{s}_{\{x,y\}}(a^{-1} \cdot \nu^{s}_{y}(a) \cdot b^{-1}), \{x,y\} \right) \displaybreak[0]\\
		&= d\big(\tau^{s}_{x}(b)^{a} \cdot \{a,b\} \cdot ba \cdot T^{s}_{f}(x,y) \cdot \sigma^{s}_{\{x,y\}}(a^{-1} \cdot \nu^{s}_{y}(a) \cdot b^{-1}) \big) \displaybreak[0]\\
		&\quad + \sigma_{\left\langle \tau^{s}_{x}(b)^{a} \cdot \{a,b\} \cdot ba \cdot T^{s}_{f}(x,y) \cdot \sigma^{s}_{\{x,y\}}(a^{-1} \cdot \nu^{s}_{y}(a) \cdot b^{-1}), 1 \right\rangle} (\eta(\{x,y\})) - h(\{x,y\}) \displaybreak[0]\\
		&= d(\tau^{s}_{x}(b)^{a}) + \sigma_{\left\langle \tau^{s}_{x}(b)^{a} , 1 \right\rangle} \big( d(\{a,b\}) \big) + \sigma_{\left\langle \tau^{s}_{x}(b)^{a} \cdot \{a,b\}, 1 \right\rangle} \big(d(ba) \big) + \sigma_{\left\langle \tau^{s}_{x}(b)^{a} \cdot \{a,b\} \cdot ba, 1 \right\rangle} \big( d(T^{s}_{f}(x,y)) \big) \displaybreak[0]\\
		&\quad + \sigma_{\left\langle \tau^{s}_{x}(b)^{a} \cdot \{a,b\} \cdot ba \cdot T^{s}_{f}(x,y), 1 \right\rangle} \big(d( \sigma^{s}_{\{x,y\}}(a^{-1} \cdot \nu^{s}_{y}(a) \cdot b^{-1})) \big) \displaybreak[0]\\
		&\quad + \sigma_{\left\langle \tau^{s}_{x}(b)^{a} \cdot \{a,b\} \cdot ba \cdot T^{s}_{f}(x,y) \cdot \sigma^{s}_{\{x,y\}}(a^{-1} \cdot \nu^{s}_{y}(a) \cdot b^{-1}), 1 \right\rangle} (\eta(\{x,y\})) - h(\{x,y\}) \displaybreak[0]\\
		&= d(\tau^{s}_{x}(b)^{a}) + \sigma_{\left\langle \tau^{s}_{x}(b)^{a} , 1 \right\rangle} \big(d(\{a,b\}) \big) + \sigma_{\left\langle \tau^{s}_{x}(b)^{a} \cdot \{a,b\}, 1 \right\rangle} \big( d(ba) \big) + \sigma_{\left\langle \tau^{s}_{x}(b)^{a} \cdot \{a,b\} \cdot ba, 1 \right\rangle} \big(d(T^{s}_{f}(x,y)) \big) \displaybreak[0]\\
		&\quad + \sigma_{\left\langle \tau^{s}_{x}(b)^{a} \cdot \{a,b\} \cdot ba \cdot T^{s}_{f}(x,y), \{x,y\} \right\rangle} \big(d(a^{-1} \cdot \nu^{s}_{y}(a) \cdot b^{-1})\big) \displaybreak[0]\\
		&\quad + \sigma_{\left\langle \tau^{s}_{x}(b)^{a} \cdot \{a,b\} \cdot ba \cdot T^{s}_{f}(x,y),1 \right\rangle} (\eta(\{x,y\})) - h(\{x,y\}). \displaybreak[0]\\
	\end{align*}
	Putting these last two derivations together we see that $\Psi$ is a derivation of the bracket operation.

\underline{$\im \partial \subseteq \ker \check{\sigma}$} : By Lemma~\ref{lem:exactnesshelp} we see that the extension $1 \rightarrow A^{I} \rightarrow N_{A \rtimes M}(H) \rightarrow M \rightarrow 1$ is the pullback of $\pi$ and $\tilde{\rho}$ in the extension $1 \rightarrow A^{I} \stackrel{i}{\rightarrow} N_{A \rtimes M}(H)/H \stackrel{\rho}{\rightarrow} Q \rightarrow 1$, and for the inclusion $A^{I} \stackrel{j}{\rightarrow} A$,we see that the extension $1 \rightarrow A \stackrel{i}{\rightarrow} A \rtimes M \rightarrow M \rightarrow 1$ is the pushout of $j$ and $i$. Then for $[d] \in H^{1}(I,A,\chi_{I})^{\square}$, Lemma~\ref{lem:201} and Lemma~\ref{lem:202} guarantees that the inflation of $\partial([d])$ is trivial; that is, $\check{\sigma} \circ \partial ([d])=0$.

\underline{$\im \partial = \ker \check{\sigma}$} : Suppose $[T'] \in \ker \check{\sigma}$. Then there exists $h: M \rightarrow A$ such that 
	\begin{align*}
		T'(\pi(a,x),\pi(b,y)) &= h(a,x) + \sigma_{\left\langle a,x \right\rangle}(h(b,y)) - h( \left\langle a,x \right\rangle \cdot \left\langle b,y \right\rangle ) \\
		&= h(a,x) + \sigma_{\left\langle a,x \right\rangle}(h(b,y))-h(a \cdot \sigma^{s}_{x}(b) \cdot T^{s}(x,y), x\cdot y)
\end{align*}
	and
	\begin{align*}
		T'_{f}(\pi(a,x),\pi(b,y)) &= h(a,x) + h(b,y) + \tau_{\left\langle a, x \right\rangle}(h(b,y)) + \sigma_{\{\left\langle a, x \right\rangle, \left\langle b, y \right\rangle\}}(\nu_{\left\langle b, y \right\rangle}(h(a,x)) - h(a,x) - h(b,y)) \\
		&\quad - h \big(\big\{\left\langle a,x \right\rangle , \left\langle b,y \right\rangle\big\} \big)
	\end{align*}
	From this we see that
	\begin{align*}
		0 = T'(1,x) = T'(\pi(a,1),\pi(1,x)) &= h(a,1) + \sigma_{\left\langle a,1 \right\rangle}(h(1,x)) - h( a,x )
	\end{align*}
	which yields 
	\begin{align}\label{eqn:4000}
		h(a,x) = h(a,1) + \sigma_{\left\langle a,1 \right\rangle}(h(1,x)); 
	\end{align}
	similarly, we observe that
	\begin{align*}
		0 = T'(1,x \cdot y) = T'(\pi(T^{s}(x,y),1),\pi(1,x \cdot y)) &= h(T^{s}(x,y),1) + \sigma_{\left\langle T^{s}(x,y),1 \right\rangle}(h(1,x \cdot y)) - h( T^{s}(x,y),x \cdot y )
	\end{align*}
	which yields 
	\begin{align}\label{eqn:4001}
		h( T^{s}(x,y),x \cdot y ) = h(T^{s}(x,y),1) + \sigma_{\left\langle T^{s}(x,y),1 \right\rangle}(h(1,x \cdot y)) . 
	\end{align}
	We also observe that
	\begin{align*}
		0 = T'(x,y) = T'(\pi(1,x),\pi(1,y)) &= h(1,x) + \sigma_{\left\langle 1,x \right\rangle}(h(1,y)) - h( T^{s}(x,y),x \cdot y )
	\end{align*}
	and
	\begin{align}\label{eqn:kernelderiv}
		0 = T'(1,1) = T'(\pi(a,1),\pi(b,1)) &= h(a,1) + \sigma_{\left\langle a,1 \right\rangle}(h(b,1)) - h( a \cdot b,1 ) .
	\end{align}
	Define $\psi : I \rightarrow A$ by $\psi(a):= h(a,1)$. Then Eq~\eqref{eqn:kernelderiv} shows $\psi$ is a derivation for the group operation. Since $I \triangleleft M$, $\pi( \left\langle 1, x \right\rangle^{-1} \cdot \left\langle a, 1 \right\rangle \cdot \left\langle 1, x \right\rangle ) = 1$ and so
	\begin{align}\label{eqn:4003}
		\begin{split}
			0 = T'(x,1) &= T' \big( \pi(1,x),\pi( \left\langle 1, x \right\rangle^{-1} \cdot \left\langle a, 1 \right\rangle \cdot \left\langle 1, x \right\rangle) \big) \\
			&= h(1,x) + \sigma_{\left\langle 1,x \right\rangle} \big( h( \left\langle 1, x \right\rangle^{-1} \cdot \left\langle a , 1 \right\rangle \cdot \left\langle 1, x \right\rangle ) \big) - h( a,x ) .
		\end{split}
	\end{align}
	Then Eq~\eqref{eqn:4000} and Eq~\eqref{eqn:4003} gives
	\begin{align*}
		h(1,x) + \sigma_{\left\langle 1,x \right\rangle} \big( h( \left\langle 1, x \right\rangle^{-1} \cdot \left\langle a , 1 \right\rangle \cdot \left\langle 1, x \right\rangle ) \big) = h(a,x) = h(a,1) + \sigma_{\left\langle a, 1 \right\rangle}(h(1,x))
	\end{align*}
	which yields 
	\begin{align}
		\sigma_{\left\langle 1, x \right\rangle} ( \left\langle 1, x \right\rangle^{-1} \cdot \left\langle a, 1 \right\rangle \cdot \left\langle 1, x \right\rangle ) - h(a,1) = \sigma_{\left\langle a, 1 \right\rangle}(h(1,x)) - h(1,x) .
	\end{align}
	For the bracket operation, we see that
	\begin{align*}
		0 = T_{f}'(1,1) &= T_{f}'(\pi(a,1),\pi(b,1)) \\
		&= h(a,1) + h(b,1) + \tau_{\left\langle a,1 \right\rangle}(h(b,1)) + \sigma_{\{\left\langle a,1 \right\rangle, \left\langle b,1 \right\rangle\}} \big( \nu_{\left\langle b,1 \right\rangle}(h(a,1)) - h(a,1) - h(b,1) \big) 
	\end{align*}
	which shows $\psi \in \mathrm{Der}(I,A,\chi_{I})$. We also see that 
	\begin{align*}
		1 = T_{f}'(x,1) &= T_{f}'( \pi(1,x),\pi(a,1)) \\
		&= h(1,x) + h(a,1) + \tau_{\left\langle 1, x \right\rangle}(h(a,1)) + \sigma_{\{\left\langle 1, x \right\rangle, \left\langle a, 1 \right\rangle\}} \big(\nu_{\left\langle a, 1 \right\rangle}(h(1,x)) - h(1,x) - h(a,1) \big) \\
		&\quad - h \big(\{ \left\langle 1, x \right\rangle, \left\langle a, 1 \right\rangle\} \big)
	\end{align*}
	which gives
	\begin{align*}
		h(\tau^{s}_{x}(a),1) - \tau_{\left\langle 1, x \right\rangle}(h(a,1)) &+ \sigma_{\{ \left\langle 1, x \right\rangle, \left\langle a, 1 \right\rangle \}}(h(a,1)) - h(a,1) \\
		&= h(1,x) + \sigma_{\{ \left\langle 1, x \right\rangle, \left\langle a, 1 \right\rangle\}} \big( \nu_{\left\langle 1, x \right\rangle}(h(1,x)) - h(1,x) \big)
	\end{align*}
	Altogether, we have shown $[\psi] \in H^{1}(I,A,\chi_{I})^{\square}$ and we can take the choice $\eta(x):= h(1,x)$.

	From the condition $[T'] \in \ker \partial$, we have for the group operation that
	\begin{align*}
		0 = T'(1, x \cdot y) &= T'( \pi(T^{s}(x,y),1) , \pi(1, x \cdot y) ) \\
		&= h(T^{s}(x,y),1) + \sigma_{\left\langle T^{s}(x,y) , 1 \right\rangle}(h(1, x \cdot y)) - h(T^{s}(x,y), x \cdot y)
	\end{align*}
	and 
	\begin{align*}
		T'(x, y) &= T'( \pi(1,x) , \pi(1, y) ) \\
		&= h(1,x) + \sigma_{\left\langle 1 , x \right\rangle}(h(1,y)) - h(T^{s}(x,y), x \cdot y)
	\end{align*}
	which yields
	\begin{align*}
		\partial([\psi])(x,y) &= \Big\{h(1,x) + \sigma_{\left\langle 1, x \right\rangle}(h(1,y)) - \sigma_{\left\langle T^{s}(x,y), 1 \right\rangle}(h(1, x \cdot y)) - h(T^{s}(x,y),1) \Big\} \\
		&= \Big\{ h(1,x) + \sigma_{\left\langle 1 , x \right\rangle}(h(1,y)) - h( T^{s}(x,y) , x \cdot y ) \Big\} \\
		&= \Big\{ T'(x,y) \Big\}.
	\end{align*}
	For the bracket operation, we see that
	\begin{align*}
		T_{f}'(x,y) &= T_{f}'(\pi(1,x),\pi(1,y)) \\
		&= h(1,x) + h(1,y) + \tau_{\left\langle 1, x \right\rangle}(h(1,y)) + \sigma_{\{\left\langle 1, x \right\rangle, \left\langle 1, y \right\rangle\}}(\nu_{\left\langle 1, y \right\rangle}(h(1,x)) - h(1,x) - h(1,y)) \\
		&\quad - h \big(\big\{ \left\langle 1,x \right\rangle , \left\langle 1,y \right\rangle \big\} \big) \\ 
		&= h(1,x) + h(1,y) + \tau_{\left\langle 1, x \right\rangle}(h(1,y)) + \sigma_{\{\left\langle 1, x \right\rangle, \left\langle 1, y \right\rangle\}}(\nu_{\left\langle 1, y \right\rangle}(h(1,x)) - h(1,x) - h(1,y)) \\
		&\quad - h \big(T^{s}_{f}(x,y), \{x,y\} \big) 
	\end{align*}
	and 
	\begin{align*}
		0 = T'(1,\{x,y\}) &= T'(\pi(T^{s}_{f}(x,y), 1), \pi(1,\{x,y\})) \\
		&= h(T^{s}_{f}(x,y),1) + \sigma_{\left\langle T^{s}_{f}(x,y), 1 \right\rangle}(h(1,\{x,y\})) - h(\left\langle T^{s}_{f}(x,y), 1 \right\rangle \cdot \left\langle 1, \{x,y\} \right\rangle ) \\
		&= h(T^{s}_{f}(x,y),1) + \sigma_{\left\langle T^{s}_{f}(x,y), 1 \right\rangle }(h(1,\{x,y\})) - h(T^{s}_{f}(x,y), \{x,y\})
	\end{align*}
	which then yields
	\begin{align*}
		\partial([\psi])_{f}(x,y) &= \Big\{h(1,x) + h(1,y) + \tau_{\left\langle 1, x \right\rangle}(h(1,y)) + \sigma_{\{\left\langle 1, x \right\rangle, \left\langle 1, y \right\rangle\}} \big( \nu_{\left\langle 1, y \right\rangle}(h(1,x)) - h(1,x) - h(1,y) \big) \\
		&\quad - \sigma_{\left\langle T^{s}_{f}(x,y), 1 \right\rangle }(h(1,\{x,y\})) - h(T^{s}_{f}(x,y), 1 ) \Big\} \\
		&= \Big\{ h(1,x) + h(1,y) + \tau_{\left\langle 1, x \right\rangle}(h(1,y)) + \sigma_{\{\left\langle 1, x \right\rangle, \left\langle 1, y \right\rangle\}} \big(\nu_{\left\langle 1, y \right\rangle}(h(1,x)) - h(1,x) - h(1,y) \big) \\
		&\quad - h(T^{s}_{f}(x,y), \{x,y\}) \Big\} \\
		&= \Big\{ T_{f}'(x,y) \Big\} .
	\end{align*}
\end{proof}

\vspace{0.2cm}

\begin{acknowledgments}
Alexander Wires was supported by NSF China Grant \#12071374. 
\end{acknowledgments}


\end{document}